\documentclass[12pt]{article} 
\usepackage{a4,fullpage,amssymb,epsf,psfrag,times}

\makeatletter\@addtoreset{equation}{section}\makeatother
\makeatletter\@addtoreset{figure}{section}\makeatother
\makeatletter\@addtoreset{table}{section}\makeatother

\newtheorem{theorem}{Theorem}[section]
\newtheorem{proposition}[theorem]{Proposition}
\newtheorem{lemma}[theorem]{Lemma}

\newtheorem{corollary}[theorem]{Corollary}

\newcommand{\R}{{\mathbb R}}
\newcommand{\C}{{\mathbb C}}
\newcommand{\Z}{{\mathbb Z}}

\newcommand{\T}{{\mathbb T}}

\newcommand{\op}[1]{\!\!\mathop{\rm ~#1}\nolimits}
\newcommand{\fop}[1]{\!\!\mathop{\mbox{\rm \footnotesize ~#1}}\nolimits}
\newcommand{\scriptop}[1]{\!\!\mathop{\mbox{\rm \scriptsize ~#1}}\nolimits}
\newcommand{\tinyop}[1]{\!\!\mathop{\mbox{\rm \tiny ~#1}}\nolimits}

\newenvironment{proof}{\par\medskip\noindent{\bf Proof}~~}{
\unskip\nobreak\hfill\hbox{$\Box$}\par \bigskip}        

\newenvironment{remark}{\refstepcounter{theorem}\par\medskip\noindent{\bf Remark~\thetheorem~~}}{\unskip\nobreak\hfill\hbox{ $\oslash$}\par\bigskip}

\newenvironment{definition}{\refstepcounter{theorem}\par\medskip\noindent{\bf Definition~\thetheorem~~}}{\unskip\nobreak\hfill\hbox{ $\oslash$}\par\bigskip}

\newcommand{\got}[1]{\mathfrak{#1}}

\title{Symplectic torus actions with coisotropic principal orbits}
\author{J.J. Duistermaat\thanks{Research stimulated 
by a KNAW professorship} ~and A. Pelayo}

\begin{document}

\maketitle
\begin{abstract}
In this paper we completely classify symplectic actions of 
a torus $T$  on a compact connected 
symplectic manifold $(M,\,\sigma )$ 
when some, hence every, principal orbit is a coisotropic   
submanifold of $(M,\,\sigma)$. That is, we construct an explicit model, 
defined in terms of 
certain invariants, of the manifold, 
the torus action and the symplectic form. 
The invariants are invariants of the topology of the manifold, 
of the torus action, or of the symplectic form. 

In order to deal with symplectic actions which are
not Hamiltonian, we develop new techniques, extending the theory of 
Atiyah, Guillemin\--Sternberg, Delzant, and Benoist. More specifically, 
we prove that there is a well-defined notion of constant vector 
fields on the orbit space $M/T$. Using a generalization of 
the Tietze\--Nakajima theorem to what we call $V$\--parallel 
spaces, we obtain that $M/T$ is isomorphic to the Cartesian 
product of a Delzant polytope with a torus. 

We then construct special lifts of 
the constant vector fields on $M/T$, in terms of which  
the model of the symplectic manifold with the torus action 
is defined.   
\end{abstract}

\section{Introduction}
\label{intrsec}
Let $(M,\,\sigma )$ be a smooth compact and connected symplectic manifold 
of dimension $2n$ and let $T$ be a torus which acts 
effectively on $(M,\,\sigma )$ by means of symplectomorphisms. 
We furthermore assume that some principal $T$\--orbit is a 
coisotropic submanifold of $(M,\,\sigma )$, which implies that 
$d_T\geq n$ if $d_T$ denotes the dimension of $T$. See Lemma 
\ref{coisotropicorbitlem} for alternative characterizations 
of our assumptions. In this paper we will 
classify the compact connected symplectic manifolds 
with such torus actions, by constructing a list of 
explicit examples to which each of our manifolds 
is equivariantly symplectomorphic. See Theorem \ref{invariantthm}, 
Theorem \ref{existencethm} and 
Corollary \ref{crowncor} for our main result. 

In many integrable systems in classical mechanics, 
we have an effective Hamiltonian action of an $n$\--dimensional torus 
on the $2n$\--dimensional symplectic manifold, but also 
non\--Hamiltonian actions occur in physics, 
see for instance Novikov \cite{novikov}.  

If the effective action of $T$ on 
$(M,\,\sigma )$ is Hamiltonian, then $d_T=n$ 
and the principal orbits are Lagrange submanifolds.  
Moreover, the image of the momentum mapping is a convex polytope $\Delta$ 
in the dual space $\got{t}^*$ of $\got{t}$, where 
$\got{t}$ denotes the Lie algebra of $T$. $\Delta$ has the 
special property that at each vertex of 
$\Delta$ there are precisely $n$ codimension one faces 
with normals which form a $\Z$\--basis of the integral lattice 
$T_{\Z}$ in ${\got t}$, where $T_{\Z}$ is defined as 
the kernel of the exponential mapping from $\got{t}$ to $T$. 
The classification of Delzant \cite{delzant} 
says that for each such polytope $\Delta$ there is a 
compact connected symplectic manifold with Hamiltonian torus 
action having $\Delta$ as image of the momentum mapping, and the 
symplectic manifold with torus action is unique up to 
equivariant symplectomorphisms. Such polytopes $\Delta$ and corresponding 
symplectic $T$\--manifolds $(M,\,\sigma ,\, T)$ are  
called {\em Delzant polytopes} and {\em Delzant manifolds} 
in the exposition of this subject 
by Guillemin \cite{guillemin}, after 
Delzant \cite{delzant}.  Each Delzant manifold 
has a $T$\--invariant K\"ahler structure such that the 
K\"ahler form is equal to $\sigma$. 

Because critical points of the Hamiltonian function correspond 
to zeros of the Hamiltonian vector field, a Hamiltonian 
action on a compact manifold always has fixed points. 
Therefore the other extreme case of a symplectic torus 
action with coisotropic principal orbits occurs if the action 
is {\em free}. In this case, $M$ is a principal torus 
bundle over a torus, hence a nilmanifold for a two\--step 
nilpotent Lie group as described in  Palais and Stewart \cite{ps}. 
If the nilpotent 
Lie group is not commutative, then $M$ does not admit 
a K\"ahler structure, cf. Benson and Gordon \cite{bg}. For 
four\--dimensional manifolds $M$, 
these were the first examples of compact symplectic manifolds 
without K\"ahler structure, introduced by Thurston \cite{thurston}. 
See the end of Remark \ref{freerem}. 

The general case is a combination of the Hamiltonian case and the free case, 
in the sense that $M$ is an associated $G$\--bundle $G\times _HM_{\scriptop{h}}$ 
over $G/H$ with a $2d_{\scriptop{h}}$\--dimensional 
Delzant submanifold $(M_{\scriptop{h}},\, \sigma _{\scriptop{h}},\, T_{\scriptop{h}})$ of 
$(M,\,\sigma ,\, T)$  
as fiber. Here $T_{\scriptop{h}}$ is the unique maximal subtorus 
of $T$ which acts in Hamiltonian fashion on $(M,\,\sigma )$. 
It has dimension $d_{\scriptop{h}}$ and its Lie algebra is 
denoted by 
$\got{t}_{\scriptop{h}}$. $G$ is a two\--step nilpotent Lie group, and $H$ 
is a commutative closed Lie subgroup of $G$, which 
acts on $M_{\scriptop{h}}$ via $T_{\scriptop{h}}\subset H$. 
The base space 
$G/H$ is a torus bundle over a torus, see Remark \ref{Horem}. 
This leads to an 
explicit model of $(M,\,\sigma ,\, T)$ in terms 
of the ingredients 1) -- 6) in Definition \ref{ingredientdef}. 
See Proposition \ref{Aprop} and Proposition \ref{omegalem}. 
The model allows explicit computations 
of many aspects of $(M,\,\sigma ,\, T)$. 
As an example 
we determine the fundamental group of $M$ in  
Proposition \ref{pi1prop}, and the Chern classes 
of the normal bundle in $M/T_{\scriptop{f}}$ of the fixed point set 
of the action of $T_{\scriptop{h}}$ on $M/T_{\scriptop{f}}$ 
in Proposition \ref{chernprop}. 
Here $T_{\scriptop{f}}$ is a complementary subtorus to  
$T_{\scriptop{h}}$ in $T$, which acts freely on $M$.  
The main result of this paper is that the compact 
connected symplectic manifolds 
with symplectic torus action with coisotropic principal 
orbits are completely classified by the ingredients 
1) -- 6) in Definition \ref{ingredientdef}, see Theorem \ref{invariantthm} 
and Theorem \ref{existencethm}. 

The proof starts with the observation that the symplectic form 
on the orbits is given by a two\--form $\sigma ^{\got{t}}$ 
on $\got{t}$, see Lemma \ref{constlem}. 
Write $\got{l}:=\op{ker}\sigma ^{\got{t}}$. 
The inner product of the symplectic form $\sigma$ with the infinitesimal 
action of $T$ defines a closed basic $\got{l}^*$\--valued one\--form 
$\widehat{\sigma}$ on $M$, which turns the orbit space $M/T$ 
into a locally convex polyhedral $\got{l}^*$\--parallel space, as 
defined in Definition \ref{cornerdef}. 
The locally convex polyhedral $\got{l}^*$\--parallel space    
$M/T$ is isomorphic to $\Delta\times (N/P)$, 
in which $\Delta$ is a Delzant polytope in $(\got{t}_{\scriptop{h}})^*$ 
and $P$ is a cocompact discrete additive subgroup 
of the space $N$ of all linear forms 
on $\got{l}$ which vanish on $\got{t}_{\scriptop{h}}$. 
See Proposition \ref{M/Tprop}. 

The main step 
in the proof of the classification is the construction of lifts 
to $M$ of the constant vector fields on the $\got{l}^*$\--parallel
manifold $M/T$ with the simplest possible 
Lie brackets and symplectic products of the lifts. See 
Proposition \ref{liftprop}. This construction   
uses calculations involving the 
de Rham cohomology of $M/T$.  

All the proofs become much simpler in the case that 
the action of $T$ on $M$ is free. We actually first 
analyzed the free case with Lagrangian principal orbits, 
meaning that the principal orbits are Lagrange submanifolds 
of $M$. Next we treated    
the case with Lagrangian principal orbits  
where $M$ is fibered by Delzant manifolds, and only after 
we became aware of the article of Benoist 
\cite{benoist}, we generalized our results to the  
case with coisotropic principal orbits.  
In \cite{benoist}, Th. 6.6 states that every compact connected symplectic 
manifold with a symplectic torus action with coisotropic 
principal orbits is isomorphic
to the Cartesian product of a Delzant manifold and a 
compact connected symplectic manifold with a free  
symplectic torus action. 
However, even in the special case that the principal orbits 
are Langrange submanifolds of $M$, this conclusion 
appears to be too strong, 
if the word ``isomorphic'' implies 
``equivariantly diffeomorphic'', see 
Remark \ref{productrem} and Benoist 
\cite{benoistcorr}. 
 
The paper is organized as follows. In Section \ref{stasec} 
we discuss the condition that some (all) principal orbits are 
coisotropic submanifolds of $(M,\,\sigma )$. In Section 
\ref{M/Tsec} we analyze the space of $T$\--orbits in 
all detail, where we use the definitions and theorems 
in the appendix Section \ref{Vsec} concerning what we 
call ``$V$\--parallel spaces''. Section \ref{twolemmasec} 
contains a lemma about basic differential forms and one about equivariant 
diffeomorphisms which preserve the orbits. In Section 
\ref{liftsec} we construct our special lifts of constant 
vector fields on the orbit space. These are used in Section \ref{delzantsec} 
in order to construct the Delzant submanifolds of 
$(M,\,\sigma )$ and in Section \ref{normalformsec} for 
the normal form of the symplectic $T$\--manifold. 
The classification is completed by means of the theorems in 
Section \ref{invariantsec}. In the first appendix, Section \ref{Vsec}, 
we prove that every complete connected locally convex $V$\--parallel space 
is isomorphic to the Cartesian product of a closed convex subset 
of a finite\--dimensional vector space and a torus. 
See Theorem \ref{cornerclassthm} for the precise statement. 
This result is 
a generalization of the theorem of  
Tietze \cite{tietze} and Nakajima \cite{nakajima},  
which states that every closed and connected locally 
convex subset of a finite\--dimensional vector space is convex. 
In the second appendix, Section \ref{sympltubes}, we describe   
the local model of Benoist \cite[Prop. 1.9]{benoist} 
and Ortega and Ratiu \cite{ortegaratiu}
for a proper symplectic action of an arbitrary Lie group on an arbitrary 
symplectic manifold. 

There are many other texts on the classification of 
symplectic torus actions 
on compact manifolds which in some way are related to ours. 
The book of Audin \cite{audin} is on Hamiltonian torus 
actions, with emphasis on the topological aspects. 
Orlik and Raymond \cite{OR} and Pao \cite{pao} classified actions 
of two\--dimensional tori on four\--dimensional compact 
connected smooth manifolds. Because they do not assume 
an invariant symplectic structure, our classification 
in the four\--dimensional case forms only a tiny part of theirs. 
On the other hand the completely integrable 
systems with local torus actions of Kogan \cite{Ko} 
form a relatively close generalization of torus actions 
with Lagrangian principal orbits. 
The classification of Hamiltonian circle actions on 
compact connected four\--dimensional manifolds in Karshon 
\cite{K}, and of centered complexity one Hamiltonian 
torus actions in arbitrary dimensions in Karshon and 
Tolman \cite{KT}, are also much richer than our 
classification in the case that $n-d_{\scriptop{h}}\leq 1$. 
McDuff \cite{M} and McDuff and Salamon \cite{MS} studied 
non\--Hamiltonian circle actions, and Ginzburg 
\cite{ginzburg} non\--Hamiltonian symplectic actions of 
compact groups under the assumption of a ``Lefschetz 
condition''. In another direction 
Symington \cite{s} and Leung and Symington \cite{ls} 
classified four\--dimensional compact connected symplectic 
manifolds which are fibered by Lagrangian tori where 
however the fibration is allowed to have elliptic or 
focus\--focus singularities. 

We are very grateful to Yael Karshon for her suggestion 
of the problem. A. Pelayo thanks her for moral and intellectual support
during this project.   

\section{Coisotropic principal orbits}
\label{stasec}
Let $(M,\,\sigma )$ be a smooth compact and connected symplectic manifold 
and let $T$ be a torus which acts effectively on 
$(M,\,\sigma )$ by means of symplectomorphisms. 
In this section we show that some principal $T$\--orbit is 
a coisotropic submanifold of $(M,\,\sigma )$ if and 
only if the Poisson brackets of any pair of smooth 
$T$\--invariant functions on $M$ vanish if and only if 
every principal $T$\--orbit is a coisotropic submanifold of 
$(M,\,\sigma )$. See Lemma \ref{coisotropicorbitlem}, 
Remark \ref{multfreerem} and Remark \ref{orbtyperem} below. 

This follows from the local model of Benoist \cite[Prop. 1.9]{benoist}, 
see Theorem \ref{Gthm}, which in the case of symplectic 
torus actions with coisotropic principal orbits assumes a 
particularly simple form, see Lemma \ref{modellem}. 

\medskip
If $X$ is an element of the Lie algebra $\got{t}$ 
of $T$, then we denote by $X_M$ the infinitesimal 
action of $X$ on $M$. It is a smooth vector field 
on $M$, and the invariance of $\sigma$ under the action of 
$T$ implies that 
\begin{equation}
\op{d}(\op{i}_{X_M}\sigma )=\op{L}_{X_M}\sigma =0.
\label{LXsigma}
\end{equation} 
Here $\op{L}_v$ denotes the Lie derivative with respect to the 
vector field $v$, and $\op{i}_v\omega$ the inner product 
of a differential form $\omega$ with $v$, obtained by inserting 
$v$ in the first slot of $\omega$. 
The first identity in (\ref{LXsigma}) 
follows from the homotopy 
identity $\op{L}_v =\op{d}\circ\op{i}_v +\op{i}_v\circ\op{d}$ 
combined with $\op{d}\!\sigma =0$.  

If $f$ is a smooth real\--valued function on $M$, then the 
unique vector field $v$ on $M$ such that $-\op{i}_v\sigma =
\op{d}\! f$ is called the {\em Hamiltonian vector field}  
of $f$, and will be denoted by $\op{Ham}_f$. Given 
$v$, the function $f$ is uniquely determined up to 
and additive constant, which implies 
that $f$ is $T$\--invariant if and only if $v$ is $T$\--invariant. 
If $X\in\got{t}$, then $X_M$ is Hamiltonian if and only if 
the closed two\--form $\op{i}_{X_M}\sigma$ is exact. 

The following lemma says that the pull\--back to the $T$\--orbits of 
the symplectic form $\sigma$ on $M$ is given by a constant 
antisymmetric bilinear form on the Lie algebra $\got{t}$ of $T$. 
\begin{lemma}
There is a unique antisymmetric bilinear form 
$\sigma^{\got{t}}$ on $\got{t}$, such that 
such that 
\[
\sigma _x(X_M(x),\, Y_M(x))=\sigma^{\got{t}}(X, Y)
\] 
for every $X,\, Y\in\got{t}$ and every $x\in M$. 
\label{constlem}
\end{lemma}
\begin{proof}
It follows from Benoist \cite[Lemme 2.1]{benoist} 
that if $u$ and $v$ are smooth 
vector fields on $M$ 
such that $\op{L}_u\sigma =0$ and $\op{L}_v\sigma =0$, 
then $[u,\, v]=\op{Ham}_{\sigma (u,\, v)}$. We repeat the proof. 
\[
\op{i}_{[u,\, v]}\sigma =\op{L}_u(\op{i}_v\sigma )
=\op{i}_u(\op{d}(\op{i}_v\,\sigma ))+\op{d}(\op{i}_u(\op{i}_v\sigma))
=\, -\op{d}(\sigma (u,\, v)).
\]
Here we used $\op{L}_u\sigma =0$ in the first equality, 
the homotopy formula for the Lie derivative in the second identity, 
and finally $\op{d}\!\sigma =0$, the homotopy identity and 
$\op{L}_v\sigma =0$ in the third equality.  
Applying this 
to $u=X_M$, $v=Y_M$ for $X,\, Y\in\got{t}$, 
and using that $[X,\, Y]=0$, hence 
$[X_M,\, Y_M]=\, -[X,\, Y]_M=0$, it follows that 
$\op{Ham}_{\sigma (X_M,\, Y_M)}=0$. Thus   
$\op{d}(\sigma (X_M,\, Y_M))=0$, and  
the function $x\mapsto\sigma _x(X_M(x),\, Y_M(x))$ is constant 
on $M$, because $M$ is connected.  
\end{proof}

In the further discussion we will need some basic facts 
about proper actions of Lie groups, see for instance 
\cite[Sec. 2.6--2.8]{dk}. 
For each $x\in M$ we write $T_x:=\{ t\in T\mid t\cdot x=x\}$ 
for the {\em stabilizer subgroup} of the $T$\--action at the point $x$. 
$T_x$ is a closed Lie subgroup of $T$, it has finitely many 
components and its identity component is a torus subgroup of $T$. 
The Lie algebra $\got{t}_x$ of $T_x$ is equal to the space 
of all $X\in\got{t}$ such that $X_M(x)=0$. In other words, 
$\got{t}_x$ is the kernel of the linear mapping 
$\alpha _x:X\mapsto X_M(x)$ from $\got{t}$ to $\op{T}_x\! M$. 
The image of $\alpha _x$ is equal to the tangent space 
at $x$ of the $T$\--orbit through $x$, and 
will be denoted by $\got{t}_M(x)$. The linear mapping 
$\alpha _x:\got{t}\to\op{T}_x\! M$ induces a linear isomorphism from 
$\got{t}/\got{t}_x$ onto $\got{t}_M(x)$. 

For each closed subgroup $H$ of $T$ which 
can occur as a stabilizer subgroup, 
the {\em orbit type} $M^H$ is defined as the set of all $x\in M$ 
such that $T_x$ is conjugate to $H$, but because $T$ is commutative 
this condition is equivalent to the equation $T_x=H$. 
Each connected component $C$ of $M^H$ is a smooth 
$T$\--invariant submanifold of $M$.  
The connected components of the orbit types in $M$ 
form a finite partition of $M$, which actually is a Whitney stratification.  
This is called the {\em orbit type stratification} of $M$. 
There is a unique open orbit type, called the {\em principal orbit type}, 
which is the orbit type of a subgroup $H$ which is contained 
in every stabilizer subgroup $T_x$, $x\in M$. Because 
the effectiveness of the action means that the intersection 
of all the $T_x$, $x\in M$ is equal to the identity element, 
this means that the principal orbit type consists of the 
points $x$ where $T_{x}=\{ 1\}$, that is where the action is free. 
If the action is free at $x$, then  
the linear mapping $X\mapsto X_M(x)$ from $\got{t}$ to $\op{T}_x\! M$ 
is injective. The points $x\in M$ at which the $T$\--action 
is free are also called the {\em regular points} of $M$, 
and the principal orbit type, the set of all regular points 
in $M$ is denoted by $M_{\scriptop{reg}}$. 
The principal orbit type $M_{\scriptop{reg}}$ is 
a dense open subset of 
$M$, and connected because $T$ is connected, 
see \cite[Th. 2.8.5]{dk}. The {\em principal orbits} 
are the orbits in $M_{\scriptop{reg}}$, the principal orbit type. 
In our situation, the principal orbits are the orbits on which 
the action of $T$ is free. 
\begin{lemma}
Let $\got{l}$ be the kernel in $\got{t}$ of the 
two\--form $\sigma ^{\got{t}}$ on $\got{t}$ defined in 
Lemma \ref{constlem}, the set of all $X\in\got{t}$ 
such that $\sigma ^{\got{t}}(X,\, Y)=0$ for every 
$Y\in\got{t}$. Then $\got{t}_x\subset\got{l}$ for 
every $x\in M$. 
\label{txinslem}
\end{lemma}
\begin{proof}
If $X\in\got{t}_x$, then $X_M(x)=0$, hence 
$\sigma ^{\got{t}}(X,\, Y)=\sigma _x(X_M(x),\, Y_M(x))=0$ 
for every $Y\in\got{t}$. 
\end{proof}
The linear subspace $\got{l}$ of $\got{t}$ 
will play an important 
part in the classification of the symplectic torus actions 
with coisotropic principal orbits. 

A submanifold $C$ of $M$ is called {\em coisotropic}, 
if for every $x\in C$, $v\in\op{T}_x\! M$, the condition that 
$\sigma _x(u,\, v)=0$ for every $u\in\op{T}_x\! C$ implies that 
$v\in\op{T}_x\! C$. In other words, if the $\sigma_x$\--orthogonal 
complement $(\op{T}_x\! C)^{\sigma _x}$ of $\op{T}_x\! C$ in $\op{T}_x\! M$ 
is contained in $\op{T}_x\! C$.  
Every symplectic manifold has an even dimension, 
say $2n$, and if $C$ is a coisotropic submanifold of dimension $k$, 
then 
\[
2n-k=\op{dim}(\op{T}_x\! C)^{\sigma _x}\leq
\op{dim}(\op{T}_x\! C)=k
\]
shows that $k\geq n$. $C$ has the minimal dimension $n$ 
if and only if $(\op{T}_x\! C)^{\sigma _x}=\op{T}_x\! C$, 
if and only if $C$ is Lagrange submanifold of $M$, 
an isotropic submanifold of $M$ of maximal dimension $n$. 
The next lemma is basically the implication 
(iv) $\Rightarrow$ (ii) in Benoist \cite[Prop. 5.1]{benoist}. 
\begin{lemma} 
Let $(M,\,\sigma )$ be a connected symplectic manifold, 
and $T$ a torus which acts effectively and 
symplectically on $(M,\,\sigma )$. 
Then every coisotropic $T$\--orbit is a principal orbit. 
Furthermore, if some $T$\--orbit is coisotropic, then 
every principal orbit is coisotropic, and   
$\op{dim}M=\op{\dim}T+\op{dim}\got{l}$. 
\label{coisotropicorbitlem}
\end{lemma}
\begin{proof}
We use Theorem \ref{Gthm} with $G=T$, where we note that 
the commutativity of $T$ implies that the adjoint action of 
$H=T_x$ 
on $\got{t}$ is trivial, which implies that the coadjoint 
action of $H$ on the component $\got{l}/\got{h}$ is trivial as well. 

Let us assume that the orbit $T\cdot x$ is coisotropic, 
which means that  
$\got{t}_M(x)^{\sigma _x}\subset\got{t}_M(x)$, 
or equivalently the subspace $W$ defined in 
(\ref{Wdef}) is equal to 
zero. This implies that the action of $H$ on $E=(\got{l}/\got{h})^*$ 
is trivial, and the vector bundle 
$T\times_HE=T\times_H(\got{l}/\got{h})^*$ is 
$T$\--equivariantly isomorphic to 
$(T/H)\times (\got{l}/\got{h})^*$, 
where $T$ acts by left multiplications on
the first factor. It follows that in the model 
all stabilizer subgroups are equal to $H$, and therefore 
$T_y=H$ for all $y$ in the 
$T$\--invariant open neighborhood $U$ of $x$ in $M$. 
Because the principal orbit type is dense in $M$, 
there are $y\in U$ such that $T_y=\{ 1\}$, and it follows 
that $T_x=H=\{ 1\}$, that is, $T\cdot x$ is a principal orbit. 
We note in passing that this 
implies that $\op{dim}M=\op{\dim}T+\op{dim}\got{l}$. 

When $W=\{ 0\}$, we read off from (\ref{Gsigma}) with 
$\sigma ^{G/H}$ given by $\sigma ^{\got{t}}$ in Lemma \ref{constlem},  
and (\ref{etadef}), that the symplectic form $\Phi ^*\sigma$ is given by 
\[
(\Phi ^*\sigma )_{(t\, H,\,\lambda )}
((X+\got{h},\,\delta\lambda ),\, 
(X'+\got{h},\,\delta '\lambda ))
=\sigma ^{\got{t}}(X,\, X')
+\delta\lambda (X'_{\got{l}})-\delta '\lambda (X_{\got{l}})
\]
for all $(t\, H,\,\lambda )\in (T/H)\times E_0$, 
and $(X+\got{h},\,\delta\lambda ),\, (X'+\got{h},\,\delta '\lambda )
\in (\got{t}/\got{h})\times (\got{l}/\got{h})^*$. 
In this model, the tangent space of the $T$\--orbit 
is the set of all $(X'+\got{h},\,\delta '\lambda )$ 
such that $\delta '\lambda =0$, of which the 
symplectic orthogonal complement is equal to the 
set of all $(X+\got{h},\,\delta \lambda )$ such that 
$X\in\got{l}$ and $\delta\lambda =0$, which implies 
that in this model every $T$\--orbit is coisotropic 
and therefore  the orbit $T\cdot y$ is coisotropic 
for every $y\in U$. This shows that the 
set of all $x\in M$ such that $T\cdot x$ is coisotropic 
is an open subset of $M$. Because for all $x\in M_{\scriptop{reg}}$ 
the tangent spaces of the orbits $T\cdot x$ 
have the same dimension, equal to 
$\op{dim}T$, the set of all $x\in M_{\scriptop{reg}}$ such 
that $T\cdot x$ is coisotropic is closed in $M_{\scriptop{reg}}$. 
Because $M_{\scriptop{reg}}$ is connected, it follows that 
$T\cdot x$ is coisotropic for all $x\in M_{\scriptop{reg}}$ as soon 
as $T\cdot x$ is coisotropic for some $x\in M_{\scriptop{reg}}$.   
\end{proof}

\begin{remark}
In the proof of Lemma \ref{coisotropicorbitlem}, linear forms 
on $\got{l}/\got{h}$ were identified with linear forms on 
$\got{l}$. For any linear subspace $F$ of a 
finite\--dimensional vector space 
$E$ we have the canonical projection $p:x\mapsto x+F:E\mapsto E/F$, 
and its dual mapping $p^*:(E/F)^*\to E^*$. Because 
$p$ is surjective, $p^*$ is injective, and 
its image $p^*((E/F)^*)$ is equal to the space 
$F^0$ of all $\varphi\in E^*$ 
such that $\varphi |_F=0$. 
This leads to a canonical 
linear isomorphism $p^*$ from $(E/F)^*$ onto $F^0$, which 
will be used throughout this paper 
to identify $(E/F)^*$ with the linear 
subspace $F^0$ of $E^*$. 
\label{E/Frem}
\end{remark}

\begin{remark}
Let $x\in M_{\scriptop{reg}}$. 
Because the principal orbit type $M_{\scriptop{reg}}$ is fibered by 
the $T$\--orbits, the tangent space $\got{t}_M(x)$ at $x$ 
of $T\cdot x$ is equal to the common kernel of the 
$\op{d}\! f(x)$, where $f$ ranges over the $T$\--invariant 
smooth functions on $M$. Because $-\op{d}\! f=\op{i}_{\fop{Ham}_f}\sigma$, 
it follows that $\got{t}_M(x)^{\sigma _x}$ is equal to the 
set of all $\op{Ham}_f(x)$, $f\in\op{C}^{\infty}(M)^T$. 
Here $\op{C}^{\infty}(M)^T$ denotes the space of all 
$T$\--invariant smooth functions on $M$. 

Suppose that the principal orbits are coisotropic and 
let $f\in\op{C}^{\infty}(M)^T$. Then we have for every 
$x\in M_{\scriptop{reg}}$ 
that $\op{Ham}_f(x)\in\got{t}_M(x)^{\sigma _x}\cap\got{t}_M(x)$, 
or $\op{Ham}_f(x)=X(x)_M(x)$ for a uniquely determined $X(x)\in\got{l}$. 
It follows that the $\op{Ham}_f$\--flow leaves every principal 
orbit invariant, and because $M_{\scriptop{reg}}$ is dense in 
$M$, the $\op{Ham}_f$\--flow leaves every $T$\--orbit invariant. 
Because a point $x\in M$ is called a {\em relative equilibrium} 
of a $T$\--invariant vector field $v$ if the $v$\--flow leaves 
$T\cdot x$ invariant, the conclusion is that all points of $M$ 
are relative equilibria of $\op{Ham}_f$, 
and the induced flow in $M/T$ is at standstill. 
Moreover the $T$\--invariance of $\op{Ham}_f$ implies 
that $x\mapsto X(x)\in\got{l}$ is constant on each 
principal $T$\--orbit, 
which implies that the $\op{Ham}_f$\--flow in $M_{\scriptop{reg}}$ 
is {\em quasiperiodic}, in  
the direction of the infinitesimal action of $\got{l}$ on 
$M_{\scriptop{reg}}$. 

If $f,\, g\in\op{C}^{\infty}(M)^T$ and $x\in M_{\scriptop{reg}}$, 
then $\op{Ham}_f(x)$ and $\op{Ham}_g(x)$ both belong to 
$\got{t}_M(x)^{\sigma _x}\cap\got{t}_M(x)$, and it follows that 
the {\em Poisson brackets} $\{ f,\, g\} :=\op{Ham}_f\, g
=\sigma (\op{Ham}_f,\,\op{Ham}_g)$ of $f$ and $g$ 
vanish at $x$. Because $M_{\scriptop{reg}}$ is dense 
in $M$, it follows that $\{ f,\, g\}\equiv 0$ for all 
$f,\, g\in\op{C}^{\infty}(M)^T$ if the principal orbits 
are coisotropic. 

If conversely $\{ f,\, g\} \equiv 0$ for all 
$f,\, g\in\op{C}^{\infty}(M)^T$, then we have for every 
$x\in M_{\scriptop{reg}}$ that $\got{t}_M(x)^{\sigma _x}
\subset (\got{t}_M(x)^{\sigma _x})^{\sigma _x}=\got{t}_M(x)$, 
which means that $T\cdot x$ is coisotropic. 
Therefore the principal orbits are coisotropic if and 
only if the Poisson brackets of all $T$\--invariant smooth 
functions vanish. 

In Guillemin and Sternberg \cite{multfree}, a symplectic manifold with a 
Hamiltonian action of an arbitrary compact Lie group is called 
a {\em multiplicity\--free space} if the Poisson brackets of any pair of 
invariant smooth functions vanish. Because in \cite{multfree} 
the emphasis is on representations of noncommutative 
compact Lie groups, which do not play a role in our paper, 
and because on the other hand we allow non\--Hamiltonian 
actions, we did not put the adjective ``multiplicity\--free'' 
in the title. 
\label{multfreerem}
\end{remark}
The next lemma is statement (1) (a) in Benoist 
\cite[Lemma 6.7]{benoist}. For general symplectic torus 
actions the stabilizer subgroups need not be connected. 
For instance, there exist symplectic torus actions 
with symplectic orbits and nontrivial finite stabilizer 
subgroups. 
\begin{lemma} 
Let $(M,\,\sigma )$ be a connected symplectic manifold, 
and $T$ a torus which acts effectively and 
symplectically on $(M,\,\sigma )$, with 
coisotropic principal orbits. Then, for every $x\in M$, 
the stabilizer group $T_x$ is connected, that is, a subtorus of $T$. 
\label{Txconnectedlem}
\end{lemma}
\begin{proof}
As in the proof of Lemma \ref{coisotropicorbitlem}, 
we use Theorem \ref{Gthm} with $G=T$, where $H$ acts 
trivially on the factor $(\got{l}/\got{h})^*$ in 
$E=(\got{l}/\got{h})^*\times W$. 
Recall that $t\in T$ acts on $T\times _HE$ by 
sending $H\cdot (t',\, e)$ to $H\cdot (t\, t',\, e)$. 
When $t=h\in H$, then 
\[
H\cdot (h\, t',\, e)=H\cdot (h\, t'\, h^{-1},\, h\cdot e)=
H\cdot (t',\, h\cdot e)
\]
because $T$ is commutative, and we see that the action 
of $H$ on $T\times _HE$ is represented by the 
linear symplectic action 
of $H$ on $W$, where $W$ is defined by (\ref{Wdef}).  

Because 
\[
\op{dim}M=(\op{dim}T+\op{dim}(\got{l}/\got{h})+\op{dim}W)-\op{dim}H 
\]
and because the assumption that the principal orbits are 
coisotropic implies that $\op{dim}M=\op{dim}T+\op{dim}\got{l}$, 
see Lemma \ref{coisotropicorbitlem}, 
it follows that $\op{dim}W=2\op{dim}H$. 

Write $m=\op{dim}H$.  
The action of the compact and commutative group $H$ 
by means of symplectic linear transformations on the 
$2m$\--dimensional symplectic vector space 
$(W,\,\sigma ^W)$ leads to a direct sum decomposition 
of $W$ into $m$ mutually $\sigma ^W$\--orthogonal 
two\--dimensional $H$\--invariant linear subspaces 
$E_j$, $1\leq j\leq m$. 

For $h\in H$ and every $1\leq j\leq m$, let $\iota _j(h)$ denote 
the restriction to $E_j\subset W\simeq \{ 0\}\times W
\subset (\got{l}/\got{h})^*\times W$ of the action of $h$ on $E$. 
Note that $\op{det}\iota _j(h)=1$, because 
$\iota _j(h)$ preserves the restriction to $E_j\times E_j$ of 
$\sigma ^W$, which is an area form on $E_j$.  
Averaging any inner product 
in each $E_j$ over $H$, we obtain an $H$\--invariant inner 
product $\beta _j$ on $E_j$, and $\iota _j$ is a homomorphism of Lie 
groups from $H$ to $\op{SO}(E_j,\,\beta _j)$, the group of 
linear transformations of 
$E_j$ which preserve both $\beta _j$ 
and the orientation. 

On the other hand, if $h\in H$ and $w\in W_{\scriptop{reg}}$, then 
\[
h\cdot w=\sum _{j=1}^m\,\iota _j(h)\, w_j
\quad\mbox{\rm if}\quad w=\sum_{j=1}^m\, w_j,\quad w_j\in E_j. 
\]
Therefore $\iota _j(h)\, w_j=w_j$ for all $1\leq j\leq m$ 
implies that $h\cdot w=w$, hence $h=1$. 
This implies that the 
homomorphism of Lie groups $\iota$, defined by
\[
\iota :h\mapsto (\iota _1(h),\,\ldots ,\,\iota _m(h)): 
H\to \prod_{j=1}^m\,\op{SO}(E_j,\,\beta _j),
\]
is injective. Because both the source group $H$ and the target group 
are $m$\--dimensional Lie groups, and the target group 
is connected, it follows that $\iota$ is 
an isomorphism of Lie groups. This implies in turn that 
$H$ is connected. 
\end{proof}

\begin{remark} 
The $H$\--invariant inner product $\beta _j$ on $E_j$, 
introduced in the proof of Lemma \ref{Txconnectedlem},  
is unique, if we also require that the symplectic 
inner product of any orthonormal basis with respect to 
$\sigma ^W$ is equal to $\pm 1$. In turn this leads to 
the existence of a unique complex structure on 
$E_j$ such that, for any unit vector 
$e_j$ in $(E_j,\,\beta _j)$, we have that 
$e_j$, $\op{i}\, e_j$ is an orthonormal basis 
in $(E_j,\,\beta _j )$ and $\sigma ^W(e_j,\,\op{i}\, e_j)=1$. 
Here $\op{i}:=\sqrt{-1}\in\C$. 
This leads to an identification of $E_j$ with $\C$, 
which is unique up to multiplication by an 
element of $\T :=\{ z\in\C\mid |z|=1\}$. 
 
In turn this leads to an identification of $W$ with $\C^ m$, 
with the symplectic form $\sigma ^W$ defined by 
\begin{equation}
\sigma ^{\C ^m}
=\sum_{j=1}^m\,\op{d}\overline{z^j}\wedge
\op{d}\! z^j/2\op{i}. 
\label{sigmaW}
\end{equation}
The element $c\in\T ^m$ acts 
on $\C ^m$ by sending $z\in\C ^m$ to the element 
$c\cdot z$ such that $(c\cdot z)^j=c^j\, z^j$ for every 
$1\leq j\leq m$. There is a unique isomorphism 
of Lie groups $\iota :H\to\T ^m$ such that 
$h\in H$ acts on $W=\C ^m$ by sending 
$z\in\C ^m$ to $\iota (h)\cdot z$. 

The identification of $W$ with $\C ^m$ is unique up to 
a permutation of the coordinates and the action of an 
element of $\T ^m$. 
\label{complexrem}
\end{remark}

In the local model of Lemma \ref{modellem} below, we will use 
that any subtorus of a torus has a complementary 
subtorus, in the following sense.  
\begin{lemma}
Let $U$ be a $d_U$\--dimensional subtorus of 
a $d_T$\--dimensional torus $T$. Let $U_{\Z}$ and 
$T_{\Z}$ denote the integral lattice, the kernel of the 
exponential mapping, in the Lie algebra $\got{u}$ and 
$\got{t}$ of $U$ and $T$, respectively. 
Let $Y_i$, $1\leq i\leq d_U$, be a $\Z$\--basis 
of $U_{\Z}$. Then there are 
$Z_j$, $1\leq j\leq d_V:=d_T-d_U$, such that 
the $Y_i$ and $Z_j$ together form a $\Z$\--basis of $T_{\Z}$. 
If we denote by $\got{v}$ the span of the $Z_j$,  
then $V=\exp\got{v}$ is a subtorus of $T$ 
with Lie algebra equal to $\got{v}$. $V$ is a 
complementary subtorus of $U$ in $T$ in the sense that 
the mapping 
$
U\times V\ni (u,\, v)\mapsto 
u\, v\in T
$
is an isomorphism from 
$U\times V$ onto $T$. 
The $Z_j$ form a $\Z$\--basis of the integral lattice 
$V_{\Z}$ in the Lie algebra $\got{v}$ of $V$. 
\label{tbaselem}
\end{lemma}
\begin{proof}
We repeat the well\--known argument. 
If $X\in T_{\Z}$, $c\in\Z$, $c\neq 0$, and 
$c\, X\in U_{\Z}$, then 
$X\in\got{u}$ and 
$\op{exp}X=1$ in $T$, hence $\op{exp}X=1$ in $U$, 
and it follows that $X\in U_{\Z}$. This means that 
the finitely generated commutative group $T_{\Z}/U_{\Z}$ 
is torsion\--free, and therefore has a $\Z$\--basis 
$\widetilde{Z}_j$, $1\leq j\leq k$, cf. 
Hungerford \cite[Th. 6.6 on p. 221]{hungerford}. 
We have that $\widetilde{Z}_j=Z_j+
U_{\Z}$ for some $Z_j\in T_{\Z}$. If $X\in T_{\Z}$, 
then there are unique $z^j\in\Z$ such that 
$X+U_{\Z}=\sum_{j=1}^k\, z^j\,\widetilde{Z}_j$, 
which means that $X-\sum_{j=1}^k\, z^j\, Z_j\in 
U_{\Z}$. But this implies that there are unique 
$y^i\in\Z$ such that $X-\sum_{j=1}^k\, z^j\, Z_j
=\sum_{i=1}^{d_U}\, y^i\, Y_i$, which shows that 
the $Y_i$ and $Z_j$ together form a $\Z$\---basis of $T_{\Z}$, 
which in turn implies that $k=d_T-d_U=d_V$. 

The last statement follows from the fact that 
the mapping 
\[
(y,\, z)\mapsto\exp\left(\sum_{i=1}^{d_U}
\, y^i\, Y_i+\sum_{j=1}^{d_V}
\, z^j\, Z_j\right) 
\]
from $\R ^{d_T}$ to $T$ induces an isomorphism from 
$(\R /\Z)^{d_T}$ onto $T$ which maps $(\R /\Z)^{d_U}
\times\{ 0\}$ onto $U$ and 
$\{ 0\}\times (\R/\Z)^{d_V}$ onto $V$. 
\end{proof}

\begin{remark}
The complementary subtorus $V$ in Lemma \ref{tbaselem} 
is by no means unique. The $Z_j$ can be replaced 
by any 
\[
Z_j'=Z_j+\sum_{i=1}^{d_U}\, c^i_j\, Y_i,\quad 1\leq j\leq d_V,
\]
in which the $c^i_j$ are integers. This leads to 
a bijective correspondence between the set of 
all complementary subtori of a given 
subtorus $U$ and the set of all  
$d_U\times d_F$\--matrices with integral 
coefficients. 
\label{tbaserem}
\end{remark}
Let $H=T_x$ be the subtorus of $T$ in Lemma \ref{Txconnectedlem}. 
Let $K$ be a complementary subtorus of $H$ in $T$ 
and, for any $t\in T$, let $t_H$ and $t_K$ 
be the unique elements in $H$ and $K$, respectively, such that 
$t=t_H\, t_K$. Let $X\mapsto X_{\got{l}}$ be a linear 
projection from $\got{t}$ onto $\got{l}$. 
We also use the identification of $W$ with $\C ^m$ as in 
Remark \ref{complexrem}. 
With these notations, we have the following local model 
for our symplectic $T$\--space with coisotropic 
principal orbits.
\begin{lemma}
Under the assumptions of Lemma \ref{Txconnectedlem}, 
there is an isomorphism of Lie groups $\iota$ 
from $H$ onto $\T ^m$, an open $\T ^m$\--invariant 
neighborhood $E_0$ of the origin in 
$E=(\got{l}/\got{h})^*\times \C ^m$,  
and a $T$\--equivariant diffeomorphism $\Phi$ from 
$K\times E_0$ onto 
an open $T$\--invariant neighborhood $U$ of $x$ in $M$, 
such that $\Phi (1,\, 0)=x$. 
Here $t\in T$ acts on $K\times (\got{l}/\got{h})^*\times \C ^m$ 
by sending $(k,\,\lambda ,\, z)$ to 
$(t_K\, k,\,\lambda ,\, \iota (t_H)\cdot z)$. 
In addition, the symplectic form $\Phi ^*\sigma$ 
on $K\times E_0$ 
is given by 
\begin{equation}
(\Phi ^*\sigma )_{(k,\,\lambda ,\, z)}
((X,\,\delta\lambda ,\,\delta z),\, (X',\,\delta '\lambda ,\,\delta 'z))
=\sigma ^{\got{t}}(X,\, X')
+\delta\lambda (X'_{\got{l}})-\delta '\lambda (X_{\got{l}})
+\sigma ^{\C ^m}(\delta z,\,\delta 'z)
\label{modelsigma}
\end{equation}
for all $(k,\,\lambda ,\, z)\in K\times (\got{l}/\got{h})^*\times\C ^m$, 
and $(X,\,\delta\lambda ,\,\delta z),\, (X',\,\delta '\lambda ,\,\delta 'z)
\in\got{k}\times (\got{l}/\got{h})^*\times\C ^m$. Here we 
identify each tangent space of the torus $K$ 
with $\got{k}$ and each tangent space of a vector space with the 
vector space itself. Finally, $\sigma ^{\C ^m}$ is the symplectic form on 
$\C ^m$ defined in (\ref{sigmaW}).
\label{modellem}
\end{lemma}
\begin{proof}
As in the proof of Lemma \ref{Txconnectedlem},  
we use Theorem \ref{Gthm} with $G=T$, where $H$ acts 
trivially on the factor $(\got{l}/\got{h})^*$ and 
$h\in H$ acts on $W=\C ^m$ by sending 
$z\in\C ^m$ to $\iota (h)\cdot z$. Here $\iota :H\to\T ^m$ 
is the isomorphism from the torus $H$ onto the standard torus 
$\T ^m$ introduced in Remark \ref{complexrem}, and  
the symplectic form $\sigma ^{\C ^m}$ on $\C ^m$ is 
given by (\ref{sigmaW}). 

Because $K$ is a complementary subtorus of $H$ in $T$, 
the manifold $K\times E$ is a global section of the 
vector bundle $\pi _K:T\times _HE\to T/H\simeq K$. Indeed, 
if $(t,\, e)\in T\times E$, then $(t_K,\, t_H\cdot e)
=(t\, {t_H}^{-1},\, t_H\cdot e)$ is the unique 
element in $(K\times E)\cap H\cdot (t,\, e)$. 
Furthermore, if $t\in T$ and $(k,\, e)\in K\times E$, 
then $(t_K\, k,\, t_H\cdot e)$ is the unique element 
in $(K\times E)\cap H\cdot (t\, k,\, e)$. 
This exhibits $T\times _HE$ as a trivial vector bundle over $K$, 
which is a homogeneous $T$\--bundle, where $t\in T$ acts 
on $K\times E$ by sending $(k,\, e)$ to $(t_K\, k,\, t_H\cdot e)$. 

Finally, if in (\ref{etadef}) we restrict ourselves to 
$X\in\got{k}$, then the right hand side 
simplifies to $\lambda (X_{\got{l}})+\sigma ^W(w,\,\delta w)/2$, 
which leads to (\ref{modelsigma}). 
\end{proof}

\begin{remark}
In the local model of Lemma \ref{modellem}, we have that 
$T_{(k,\,\lambda ,\, z)}=H$ if and only if $z$ is a fixed 
point of $\iota (H)=\T ^m$ if and only if $z=0$. 
Because $K\times (\got{l}/\got{h})^*\times\{ 0\}$ is a 
symplectic submanifold of $K\times (\got{l}/\got{h})^*\times\C ^m$, 
it follows that every orbit type is a smooth symplectic submanifold 
of $(M,\,\sigma )$. 

Moreover, $T\cdot (k,\,\lambda ,\, 0)=K\times\{\lambda\}\times\{ 0\}$ 
is a coisotropic submanifold of $K\times (\got{l}/\got{h})^*\times\C ^m$, 
and we conclude that every $T$\--orbit is a coisotropic submanifold 
of its orbit type. 

The discussion of the relative equilibria in 
Remark \ref{multfreerem}, with $M_{\scriptop{reg}}$ replaced 
by any orbit type $M^H$, 
leads to the conclusion that for every $f\in\op{C}^{\infty}(M)^T$ the flow of 
the Hamiltonian vector field $\op{Ham}_f$ 
in $M^H$ is quasiperiodic, in the direction of the 
infinitesimal action of $\got{l}/\got{h}$ in $M^H$.  
\label{orbtyperem}
\end{remark}
We conclude this section with a discussion of the 
special case that the two\--form $\sigma^{\got t}$ in Lemma 
\ref{constlem} is equal to zero.
\begin{lemma}
We have 
$\sigma ^{\got{t}}=0$  
if and only if $\got{l}:=\op{ker}\sigma ^{\got{t}}=\got{t}$  
if and only if some $T$\--orbit is isotropic 
if and only if  
every $T$\--orbit is isotropic.

Also, every principal orbit is a Lagrange submanifold 
of $(M,\,\sigma )$  
if and only if some principal orbit 
is a Lagrange submanifold of $(M,\,\sigma )$  
if and only if $\op{dim}M=2\op{dim}T$ and $\sigma^{\got t}=0$. 
\label{isotropicorbitlem}
\end{lemma}
\begin{proof}
The equivalence of $\sigma ^{\got{t}}=0$ and 
$\op{ker}\sigma ^{\got{t}}=\got{t}$  is obvious, whereas the 
equivalence between $\sigma ^{\got{t}}=0$ and the isotropy 
of some (every) $T$\--orbit follows from Lemma 
\ref{constlem}. 

If $x\in M_{\scriptop{reg}}$ and $T\cdot x$ is a Lagrange 
submanifold of $(M,\,\sigma )$, then 
$\op{dim}M=2\op{dim}(T\cdot x)=2\op{dim}T$, and 
$\sigma^{\got t}=0$ follows in view of the first statement 
in the lemma. 

Conversely, if $\op{dim}M=2\op{dim}T$ and $\sigma^{\got t}=0$, 
then every orbit is isotropic and for every 
$x\in M_{\scriptop{reg}}$ we have 
$\op{dim}M=2\op{dim}T=2\op{dim}(T\cdot x)$, 
which implies that $T\cdot x$ is a Lagrange submanifold 
of $(M,\,\sigma )$. 
\end{proof}

\section{The orbit space}
\label{M/Tsec}
In this section we investigate the orbit space of our 
action of the torus $T$ on the compact connected symplectic manifold 
$(M,\,\sigma )$ with coisotropic principal orbits. 
The main results are that 
the closed basic one\--form $\widehat{\sigma}$ of 
Lemma \ref{alphalem} exhibits the orbit space as a 
locally convex polyhedral $\got{l}^*$\--parallel space, 
see Definition 
\ref{cornerdef} and Lemma \ref{parallellem}, and that 
as such $M/T$ is isomorphic to the Cartesian product of 
a Delzant polytope and a torus, see Proposition \ref{M/Tprop}. 
The assumption that the principal orbits are coisotropic 
will be assumed throughout this section, unless explicitly 
stated otherwise. 

\subsection{Canonical local charts on the orbit space}
\label{canchartss}
In this subsection we exhibit the space of $T$\--orbits as an 
$\got{l}^*$\--parallel space in the sense of Definition 
\ref{cornerdef}.

\medskip
We denote the space of all orbits in $M$ of the $T$\--action 
by $M/T$, and by $\pi :M\to M/T$ the 
canonical projection which assigns to each $x\in M$ the 
orbit $T\cdot x$ through the point $x$. The orbit space 
is provided with the maximal topology for which the 
canonical projection is continuous; this topology is 
Hausdorff. 

For each connected component $C$ of an orbit type $M^H$ in $M$ 
of the subgroup $H$ of $T$, as introduced in the 
paragraphs preceding Lemma \ref{txinslem}, the action of $T$ on $C$ induces a 
proper and free action of the torus $T/H$ on $C$, and $\pi(C)$ 
has a unique structure of a smooth manifold such that 
$\pi :C\to\pi(C)$ is a principal $T/H$\--bundle. 
$(M/T)^H:=\pi (M^H)$ is called the {\em orbit type of $H$ in $M/T$} 
and $\pi(C)$ is a connected component of $(M/T)^H$. 
The connected components of the orbit types in the orbit space 
form a finite stratification of the orbit space, cf. 
\cite[Sec. 2.7]{dk}. 

Although $M/T$ is equal to the union of the finitely many 
strata of the orbit type stratification in $M/T$, where 
each of these strata is a smooth manifold, the orbit 
space $M/T$ is not a smooth manifold, unless the action 
of $T$ on $M$ is free. In general the principal orbit type  
$(M/T)_{\scriptop{reg}}=M_{\scriptop{reg}}/T$ is a smooth 
manifold of dimension $\op{dim}M-\op{dim}T$, which is 
an open and dense subset of $M/T$,and $M/T$ will 
have singularities at the lower dimensional strata, 
the strata in the complement of $(M/T)_{\scriptop{reg}}$ 
in $M/T$. However, in this section we will obtain a much more 
explicit description of the orbit space $M/T$. 

A smooth differential form $\omega$ on $M$ is called 
{\em basic} with respect to the $T$\--action if 
it is $T$\--invariant, that is $\op{L}_{X_M}\omega = 0$ 
for every $X\in\got{t}$, and if $\op{i}_{X_M}\omega = 0$ 
for every $X\in\got{t}$. The basic differential forms constitute 
a module over the algebra $\op{C}^{\infty}(M)^T$ of 
$T$\--invariant smooth functions on $M$, the basic forms of 
degree zero on $M$. A smooth 
differential form $\omega$ on $M$ is basic if and only if 
the restriction of $\omega$ to the principal orbit type is equal to 
$\pi^*\nu$ for a smooth differential form $\nu$ 
on the principal orbit type in $M/T$. 

A theorem of Koszul \cite{koszul} says that the \v{C}ech 
(= sheaf) 
cohomology group $\op{H}^k(M/T,\,\R )$ of $M/T$ is canonically 
isomorphic to the de Rham cohomology of the basic forms 
on $M$, that is, the space of closed basic $k$\--forms on $M$ modulo 
its subspace consisting of the $\op{d}\!\nu$ in which 
$\nu$ ranges over the basic $(k-1)$\--forms on $M$. 
This theorem holds for any proper action of a Lie group on any 
smooth manifold, and in particularly it does not need the 
compactness of $M$. 

\begin{lemma}
Recall that $\got{l}$ is the kernel of the antisymmetric 
bilinear form $\sigma ^{\got{t}}$ which had been introduced 
in Lemma \ref{constlem}. For each $X\in\got{l}$,
$\widehat{\sigma}(X) := \, -\op{i}_{X_M}\sigma$ is a closed 
basic one\--form on $M$. 
\label{alphalem}
\end{lemma}
\begin{proof}
That $\widehat{\sigma}(X)$ is closed follows from (\ref{LXsigma}).  
Because $X\in\got{l}$, we have for each $Y\in {\got t}$ that  
$-\op{i}_{Y_M}(\widehat{\sigma} (X)) =\sigma ^{\got{t}}(X,\, Y) = 0$. 
Also we have for every $Y\in\got{t}$ that 
\[
-\op{L}_{Y_M}\widehat{\sigma} (X) =\op{i}_{[Y_M,\, X_M]}\sigma + 
\op{i}_{X_M}(\op{L}_{Y_M}\sigma) = 0.
\]
Here we have used the Leibniz identity for the Lie derivative, 
the commutativity of ${\got t}$, and the $T$\--invariance of $\sigma$ 
which implies that $\op{L}_{Y_M}\sigma =0$. 
\end{proof}

For each $x\in M$, $\widehat{\sigma} (X)_x$ is a linear form 
on $\op{T}_x\! M$ which depends linearly on $X\in\got{l}$, and 
therefore $X\mapsto\widehat{\sigma} (X)_x$ is an $\got{l}^*$\--valued 
linear form on $\op{T}_x\! M$, which we denote by 
$\widehat{\sigma}_x$. In this way $x\mapsto\widehat{\sigma}_x$ 
is an $\got{l}^*$\--valued one\--form on $M$, which we 
denote by $\widehat{\sigma}$. With these conventions, we have 
\begin{equation}
\widehat\sigma _x(v)(X)=\widehat{\sigma} (X)_x(v)
=\sigma _x(v,\, X_M(x)),
\quad x\in M,\; v\in\op{T}_x\! M,\; X\in\got{l}.
\label{widehatsigma}
\end{equation}
Note that the $\got{l}^*$\--valued one\--form 
$\widehat{\sigma}$ on $M$ is basic and closed. 

Let $X\in\got{t}$ and suppose that 
$X_M=\op{Ham}_f$ for some $f\in\op{C}^{\infty}(M)$.  
Then we have for every 
$Y\in\got{t}$ that 
\[
Y_M(f)=\op{i}_{Y_M}(\op{d}\! f)=-\op{i}_{Y_M}(\op{i}_{X_M}\sigma )
=\sigma (Y_M,\, X_M)=\sigma ^{\got{t}}(Y,\, X),
\]
and it follows that $f\in\op{C}^{\infty}(M)^T$ if and only if 
$X\in\got{l}:=\op{ker}\sigma ^{\got{t}}$. 
The $T$\--action on $(M,\,\sigma )$ is called a 
{\em Hamiltonian $T$\--action} if for every $X\in\got{t}$ 
there exists an $f\in\op{C}^{\infty}(M)^T$ 
such that $X_M=\op{Ham}_f$. Note that if $\got{l}=\got{t}$, 
that is, if $\sigma ^{\got{t}}=0$, then the $T$\--action is 
Hamiltonian if and only if if for every $X\in\got{t}$ 
there exists an $f\in\op{C}^{\infty}{M}$ such that $X_M=\op{Ham}_f$. 

We recall the  Delzant manifolds, mentioned in Section \ref{intrsec}. 
Koszul's theorem now implies the following. 
\begin{corollary}
We do not assume that the principal orbits 
are coisotropic. Let $X\in\got{t}$. Then $X_M=\op{Ham}_f$ for some 
$f\in\op{C}^{\infty}(M)^T$, 
if and only if $X\in\got{l}:=\op{ker}{\sigma}^{\got{t}}$ 
and the cohomology class 
$[\widehat{\sigma} (X)]\in
\op{H}^1(M/T,\,\R )$ is equal to zero. 
If the $T$\--action is Hamiltonian, then 
$\sigma ^{\got{t}}=0$. Finally, if  
$\sigma ^{\got{t}}=0$ and $\op{H}^1(M/T,\,\R )=0$, 
then the $T$\--action is Hamiltonian  
and $(M,\,\sigma ,\, T)$ is a Delzant manifold.   
\label{Hamcor}
\end{corollary}

\begin{remark} 
In the local model of Lemma \ref{modellem}, the $T$\--orbit space 
of $K\times E_0$ is equal to $E_0/\T ^m$, which is contractible 
by using the radial contractions in $E_0$. It follows 
that for every $x_0\in M$ there is a $T$\--invariant open neighborhood 
$U$ of $x_0$ in $M$ such that the open subset $\pi(U)$ of 
the orbit space $M/T$ is contractible. Because of Koszul's theorem, 
and because the \v{C}ech cohomology of $\pi(U)$ is trivial, 
it follows that the infinitesimal action of $\got{l}$ on $U$ is Hamiltonian. 
Therefore, if $\sigma ^{\got{t}}=0$, then the 
$T$\--action is {\em locally Hamiltonian} in the sense that 
every element in $M$ has a $T$\--invariant open neighborhood 
in $M$ on which the $T$\--action is Hamiltonian. 
\label{lochamrem}
\end{remark}

In the local model of Lemma \ref{modellem}, we write 
$z^j=|z^j|\,\op{e}^{\fop{i}\theta ^j}$ with 
$\theta ^j\in\R/2\pi\Z$ for each $1\leq j\leq m$. 
Then the symplectic form $\sigma ^W$ with $W=\C ^m$ in 
(\ref{sigmaW}) is equal to 
\begin{equation}
\sigma ^{\C ^m}=\sum_{j=1}^m\,\op{d}\!\rho _j\wedge\op{d}\theta ^j,
\quad\mbox{\rm in which}\quad \rho_j:=|z^j|^2/2. 
\label{rhodef}
\end{equation}
The mapping $(\lambda ,\,\rho): 
\overline{M}
:=K\times (\got{l}/\got{h})^*\times\C ^m\to (\got{l}/\got{h})^*\times\R ^m$ 
induces a homeomorphism from the $T$\--orbit space 
$\overline{M}/T\simeq 
(\got{l}/\got{h})^*\times (\C ^m/\T ^m)$ onto 
$(\got{l}/\got{h})^*\times\R _+^m$, in which 
\[
\R _+^m:=\{\rho\in\R ^m\mid \rho_j\geq 0\;
\mbox{\rm for every}\; 1\leq j\leq m\} .
\]

Note that $(\op{e}^{\fop{i}\alpha ^1},\,\ldots ,
\,\op{e}^{\fop{i}\alpha ^m})\in\T ^m$ acts on $\C ^m$ 
by sending $\theta$ to $\theta +\alpha$ and leaving 
$\rho$ fixed. If we identify the Lie algebra 
of $\T ^m$ with $(\op{i}\R )^m$, then the infinitesimal action 
of $\beta\in (\op{i}\R )^m$ in $(\theta ,\,\rho)$\--coordinates 
is equal to the constant vector field $(\beta ,\, 0)$. 
The tangent mapping at $1$ of the isomorphism $\iota :H\to\T ^m$ 
is a linear isomorphism from $\got{h}$ onto 
$(\op{i}\R )^m$, which we we also denote by $\iota$. 

For every $Y\in\got{t}$, the infinitesimal action 
$Y_{\overline{M}}$ of $Y$ 
on $\overline{M}$ is equal to 
the vector field $(Y_{\got{k}},\, 0,\, \iota (Y_{\got{h}})\cdot z)$, 
see the description of the action of $T$ on the 
model in Lemma \ref{modellem}. 
Write $Y=Y_{\got{h}}+Y_{\got{k}}$ with 
$Y_{\got{h}}\in\got{h}$ and $Y_{\got{k}}\in\got{k}$. Here 
$\got{k}$ denotes the Lie algebra of the complementary 
torus $K$ to $H$ in $T$, which implies that $\got{t}
=\got{h}\oplus\got{k}$.  
Because $\got{h}\subset\got{l}$, we have 
$(Y_{\got{h}})_{\got{l}}=Y_{\got{h}}$ and  
$Y_{\got{l}}=Y_{\got{h}}+(Y_{\got{k}})_{\got{l}}
=Y_{\got{h}}+(Y_{\got{l}})_{\got{k}}$. 
Because $\delta\lambda\in (\got{l}/\got{h})^*$ is a linear 
form on $\got{l}$ which is equal to zero on $\got{h}$, 
it follows that 
\begin{equation}
\delta\lambda ((Y_{\got{k}})_{\got{l}})=\delta\lambda 
(Y_{\got{l}}),\quad \delta\lambda\in (\got{l}/\got{h})^*, 
\quad Y\in\got{t}.
\label{gotkgotl}
\end{equation} 
Therefore, 
if in (\ref{modelsigma}) we substitute 
\[
(X',\,\delta '\lambda ,\,\delta 'z)=Y_{\overline{M}}(k,\,\lambda ,\, z)
=(Y_{\got{k}},\, 0,\,\iota (Y_{\got{h}})\cdot z) 
\]
with $Y\in\got{l}$, then we obtain 
\begin{equation}
\sigma ^{\got{t}}(X,\, Y_{\got{k}}) 
+\delta\lambda ((Y_{\got{k}})_{\got{l}}) 
+\sigma ^{\C ^m}(\delta z,\,\iota (Y_{\got{h}})\cdot z)
=\delta\lambda (Y)+\sum_{j=1}^m\,
\,\iota (Y_{\got{h}})^j\,\delta\rho _j/\op{i}.
\label{iYsigma}
\end{equation}
Here we have used that $Y_{\got{k}}=
(Y_{\got{l}})_{\got{k}}=(Y_{\got{k}})_{\got{l}}\in\got{l}
:=\op{ker}\sigma
^{\got{t}}$ implies that $\sigma ^{\got{t}}(X,\, Y_{\got{k}})=0$. 
Furthermore (\ref{gotkgotl}) with $Y\in\got{l}$ implies that  
$\delta\lambda ((Y_{\got{k}})_{\got{l}})=\delta\lambda (Y)$.  
Finally the formula for the $\sigma ^{\C ^m}$\--term follows 
from (\ref{rhodef}), as $(\op{d}\!\rho _j)(\iota (Y_{\got{h}})\cdot z)=0$, 
$(\op{d}\!\rho _j)(\delta z)=\delta\rho _j$, 
and $(\op{d}\!\theta ^j)(\iota (Y_{\got{h}})\cdot z)
=\iota (Y_{\got{h}})^j/\op{i}$ because 
the infinitesimal action of $\beta :=\iota (Y_{\got{h}})\in 
\op{i}\R ^m$ is equal to 
$(\sum_j\,\beta ^j\partial /\partial\theta ^j)/\op{i}$. 

Consider the linear mapping
\begin{equation}
A:(\delta\lambda ,\,\delta\rho )\mapsto [Y\mapsto\delta\lambda (Y)
+\sum_{j=1}^m\,\iota (Y_{\got{h}})^j\,\delta\rho _j/\op{i}]
\label{Aiota}
\end{equation}
from $(\got{l}/\got{h})^*\times\R ^m$ onto $\got{l}^*$. 
$A$ is a linear isomorphism, because the source space and the 
target space have the same dimension, and $\op{ker}A=0$: 
testing with arbitrary $Y\in\got{h}$ yields that 
$\delta\rho =0$, and then testing with 
arbitrary $Y\in\got{l}$ yields that $\delta\lambda =0$.  

Let $X_j$ denote the element of $\got{h}\subset\got{l}$ 
such that $\iota (X_j)=2\pi\op{i}\, e_j$, in which $e_j$ denotes 
the $j$\--th standard basis vector in $\R ^m$. \
Note that the $2\pi\op{i}\, e_j$, $1\leq j\leq m$, form a 
$\Z$\--basis of the integral lattice of the Lie algebra 
of $\T ^m$, and because $\iota :H\to\T ^m$ is an isomorphism 
of tori, it follows that the $X_j$, $1\leq j\leq m$, form 
a $\Z$\--basis of the integral lattice of the Lie algebra 
$\got{h}$ of $H$. Also note that 
\[
(A(\lambda ,\,\rho ))(X_j)=\sum_{k=1}^m\,\iota (X_j)^k\,\rho _k/\op{i}
=2\pi\rho _j,\quad 1\leq j\leq m. 
\]
Here it is essential that we use the coordinates $\rho _j$ 
instead of their infinitesimal displacements 
$\delta\rho _j$, because in Lemma \ref{chartlem} 
below we are interested in the the consequences of the 
inequalities $\rho _j\geq 0$. 
This leads to the following conclusion. 
\begin{lemma}
Let $\Phi$ be the $T$\--equivariant 
symplectomorphism from $K\times E_0\subset\overline{M}$ onto the 
open $T$\--invariant neighborhood $U$ of $x$ in $M$ as 
introduced in Lemma \ref{modellem}.  
Then the smooth mapping $\Psi :U\to\got{l}^*$, which consists 
of $\Phi ^{-1}:U\to\overline{M}$, followed by 
the $(\lambda ,\,\rho )$\--map and then $A$, induces 
a homeomorphism $\chi$ from $U/T$ onto an open neighborhood 
of $0$ in the corner 
\[
\{\xi\in\got{l}^*\mid \xi (X_j)\geq 0\;
\mbox{\rm for every}\; 1\leq j\leq m\} 
\]
in $\got{l}^*$, such that $\widehat{\sigma}=\op{d}\!\Psi$. 
Here the $X_j$, $1\leq j\leq m$, form a $\Z$\--basis 
of the integral lattice of the Lie algebra $\got{h}\subset 
\got{l}$ of $H$. 
\label{chartlem}
\end{lemma}
\begin{proof}
For every $Y\in\got{l}$, the right hand side of (\ref{iYsigma})  
is equal to $-\op{i}_{Y_M}\sigma$. Combined with 
the definitions of $\Psi$ and $A$,  this yields 
that $\widehat{\sigma}=\op{d}\!\Psi$. 
Because $(k,\, \lambda ,\, z)\mapsto 
(\lambda ,\rho )$ is a homeomorphism 
from $(K\times E)/T$ onto $(\got{l}/\got{h})^*\times 
(\R _+)^m$, $\chi$ is a homeomorphism from $U/T$ onto an open neighborhood 
of 0 in the corner in $\got{l}^*$ which is determined 
by the inequalities $\xi (X_j)\geq 0$, $1\leq j\leq m$. 
\end{proof}

If $\widetilde{\Psi}:\widetilde{U}\to\got{l}^*$ is mapping as in Lemma 
\ref{chartlem}, with corresponding chart $\widetilde{\chi}:
\widetilde{U}/T\to\got{l}^*$, then 
\[
\op{d}(\Psi -\widetilde{\Psi})=\op{d}\!\Psi -\op{d}\!\widetilde{\Psi} 
=\widehat{\sigma}-\widehat{\sigma}=0
\]
shows that $\Psi -\widetilde{\Psi}$ is locally 
constant on $U\cap\widetilde{U}$, which implies that 
$\chi -\widetilde{\chi}$ is locally constant on 
$(U/T)\cap (\widetilde{U}/T)$. In terms of Definition \ref{cornerdef}, 
we have proved 
\begin{lemma}
With the $\chi$ of Lemma \ref{chartlem} 
as local charts on $M/T$, the orbit space $M/T$ is 
a locally convex polyhedral $\got{l}^*$\--parallel space. 
The linear forms $v^*_{\alpha ,\, j}$, $1\leq j\leq m$, 
in Definition \ref{cornerdef} 
are the $\got{l}^*\ni\xi\mapsto \xi (X_j)$, where 
the $X_j$, $1\leq j\leq m$, form a $\Z$\--basis of the 
integral lattice of the Lie algebra $\got{h}\subset\got{l}$ 
of a stabilizer group $H=T_x$ of an element $x\in M$. 
\label{parallellem} 
\end{lemma}

In the next lemma we will introduce the subtorus $T_{\scriptop{h}}$ 
of $T$ which later will turn out to be the unique maximal 
subtorus of $T$ which acts on $M$ in a Hamiltonian fashion. 
For this reason $T_{\scriptop{h}}$ will be called the 
{\em Hamiltonian torus}.

\begin{lemma}
There are only finitely many different stabilizer subgroups of $T$, 
each of which 
is a subtorus of $T$. 
The product $T_{\scriptop{h}}$ of all the different stabilizer subgroups 
is a subtorus of $T$, and the Lie algebra $\got{t}_{\scriptop{h}}$ of 
$T_{\scriptop{h}}$ is 
equal to the sum of the Lie algebras of all the different stabilizer 
subgroups of $T$. 
It follows from Lemma \ref{txinslem} that $\got{t}_{\scriptop{h}}\subset 
\got{l}:=\op{ker}\sigma ^{\got{t}}$.  
\label{stabilizerlem} 
\end{lemma}
\begin{proof}
In the local model of Lemma \ref{modellem}, 
the stabilizer subgroup of $(k,\,\lambda ,\, z)$ is equal 
to the set of all $h\in H$ such that $\iota (h)^j=1$ for every 
$j$ such that $z^j\neq 0$. It follows that we have 
$2^m$ different stabilizer subgroups $T_y$, $y\in U$, 
namely one for each subset of $\{ 1,\,\ldots ,\, m\}$. 
Because $M$ is compact, is follows that there 
are only finitely many different stabilizer subgroups of $T$. 
For the last statement we observe that the product of 
finitely many subtori is a compact and connected subgroup of 
$T$ and therefore a subtorus of $T$. Also the image under the 
exponential mapping of the sum of the finitely many 
different Lie algebras of the stabilizer 
subgroups of $T$ is equal to $T_{\scriptop{h}}$, which proves that 
the Lie algebra of $T_{\scriptop{h}}$ is equal to the sum of the 
finitely many different $\got{t}_x$, $x\in M$. 
See for instance \cite[Sec. 1.12]{dk} for the general facts 
about Lie subgroups of tori, which we have used here. 
\end{proof}

\begin{remark}
The orbit $\pi (x)=T\cdot x\in M/T$ of any $x\in M_{\scriptop{reg}}$ 
is called a {\em regular point of $M/T$}. 
Recall from the paragraph preceding Lemma \ref{txinslem} 
that $x\in M_{\scriptop{reg}}$ if and only if 
$T_x=\{ 1\}$ if and only if ${\got t}_x=\{ 0\}$. 
Therefore the set $(M/T)_{\scriptop{reg}}$ of all regular points in 
$M/T$ is just the principal orbit type, which is a smooth 
manifold of dimension $\op{dim}M-\op{dim}T$. 
In the local model of Lemma \ref{modellem} 
with $x\in M_{\scriptop{reg}}$, 
where $\got{h}=\got{t}_x=\{ 0\}$ and $m=0$, 
at each point the $\got{l}^*$\--valued 
one\--form $\widehat{\sigma}$ corresponds to the projection 
$(\delta t,\,\delta\lambda )\mapsto\delta\lambda :
\got{t}\times\got{l}^*\to\got{l}^*$, and 
$\got{t}\times\{ 0\}$ is equal to the tangent space 
of the $T$\--orbit. It follows that   
for every $p\in (M/T)_{\scriptop{reg}}$ the induced 
linear mapping 
$\widehat{\sigma}_p: \op{T}_p(M/T)_{\scriptop{reg}}\to\got{l}^*$ 
is a linear isomorphism. 

More generally, the orbit type stratification, introduced 
in the paragraph preceding Lemma \ref{txinslem}, leads to a 
corresponding decomposition of $M/T$. 
The strata for the $T$\--action in 
$\overline{M}=K\times (\got{l}/\got{h})^*\times\C ^m$  
are of the form $\overline{M}^J$ in which $J$ is a subset of 
$\{ 1,\,\ldots,\, m\}$ 
and $\overline{M}^J$ is the set of all 
$(k,\,\lambda ,\, z)$ such that $z^j=0$ if and only if $j\in J$. 
In terms of the $(\theta ,\,\rho )$\--coordinates, 
this corresponds to $\rho _j = 0$ for all $j\in J$ and 
$\rho _k>0$ for $k\notin J$. 
The Lie algebra of the corresponding stabilizer subgroup of $T$ 
corresponds to the span of the vector fields 
$\partial/\partial\theta ^j$ with 
$j\in J$. Therefore, if $\Sigma$ is a connected component of 
the orbit type in $M/T$ defined by the subtorus $H$ of $T$ 
with Lie algebra 
$\got{h}$, then for each $p\in\Sigma$ we have 
$\widehat{\sigma}_p(X) = 0$ for 
all $X\in {\got h}$, and $\widehat{\sigma}_p$ may be viewed as an 
element of $(\got{l}/\got{h})^*={\got h}^0$, 
the set of all linear forms on $\got{l}$ 
which vanish on $\got{h}$, see Remark \ref{E/Frem}. 
The linear mapping $\widehat{\sigma}_p : 
\op{T}_p\!\Sigma\to (\got{l}/\got{h})^*$ is a linear isomorphism. 
\label{regrem}
\end{remark}

\subsection{$M/T$ is the Cartesian product of a Delzant polytope 
and a torus}
\label{prodss}
In the following Proposition \ref{M/Tprop}, 
the orbit space $M/T$ is viewed as a locally convex polyhedral 
$\got{l}^*$\--parallel 
space, as in Definition \ref{cornerdef} 
with $Q=M/T$ and $V=\got{l}^*$. See Lemma \ref{parallellem}. 
Let the subset $D$ of $\got{l}^*\times (M/T)$ and the mapping 
$(\xi ,\, p)\mapsto p+\xi$ from $D$ to $M/T$ 
be defined as in Definition \ref{Ddef}. 
We have the linear subspace $N$ of $V=\got{l}^*$, 
which acts on $Q=M/T$ by means of 
translations, and the period group $P$ of the $N$\--action 
on $Q$, as defined in Lemma \ref{Nlem} and Lemma \ref{perlem}, 
respectively. With the choice of a base point $p\in M/T$, we write 
$D_p=\{\xi\in\got{l}^*\mid (\xi ,\, p)\in D\}$. 
Let $\got{t}_{\scriptop{h}}'$ be a linear complement of $\got{t}_{\scriptop{h}}$ in 
$\got{t}$ and let $p\in M/T$. With these definitions, 
and the identification of $(\got{l}/\got{t}_{\scriptop{h}})^*$ 
with the space of linear forms on $\got{l}$ which 
vanish on $\got{t}_{\scriptop{h}}$, see Remark \ref{E/Frem}, we have
the following conclusions. 

\begin{proposition}
Let $C$ be a linear complement of 
$(\got{l}/\got{t}_{\scriptop{h}})^*$ in $\got{l}^*$. 
\begin{itemize}
\item[i)] $N=(\got{l}/\got{t}_{\scriptop{h}})^*$, $P$ is a cocompact discrete 
subgroup of the additive group $N$, and 
$N/P$ is a $\op{dim}N$\--dimensional torus. 
\item[ii)] 
There is a Delzant polytope $\Delta$ in $C\simeq (\got{t}_{\scriptop{h}})^*$, 
such that $D_p=\Delta +N$. 
\item[iii)] The  
mapping $\Phi_p :(\eta ,\,\zeta )\mapsto p+(\eta +\zeta )$ 
is an isomorphism of locally convex polyhedral $\got{l}^*$\--parallel spaces 
from $\Delta\times (N/P)$ onto $M/T$. 
\end{itemize}
\label{M/Tprop}
\end{proposition}
\begin{proof}
The linear forms $v^*_j$ which appear in the 
characterization of $N$ in Theorem \ref{cornerclassthm} 
are equal to the collection of all the 
$X_i\in\got{h}\subset\got{l}=(\got{l}^*)^*$ 
which appear in $\Z$\--bases of integral lattices 
of Lie algebras $\got{h}$ of stabilizer subgroups $H$ of $T$. 
Because $N$ is equal to the 
common kernel of all the $v^*_j$, $N$ is equal to the set 
$(\got{l}/\got{t}_{\scriptop{h}})^*$ of all elements 
of $\got{l}^*$ which vanish on the sum $\got{t}_{\scriptop{h}}\subset\got{l}$ 
of the finitely many different Lie algebras $\got{h}$ 
of stabilizer subgroups of $T$. 

Because $C$ is a linear complement of $(\got{l}/\got{t}_{\scriptop{h}})^*$ 
in $\got{l}^*$, the mapping $\xi\mapsto\xi |_{\got{t}_{\scriptop{h}}}$ 
induces an isomorphism from $C$ onto $(\got{t}_{\scriptop{h}})^*$. 
$\Delta$ is a Delzant polytope in $(\got{t}_{\scriptop{h}})^*$ in the sense 
of Guillemin \cite[p. 8]{guillemin}, because 
each $\Z$\--basis of the integral lattice of $\got{t}_x$ can 
be extended to a $\Z$\--basis of the integral lattice of 
$\got{t}_{\scriptop{h}}$, see Lemma \ref{tbaselem}. 

Because $C$ is a linear complement of 
$(\got{l}/\got{t}_{\scriptop{h}})^*=N=\R\, P$ in 
$\got{l}^*$, Proposition \ref{M/Tprop} 
now follows from Lemma \ref{parallellem} and Theorem \ref{cornerclassthm}. 
\end{proof}

\begin{corollary}
Let $(M,\,\sigma )$ be a compact connected 
$2n$\--dimensional symplectic manifold and suppose that 
we have an effective symplectic action of an $n$\--dimensional torus $T$ 
on $(M,\,\sigma )$, where we do not assume that the 
principal orbits are coisotropic. 
Then the following conditions are equivalent. 
\begin{itemize}
\item[i)] The action of $T$ has a fixed point in $M$. 
\item[ii)] The sum of the Lie algebras of all the 
different stabilizer subgroups of $T$ is equal to the Lie algebra of $T$. 
\item[iii)] $\sigma ^{\got{t}}=0$ and 
$M/T$ is homeomorphic to a convex polytope.
\item[iv)] $\sigma ^{\got{t}}=0$ and $\op{H}^1(M/T,\,\R )=0$. 
\item[v)] The action of $T$ is Hamiltonian. 
\end{itemize}
\label{delzantprop}
\end{corollary}
\begin{proof}
If $x$ is a fixed point, then $T_x=T$, hence $\got{t}_x=\got{t}$, 
which implies ii). 

Write $\got{t}'$ for the sum of the Lie algebras of all the 
different stabilizer subgroups of $T$. If $X\in\got{t}_x$, 
then $X_M(x)=0$ and 
it follows from Lemma \ref{constlem} that 
$\sigma^{\got{t}}(X,\, Y)=0$ for every $Y\in\got{t}$. 
This shows that $\sigma^{\got{t}}(X,\, Y)=0$ for every 
$X\in\got{t}'$ and every $Y\in\got{t}$. 
Now ii) means that ${\got t}'
=\got{t}$, hence $\sigma^{\got{t}}=0$, and  
Lemma \ref{isotropicorbitlem} implies that every 
principal orbit is a Lagrange submanifold of 
$(M,\,\sigma )$, and therefore coisotropic. 
It follows that we may apply Proposition \ref{M/Tprop} 
with $\got{t}_{\scriptop{h}}=\got{t}'=\got{t}$, 
and conclude that 
$\Phi _p$ is a homeomorphism from the Delzant polytope  
$\Delta$ onto $M/T$. 

iii) $\Rightarrow$ iv) because any convex polytope  
is contractible. 

iv) $\Rightarrow$ v) follows from  
Corollary \ref{Hamcor}. 

Finally v) $\Rightarrow$ i) follows from 
the fact that the image of the momentum mapping is equal to the 
convex hull of the images under the momentum mapping of the 
fixed points, cf. Atiyah \cite[Th. 1]{atiyah} or 
Guillemin and Sternberg \cite[Th. 4]{gs}. 
\end{proof}

The implication i) $\Rightarrow$ v) has also been obtained 
by Giacobbe \cite[Th. 3.13]{giacobbe}. 

Note that if the conditions i) -- v) in Corollary \ref{delzantprop} 
hold, then $(M,\,\sigma )$ together with the $T$\--action on $M$ 
is a Delzant manifold, and $M/T$ is the corresponding Delzant polytope. 
If a compact Lie group $K$ acts linearly and continuously 
on a vector space $V$, then the {\em average} of $v\in V$ over 
$K$ is defined as 
\[
\int _Kk\cdot v\, m(\op{d}\! k)/m(K),
\]
in which $m$ denotes any Haar measure on $K$. 

\begin{corollary}
With the notation of Proposition {\em \ref{M/Tprop}}, let 
$\pi _{N/P}:M/T\to N/P$ be the mapping 
$\Phi_p^{-1}$ followed by the projection from 
$\Delta\times (N/P)$ onto the 
second factor. Let  
$\iota_p:N/P\to M/T$ 
be defined by $\iota_p(\zeta +P)=p+\zeta$. 
Then we have the following conclusions. 

For each nonnegative integer $k$, 
the mapping $\pi _{N/P}^*:\op{H}^k(N/P,\,\R )\to 
\op{H}^k(M/T,\,\R )$ is an isomorphism, with inverse 
equal to $\iota _p^*$.  

The mapping which assigns to any 
$\lambda\in\Lambda ^kN^*$ 
the cohomology class of the constant $k$\--form $\lambda$ 
on $N/P$ is an isomorphism 
from $\Lambda ^kN^*$ onto 
$\op{H}^k(N/P,\,\R )$, and every 
closed $k$\--form on $N/P$ is 
cohomologous to its average over the torus $N/P$.  
\label{cohomcor}
\end{corollary}
\begin{proof}
The first statement follows because $\Delta$ is a convex subset 
of $\got{t}^*$ and hence it is contractible. 
The second statement is a well\--known characterization of the 
cohomology of tori. The fact that a closed differential form on a compact 
connected Lie group is cohomologous to its average goes back to 
\'Elie Cartan \cite{cartan}. 
\end{proof}

Any finite\--dimensional vector space $W$ carries a 
positive translation\--invariant measure $m$, which 
is unique up to a positive factor. For any non\--negligible 
compact subset 
$A$ of $W$, the {\em center of mass} of $A$ is defined as 
\[
\int_A\, x\, m(\op{d}\! x)/m(A)\in W,  
\]
which is independent of the choice of the 
positive translation\--invariant measure $m$ on $W$. 
\begin{corollary}
Let $X\in\got{t}$. Then $X_M$ is Hamiltonian if and 
only if $X\in\got{t}_{\scriptop{h}}$. Furthermore, 
the image of any momentum mapping of the 
Hamiltonian action of $T_{\scriptop{h}}$ on $M$ 
is equal to a translate of the Delzant polytope $\Delta$ in 
Proposition {\em \ref{M/Tprop}}, where we note that any two 
momentum mappings for the same torus action differ by a 
constant element of ${\got{t}_{\scriptop{h}}}^*$. The translational ambiguity of 
$\Delta$ can be removed by putting the center of mass of 
$\Delta$ at the origin. Here a momentum mapping 
for the Hamiltonian action of $T_{\scriptop{h}}$ is a smooth 
${\got{t}_{\scriptop{h}}}^*$\--valued function $\mu$ on $M$ such that 
for every $X\in\got{t}_{\scriptop{h}}$ the $X$\--component of $\op{d}\!\mu$ 
is equal to $-\op{i}_{X_M}\sigma$.  
\label{Tstabcor}
\end{corollary}
\begin{proof}
It follows from Corollary \ref{Hamcor} and 
Corollary \ref{cohomcor} that 
the vector field $X_M$ is Hamiltonian if and 
only if $[\widehat{\sigma} (X)]=0$ if and only if 
$[\iota_p^*(\widehat{\sigma}(X))]=\iota_p^*[\widehat{\sigma}(X)]=0$. 
Now constant one\--forms on $N/P$ 
are canonically identified with linear forms on 
$N=(\got{l}/\got{t}_{\scriptop{h}})^*$, 
which are identified with elements of $\got{l}/\got{t}_{\scriptop{h}}$. 
With this identification, $\iota_p^*(\widehat{\sigma}(X))$ corresponds 
to $X+\got{t}_{\scriptop{h}}$, which is equal to zero 
if and only if $X\in\got{t}_{\scriptop{h}}$. 

The second statement in the corollary follows from the fact that 
if $\mu$ is a momentum mapping for the Hamiltonian $T_{\scriptop{h}}$\--action, 
then $\op{d}(\mu (X))= \widehat{\sigma}(X)$ for every 
$X\in\got{t}_{\scriptop{h}}$. In other words, $\mu$ differs from the  
$\got{t}_{\scriptop{h}}$\--component of any canonical local chart on $M/T$ 
by a constant vector in ${\got{t}_{\scriptop{h}}}^*$.  
Therefore the image of $\mu$ 
corresponds to $\Delta\simeq(M/T)/N$, 
the orbit space  of the translational $N$\--action 
on $M/T$. Here we use that restriction to $\got{t}_{\scriptop{h}}$ of 
linear forms on $\got{l}$ leads to a canonical identification 
of $\got{l}^*/(\got{l}/\got{t}_{\scriptop{h}})^*$ with 
${\got{t}_{\scriptop{h}}}^*$.  
\end{proof}

McDuff \cite{M} proved that a symplectic 
circle action on a four\--dimensional compact connected 
symplectic manifold is Hamiltonian, if and only if it has 
a fixed point, but that in higher dimensions there exist 
non\--Hamiltonian symplectic circle actions with fixed points. 
Corollary \ref{Tstabcor} follows from \cite{M} 
if $\op{dim}M=4$, but not if $\op{dim}M=2n>4$. 
Our proof of Corollary \ref{Tstabcor} uses in an essential 
way that $X_M$ is an infinitesimal action of a symplectic  
action of an $n$\--dimensional torus with a Lagrange orbit. 
\begin{remark}
Because a Hamiltonian torus action has fixed points, it follows 
from Corollary \ref{Tstabcor} that the action of $T_{\scriptop{h}}$ on 
$M$ has fixed points, that is, there exist $x\in M$ such that 
$T_{\scriptop{h}}\subset T_x$, hence $T_{\scriptop{h}}=T_x$ because 
the definition of $T_{\scriptop{h}}$ in Lemma \ref{stabilizerlem} implies that 
$T_x\subset T_{\scriptop{h}}$ 
for every $x\in M$. In other words, $T_{\scriptop{h}}$ can also be characterized 
as the unique maximal stabilizer subgroup of $T$. 

Actually the fixed points in $M$ for the action of $T_{\scriptop{h}}$
are the $x\in M$ such that $\mu (x)$ is a vertex of 
the Delzant polytope $\Delta$, where $\mu :M\to\Delta\subset 
{\got{t}_{\scriptop{h}}}^*$ denotes the momentum map of the Hamiltonian 
$T_{\scriptop{h}}$\--action. 
\label{T1rem}
\end{remark}
\begin{remark}
Let $\got{t}_{\scriptop{h}}\neq\got{t}$. It follows from 
Lemma \ref{stabilizerlem} that for every 
$X\in\got{t}\setminus\got{t}_{\scriptop{h}}$ 
the vector field $X_M$ has no zeros in $M$, and we conclude 
that the Euler characteristic $\chi (M)$ of $M$ 
is equal to zero. 

Furthermore the localization formula 
of Berline\--Vergne and Atiyah\--Bott in equivariant 
cohomology, in the form of \cite[(4.13)]{duis}, yields 
for every $T$\--equivariantly closed $T$\--equivariant
differential form $\omega$ on $M$ that the integral of $\omega$ 
over $M$ is equal to zero, when evaluated at 
$X\in\got{t}\setminus\got{t}_{\scriptop{h}}$. 
Because $\got{t}\setminus\got{t}_{\scriptop{h}}$ is dense in $\got{t}$, 
it follows that the integral over $M$ of each $T$\--equivariantly closed 
$T$\--equivariant differential form is identically equal to zero. 
If $X\in\got{t}_{\scriptop{h}}$, then Lemma \ref{Tstabcor} implies that 
$X_M$ is Hamiltonian, and the zeros of $X_M$ are the critical points of 
its Hamiltonian function, which form a non\--empty subset of $M$. 
In this case the localization formula \cite[(4.13)]{duis} 
yields that the sum over the connected components $F$ of the 
zeroset of $X_M$ of the integrals over $F$ of 
$\omega (X)/\varepsilon (X)$ is equal to zero. 
The generalization of Ginzburg \cite[Th. 6.1]{ginzburg} 
of the Duistermaat\--Heckman formula is related to these 
observations. On the other hand the integral over $M$ of 
a $T_{\scriptop{h}}$\--equivariantly closed 
$T_{\scriptop{h}}$\--equivariant differential 
form, such as $\got{t}_{\scriptop{h}}\ni X\mapsto\op{e}^{\op{i}(\mu (X)-\sigma )}$,  
is usually nonzero. 

If $\got{t}_{\scriptop{h}}=\got{t}$, then 
it follows from Corollary \ref{delzantprop} and 
Corollary \ref{Hamcor} 
that $(M,\,\sigma ,\, T)$ is a Delzant manifold, and  
$\chi (M)$ is equal to the number of vertices of the 
Delzant polytope $\Delta$. This can be proved by observing 
that for a generic $X\in\got{t}$ 
the momentum map is bijective from 
the zeroset of $X_M$ to the set of vertices of $\Delta$, 
and each zero of $X_M$ has Poincar\'e index equal to one. 
See also Guillemin \cite[Exerc. 4.15]{guillemin}. 
\label{eulerrem}
\end{remark}

\section{Two lemmas}
\label{twolemmasec}
The following lemmas will be used later in the paper. 
Lemma \ref{basicformlem} is used in the proof of 
Proposition \ref{liftprop}, whereas Lemma \ref{hslem} 
is used in the proof of Lemma \ref{LPhilem} and Lemma \ref{TPlem}. 
The proofs of Lemma \ref{basicformlem} and Lemma \ref{TPlem} 
are based on the local models of Lemma \ref{modellem}. 

Throughout this section , $(M,\,\sigma )$ is 
a symplectic manifold with an effective symplectic 
action of a torus $T$ with coisotropic principal orbits. 
\begin{lemma} 
Let $X_j$, $1\leq j\leq\op{dim}\got{l}$, be a basis of $\got{l}$. 
The basic $k$\--forms on $M$ are the $k$\--forms 
\begin{equation}
\omega =\sum_{j_1<\ldots <j_k}
\, f_{j_1,\,\ldots ,\, j_k}\,\widehat{\sigma}(X_{j_1})\wedge\ldots\wedge 
\widehat{\sigma} (X_{j_k})
\label{basicform}
\end{equation}
in which $f_{j_1,\,\ldots ,\, j_k}\in\op{C}^{\infty}(M)^T$ 
and $\widehat{\sigma}$ is defined as in Lemma \ref{alphalem}. 
\label{basicformlem}
\end{lemma}
\begin{proof}
Because the one\--forms $\widehat{\sigma} (X_j)$ are basic, any 
$\omega$ as in (\ref{basicform}) with 
$f_{j_1,\,\ldots ,\, j_k}\in\op{C}^{\infty}(M)^T$ is a basic form. 

Using partitions of unity with elements of $\op{C}^{\infty}(M)^T$, 
it is sufficient to prove the converse statement in a 
local model as in Lemma \ref{modellem}. 
Let $\omega$ be a basic $k$\--form. 
In the principal stratum where $\rho _j>0$ for every 
$1\leq j\leq m$, we have that 
\[
\omega =\sum_{l=0}^k\,\sum_{j_1<\ldots <j_l}
\, f^{k-l}_{j_1,\,\ldots ,\, j_l}(\rho )\, 
\op{d}\!\rho _{j_1}\wedge\ldots\wedge 
\op{d}\!\rho _{j_l},
\]
in which the $f^{k-l}_{j_1,\,\ldots ,\, j_l}(\rho )$ are uniquely determined 
smooth $(k-l)$\--forms on $(\got{l}/\got{h})^*$, depending smoothly on 
$\rho$. We are done if we can prove that 
the $f^{k-l}_{j_1,\,\ldots ,\, j_l}$ extend smoothly 
over the boundary where some of the $\rho _j$ are equal to zero. 

Recall that $\rho ^j=((p^j)^2+(q^j)^2)/2$, if $z^j=p^j+\op{i}q^j$ 
with $p^j,\, q^j\in\R$.  
Then $\op{d}\! \rho _j=p^j\,\op{d}\! p^j+q^j\,\op{d}\! q^j$ shows that 
$\op{d}\!\rho _{j_1}\wedge\ldots\wedge 
\op{d}\!\rho _{j_l}$ 
has the component 
\[
(\prod_{i=1}^l\, q^{j_i})
\,\op{d}\! q^{j_1}\wedge\ldots\wedge 
\op{d}\! q^{j_l},   
\]
and therefore the smoothness of $\omega$ implies that 
for each $0\leq l\leq k$ and each sequence 
$j_1,\,\ldots ,\, j_l$ with 
$j_1<\ldots <j_l$ the form  
\[
(\prod_{i=1}^l\, q^{j_i})\, f^{k-l}_{j_1,\,\ldots ,\, j_l}(\rho )\, 
\]
depends smoothly on $(p,\, q)$. Applying the 
differential operator $\partial ^l/\partial q^{j_1}\,\ldots 
\partial q^{j_l}$ and putting $q=0$, we obtain that 
\[
f^{k-l}_{j_1,\,\ldots ,\, j_l}((p^1)^2/2,\,\ldots ,\, (p^l)^2/2)
\]
depends smoothly on $p$, and moreover is 
invariant under each of the reflections $p^j\mapsto -p^j$. 
Whitney \cite{whitney} proved that this implies 
that the function $f^{k-l}_{j_1,\,\ldots ,\, j_l}$ extends smoothly 
over the boundary where some of the $\rho _j$ are equal to zero. 
\end{proof}

For any smooth mapping $f$ from a smooth manifold $M$ to a 
smooth manifold $N$, the tangent mapping $\op{T}_x\! f$ is the 
linear mapping from $\op{T}_xM$ to $\op{T}_{f(x)}N$ which in 
local coordinates corresponds to the Jacobi matrix of $f$ 
at the point $x$. 
\begin{lemma}
Let $\Phi :M\to M$ be a $T$\--equivariant diffeomorphism 
which preserves the $T$\--orbits. Then there is a unique smooth 
$T$\--invariant mapping $\tau :M\to T$ such that 
$\Phi (x)=\tau (x)\cdot x$ for every $x\in M$. 

If $\Phi$ preserves the symplectic form $\sigma$, then 
$(\op{T}_x\!\tau )(v)\in\got{l}$ for each $x\in M$ and $v\in\op{T}_x\! M$. 
Here $\got{l}$ is the kernel of the antisymmetric 
bilinear form $\sigma ^{\got{t}}$ introduced in 
Lemma \ref{constlem}, and we identify  
each tangent space of $T$ with $\got{t}$. 
\label{hslem}
\end{lemma}
\begin{proof}
The first statement has been proved for arbitrary torus actions on orbifolds 
by Haefliger and Salem \cite[Th. 3.1]{hs}, but in our case the proof is 
elementary. The statement is obvious if we replace $M$ by 
the set $M_{\scriptop{reg}}$ on which the action is free and 
defines a principal $T$\--fibration, and it remains to be proved that 
$\tau$ has a smooth extension to $M$. In the local model of Lemma 
\ref{modellem}, we have 
\[
\Phi :(k,\,\lambda ,\, z)
\mapsto (\tau (\lambda ,\,\rho )_K\, k,\,\lambda ,\, 
\iota (\tau (\lambda ,\,\rho )_H)\cdot z). 
\]
The smoothness of $\Phi$ implies that 
$(\lambda ,\,\rho )\mapsto \tau (\lambda ,\,\rho )_K$ 
has a smooth extension. Write $\widehat{\tau}(\lambda ,\rho )
=\iota (\tau (\lambda ,\,\rho )_H)\in\T ^m$. It remains to 
prove that the fact that 
$\Psi :(\lambda ,\, z)\mapsto 
\widehat{\tau}(\lambda ,\,\rho)\cdot z$ has a smooth extension,
implies that $\widehat{\tau}$ has a smooth extension, because 
the fact that $\iota :H\to\T ^m$ is an isomorphism of Lie 
groups then implies that $\tau _H$ has a smooth extension. 
 
Now the function 
\[
f^j(\lambda ,\, z):=\Psi (\lambda ,\, z)^j\,\overline{z^j}=
\widehat{\tau}(\lambda ,\,\rho)^j\, |z^j|^2
\]
has a smooth extension, of which the restriction to the ``real 
domain'' $q=0$ is an even function in each of the variables 
$p^j$. It therefore follows from Whitney \cite{whitney} 
that there is a smooth function $g^j$ such that 
$f^j(\lambda ,\, z)=g^j(\lambda ,\,\rho )$. However 
$g^j(\lambda ,\,\rho )=0$ when 
$\rho _j=0$, and it follows that 
\[
g^j(\lambda ,\,\rho )=\int_0^1\,\partial 
g^j(\lambda ,\,\rho _1,\,\ldots ,\, t\,\rho _j,\,\ldots ,\,\rho _n)/
\partial t\,\op{d}\! t=h^j(\lambda ,\,\rho )\,\rho _j,
\]
in which 
\[
(\lambda ,\,\rho )\mapsto h^j(\lambda ,\,\rho )=\int_0^1\,\left.\partial 
g^j(\lambda ,\,\rho _1,\,\ldots ,\, r_j,\,\ldots ,\,\rho _n)
/\partial r_j\right| _{r_j=t\,\rho _j}\,\op{d}\! t
\]
is smooth. Because 
\[
h^j(\lambda ,\,\rho )\,\rho _j=g^j(\lambda ,\,\rho )=f^j(\lambda ,\, z)
=2\, \widehat{\tau}(\lambda ,\,\rho)^j\,\rho _j,
\]
it follows that $\widehat{\tau}
=h^j/2$ when $\rho _j>0$, which extends smoothly 
over the boundary $\rho _j=0$. 

Write, for each $x\in M$,  
$\tau '_x:=\op{T}_x\!\tau$, viewed as a linear mapping from 
$\op{T}_x\! M$ to $\got{t}$, and 
${\tau (x)_M}':=\op{T}_x(\tau (x)_M)$, which is a 
symplectic linear 
mapping from $\op{T}_x\! M$ to $\op{T}_{\Phi (x)}\! M$. 
Then it follows from the sum rule for differentiation 
of an expression in which a variable occurs at several places, that 
\begin{equation}
(\op{T}_x\!\Phi ) \, v={\tau (x)_M}'\, v+(\tau '_x\, v)_M(\Phi (x)),
\quad v\in\op{T}_x\! M.
\label{TPhitau}
\end{equation}
If $X\in\got{t}$, then the $T$\--equivariance of $\Phi$ 
implies that $(\op{T}_x\!\Phi )\, X_M(x)=X_M(\Phi (x))$. 
On the other hand, the commutativity of $T$ implies that 
$\tau (x)_M(t\cdot x)=t\cdot\tau (x)_M(x)=t\cdot\Phi (x)$ 
for every $t\in T$, 
and differentiating this with respect to $t$ at $t=1$ in the 
direction of $X$, we obtain ${\tau (x)_M}'(X_M(x))=X_M(\Phi (x))$. 
The condition $\sigma =\Phi ^*\,\sigma$ implies that we have, for every 
$x\in M$, $v\in\op{T}_x\! M$, and $X\in\got{t}$,   
\begin{eqnarray*}
\sigma _x(v,\, X_M(x))
&=&\sigma _{\Phi (x)}((\op{T}_x\!\Phi ) \, v,\, (\op{T}_x\!\Phi )\, X_M(x))\\
&=&\sigma _{\Phi (x)}({\tau (x)_M}'\, v+(\tau '_x\, v)_M(\Phi (x))
,\, X_M(\Phi (x)))\\ 
&=&\sigma _{\Phi (x)}({\tau (x)_M}'\, v,\,{\tau (x)_M}'\, X_M(x))
+\sigma ^{\got{t}}(\tau '_x\, v,\, X),
\end{eqnarray*}
which implies that $\sigma ^{\got{t}}(\tau '_x\, v,\, X)=0$ 
because ${\tau (x)_M}'$ is symplectic. Because 
$\sigma ^{\got{t}}(X,\,\tau '_x\, v)=0$ for every $X\in\got{t}$, 
it follows that $\tau '_x\, v\in\got{l}:=\op{ker}\sigma ^{\got{t}}$. 
\end{proof}
\begin{remark}
One can prove that $\Phi$ is a $T$\--equivariant symplectomorphism 
of $(M,\,\sigma )$ which preserves the $T$\--orbits, if and 
only if for every $x\in M$ there exists a $T$\--invariant 
open neighborhood $U$ of $x$ in $M$, a $T$\--invariant smooth function 
$f$ on $U$, and an element $t\in T$, such that 
$\Phi =\op{e}^{\fop{Ham}_f}\circ t_M$ on $U$. 
The ``if'' part follows from 
Remark \ref{multfreerem}. 
\label{hsrem}
\end{remark}

\section{Lifts}
\label{liftsec}
If we identify each of the tangent spaces of $(M/T)_{\scriptop{reg}}$ 
with $\got{l}^*$ as in Remark \ref{regrem}, then any 
$\xi\in\got{l}^*$ can be viewed as a constant vector 
field on $(M/T)_{\scriptop{reg}}$. A vector field $L_{\xi}$ 
in $M_{\scriptop{reg}}$ is called a {\em lift} of $\xi$, 
if $\op{T}_x\!\pi (L_{\xi}(x))=\xi$ for all $x\in M_{\scriptop{reg}}$. 
Here the tangent mapping $\op{T}_x\pi :\op{T}_x\! M_{\scriptop{reg}}\to
\op{T}_{\pi (x)}(M/T)_{\scriptop{reg}}$ of $\pi$ is identified 
with the linear mapping $\widehat{\sigma}_x:\op{T}_x\! M\to
\got{l}^*$, defined by the $\got{l}^*$\--valued one\--form 
$\widehat{\sigma}$ on $M$. In view of the definition 
of $\widehat{\sigma}$ in Lemma \ref{alphalem} and 
(\ref{widehatsigma}), the condition that $L_{\xi}$ is a lift of 
$\xi$ therefore is equivalent to 
\begin{equation}
\sigma (L_{\xi},\, X_M)=\xi (X),
\quad \xi\in\got{l}^*,\quad X\in\got{l}. 
\label{lift}
\end{equation} 

If $L_{\xi}$, $\xi\in\got{l}^*$, is a family of 
smooth $T$\--invariant vector fields on $M_{\scriptop{reg}}$, 
which depends linearly on $\xi$ and are lifts in the 
sense of (\ref{lift}), then for each $x\in M_{\scriptop{reg}}$ 
the vectors $L_{\xi}(x)$, $\xi\in\got{l}^*$, span a linear subspace 
$H_x$ of $\op{T}_x\! M$ which is complementary to the 
tangent space $\got{t}_M(x)$ at $x$ of the orbit $T\cdot x$.  
The $H_x$, $x\in M_{\scriptop{reg}}$, are the horizontal 
spaces for a unique $T$\--invariant 
{\em infinitesimal connection} $\nabla$ for the principal $T$\--bundle 
$\pi :M_{\scriptop{reg}}\to (M/T)_{\scriptop{reg}}$. 
This connection is $T$\--invariant, if and only if each of the 
lifts $L_{\xi}$, $\xi\in\got{l}^*$, is $T$\--invariant. 

Conversely, if we have given a $T$\--invariant infinitesimal 
connection $\nabla$ for the principal $T$\--fibration in 
$M_{\scriptop{reg}}$, with horizontal spaces 
$H_x=H^{\nabla}_x$, $x\in M_{\scriptop{reg}}$, 
then we have for each $\xi\in\got{l}^*$ a unique lift 
$L_{\xi}$ of $\xi$ such that $L_{\xi}(x)\in H_x$ 
for every $x\in M_{\scriptop{reg}}$. $L_{\xi}$ is 
called the {\em horizontal lift} of $\xi$ defined by 
the connection $\nabla$, and denoted by 
$\xi _{\scriptop{hor}}^{\nabla}$ in the literature 
on connections. Because the mapping 
$\xi\mapsto\xi _{\scriptop{hor}}^{\nabla}$ is linear, 
``lifts $L_{\xi}$ which depend linearly on $\xi$''  
and ``connections'' are equivalent objects. 
We will use the somewhat simpler notation $L_{\xi}$ 
instead of $\xi _{\scriptop{hor}}^{\nabla}$, because it is 
the lifts which we will be using to construct our global model.  

In this section we construct lifts $L_{\xi}$, $\xi\in\got{l}^*$,  
depending linearly on $\xi$, which are  admissible in the sense 
of Definition \ref{liftdef}, 
and have Lie brackets and symplectic products which are as simple 
as we can get them. See Proposition \ref{liftprop} below. 
This construction is based on a computation in the 
cohomology of the closed basic differential 
forms on $M$, which according 
to the theorem of Koszul \cite{koszul} is canonically 
isomorphic to the sheaf (= \v{C}ech) cohomology of the orbit space $M/T$ 
with values in $\R$. 
The lifts in Proposition \ref{liftprop} form the core of 
the construction of the model for the symplectic $T$\--manifold 
$(M,\,\sigma ,\, T)$, given 
in Proposition \ref{Aprop} and Proposition \ref{omegalem}.

\subsection{Admissible connections}
\label{admliftss}
\begin{definition}
Let, in the local model of Lemma \ref{modellem} with the diffeomorphism 
$\Phi$, the lift $L_{\xi}$ 
be equal to the image under $\op{T}\!\Phi$ 
of the vector field 
$(X_{\xi},\,\delta\lambda _{\xi},\,\delta z_{\xi})$. 
Then, in terms of the $(\theta ,\,\rho)$\--coordinates in 
$\C ^m$, we obtain in view of (\ref{modelsigma}) 
and (\ref{Aiota}) that the equation (\ref{lift}) 
is equivalent to $A(\delta\lambda _{\xi},\,\delta\rho _{\xi})=\xi$. 
Let $(\delta\lambda _{\xi},\,\delta\rho _{\xi})=A^{-1}(\xi )$, 
and let $L_{\xi}^{\Phi}$ be the image under $\op{T}\!\Phi$ 
of the ``constant'' vector field 
$(0,\,\delta\lambda _{\xi},\, (0,\,\delta\rho _{\xi}))$, 
where we use the $(\theta ,\,\rho)$\--coordinates in $\C ^m$. 
Then $L_{\xi}^{\Phi}$ is a smooth $T$\--invariant 
vector field on $U\cap M_{\scriptop{reg}}$, and a lift of $\xi$. 
We call $L_{\xi}^{\Phi}$ the {\em local model lift defined 
by the local model with the diffeomorphism $\Phi$}. 

The local model lift $L_{\xi}^{\Phi}$ extends to a smooth $T$\--invariant 
vector field on $U$ when $\delta\rho _{\xi}=0$, that is, when 
$\xi =0$ on $\got{h}$. 
On the other hand, if we write $r_j=|z_j|$, then 
$\partial/\partial\rho _j = (1/r_j)\, \partial/\partial r _j$.  
This shows that $L^{\Phi}_{\xi}$ has 
a pole singularity at any point $(k,\,\lambda ,\, z)$ for which there 
exists a $1\leq j\leq m$ such that $z^j = 0$ and $\xi (X_j)\neq 0$. 
\label{localmodelliftdef}
\end{definition} 

\begin{lemma}
Let $\widetilde{\Phi}:\widetilde{K}\times\widetilde{E}_0\to
\widetilde{U}$ be another local model as in Lemma \ref{modellem}, 
where we use the same projection $X\mapsto X_{\got{l}}:\got{t}
\to\got{l}$. Then there is a smooth $T$\--invariant mapping 
$\alpha :U\cap\widetilde{U}\to\got{l}$, such that 
$L^{\widetilde{\Phi}}_{\xi}(x)=L^{\Phi}_{\xi}(x)+
\alpha (x)_M(x)$ for every 
$x\in U\cap\widetilde{U}\cap M_{\scriptop{reg}}$. 
Here $L^{\Phi}_{\xi}$ and $L^{\widetilde{\Phi}}_{\xi}$ 
are the local model lifts Definition \ref{localmodelliftdef}. 
\label{LPhilem}
\end{lemma}
\begin{proof}
Let $x_0\in U\cap\widetilde{U}$ and 
write $(k_0,\,\lambda _0,\, (\theta _0,\,\rho _0))=\Phi ^{-1}(x_0)$, 
where we use the $(\theta ,\,\rho )$\--coordinates in $\C ^m$.  
By permuting the coordinates in $\C ^m$, we can arrange 
that $(\rho _0)_j=0$ for $1\leq j\leq m_0$ and 
$(\rho _0)_j>0$ for $m_0<j\leq m$. Then $H_0:=T_{x_0}$ 
is equal to the subgroup $\iota ^{-1}(\T ^{m_0}\times\{ 1\})$ 
of $H$, where $\iota$ denotes the isomorphism 
from $H$ onto $\T ^m$, introduced in 
Remark \ref{complexrem}. Here $H=T_x$ as in 
Lemma \ref{modellem}. 
Let $H_0':=\iota ^{-1}(\{ 1\}\times\T ^{m-m_0})$. 
Then $H_0'$ is a complementary subtorus to $H_0$ in $H$,  
and $K_0:=H_0'\, K$ is a complementary subtorus to 
$H_0$ in $T$ which contains $K$. Let $(\theta',\,\rho')$ 
and $(\theta'',\,\rho'')$ be the 
first $m_0$ and the last $m-m_0$ of the $(\theta,\,\rho)$-coordinates, 
respectively. Then the rotation of $z^j$ over $(\theta '')^j$, 
for each $m_0<j\leq m$, defines an element $R(\theta '')$ 
of $\{ 1\}\times\T ^{m-m_0}$, and $\iota ^{-1}(R(\theta ''))\in H_0'$. 

On the other hand 
\[
\Lambda  _0(\rho''):X\mapsto\sum_{j=m_0+1}^m\,\rho _j\, 
\iota (X)/\op{i}
\]
is a linear form on $\got{h}$ which is equal to zero 
on the Lie algebra $\got{h}_0$ of $H_0$. This linear from 
has a unique extension 
to a linear form $\Lambda (\rho'')$ on $\got{l}$ which 
is equal to zero on $\got{l}\cap\got{k}$. In this way we 
obtain an element 
$\Lambda(\rho'')\in (\got{l}/\got{h}_0)^*$. A 
straightforward computation shows that the mapping 
\[
\Psi :(k,\,\lambda ,\, (\theta ,\,\rho ))\mapsto 
(\iota ^{-1}(R(\theta ''))\, k,\,\lambda -\lambda _0
+\Lambda (\rho ''-\rho ''_0),
\, (\theta ',\,\rho '))
\]
when restricted to the the domain where $\rho _j>0$ for all 
$m_0<j\leq m$, defines a smooth $T$\--equivariant 
symplectomorphism from $K\times (\got{l}/\got{h})^*\times 
\C ^m$ to $K_0\times (\got{l}/\got{h}_0)^*\times\C ^{m_0}$. 
Moreover, $\Psi\circ\Phi ^{-1}(x_0)$ belongs to the 
$T$\--orbit of $(1,\, 0,\, 0)$ in 
$K_0\times (\got{l}/\got{h}_0)^*\times\C ^{m_0}$. 
Because the tangent mapping of $\Psi$ maps 
$(0,\,\delta\lambda ,\, (0,\,\delta\rho ))$ to 
$(0,\,\delta\lambda +\Lambda (\delta\rho ''),\, (0,\,\delta\rho '))$, 
we have $L^{\Phi\circ\Psi ^{-1}}_{\xi}=L^{\Phi}_{\xi}$.  

Similarly we have a smooth $T$\--equivariant symplectomorphism 
$\widetilde{\Psi}$ from a $T$\--invariant open neighborhood 
of $\widetilde{\Phi}^{-1}(x_0)$ in $\widetilde{K}\times 
(\got{l}/\widetilde{\got{h}})^*\times\C ^{\widetilde{m}}$ 
onto a $T$\--invariant open neighborhood of 
$(1,\, 0,\, 0)$ in $\widetilde{K}_0\times (\got{l}/\got{h}_0)^*
\times\C ^{m_0}$, such that 
$\widetilde{\Psi}\circ\widetilde{\Phi}^{-1}(x_0)\in T\cdot (1,\, 0,\, 0)$ and 
$L^{\widetilde{\Phi}\circ\widetilde{\Psi}^{-1}}_{\xi}
=L^{\widetilde{\Phi}}_{\xi}$. Here $\widetilde{K}_0$ is 
another complementary subtorus to $H_0$ in $T$. 

The mapping  
\[
\Xi :(k,\,\lambda, z)\mapsto (k_{\widetilde{K}},\,\lambda ,\, 
\iota (k^{\widetilde{K}_0}_{H_0})\cdot z):
K_0\times (\got{l}/\got{h}_0)^*\times\C ^{m_0}
\to\widetilde{K}_0\times (\got{l}/\got{h}_0)^*\times\C ^{m_0}
\]
is a $T$\--equivariant symplectomorphism which maps 
$(1,\, 0,\, 0)$ to $(1,\, 0,\, 0)$. Here we have written, 
for each $k\in K_0$, 
$k=k_{\widetilde{K}_0}\, k_{H_0}^{\widetilde{K}_0}$ with 
$k_{H_0}^{\widetilde{K}_0}\in H_0$ and 
$k_{\widetilde{K}_0}\in \widetilde{K}_0$. 
Because $h=k_{H_0}^{\widetilde{K}_0}$ is the unique 
element in $H_0$ such that 
$k_{\widetilde{K}_0}:=k\, h^{-1}\in\widetilde{K}_0$, the fact that 
$\Xi$ is a $T$\--equivariant symplectomorphism follows from  
the proof of Lemma \ref{modellem}, with $(H,\, K)$ replaced by 
$(H_0,\, K_0)$ and by $(H_0,\, \widetilde{K}_0)$, respectively.  
Because the tangent mapping of $\Xi$ maps 
$(0,\,\delta\lambda ,\,\delta z)$ to 
$(0,\,\delta\lambda ,\,\delta z)$, we have 
that 
\[
L_{\xi}^{\widetilde{\Phi}\circ\widetilde{\Psi}^{-1}\circ\Xi} 
=L_{\xi}^{\widetilde{\Phi}\circ\widetilde{\Psi}^{-1}}
=L_{\xi}^{\widetilde\Phi}.
\]

The mapping $\Theta :=
\Psi\circ\Phi ^{-1}\circ\widetilde{\Phi}
\circ\widetilde{\Psi}^{-1}\circ\Xi$
is a smooth $T$\--equivariant symplectomorphism from an 
open $T$\--invariant neighborhood of $(1,\, 0,\, 0)$ 
in $K_0\times (\got{l}/\got{h}_0)^*\times\C ^{m_0}$,  
onto an 
open $T$\--invariant neighborhood of $(1,\, 0,\, 0)$ 
in $K_0\times (\got{l}/\got{h}_0)^*\times\C ^{m_0}$, 
which preserves the $T$\--orbit of $(1,\, 0,\, 0)$. 
Recall the $\got{l}^*$\--valued one form 
$\widehat{\sigma}$ defined in Lemma \ref{alphalem}, 
which we used to identify 
all tangent spaces of the orbit space with $\got{l}^*$. 
Because every $T$\--equivariant symplectomorphism 
preserves $\widehat{\sigma}$, 
its induced transformation 
of the orbit space has derivative equal to the 
identity at every point. Therefore $\Theta$ 
is a translation on each connected open subset of 
the $T$\--orbit space by means of a constant element $v$ of 
$\got{l}^*$. Because $\Theta$ preserves the 
$T$\--orbit of $(1,\, 0,\, 0)$, we have $v=0$ on the 
connected component of $(1,\, 0,\, 0)$ 
of the domain of definition $\Upsilon$ of $\Theta$. That is, 
$\Theta$ preserves all the $T$\--orbits in a $T$\--invariant 
open neighborhood of $(1,\, 0,\, 0)$. 

It now follows from Lemma \ref{hslem} that there there is a smooth 
$T$\--invariant $\got{l}$\--valued function $\tau$ 
on $\Upsilon$, such that 
$\Theta (\upsilon )=\tau (\upsilon )\cdot\upsilon$ 
and $\op{T}_{\upsilon}\!\tau (\delta\upsilon )\in\got{l}$ 
for every $\upsilon\in\Upsilon$ and 
$\delta\upsilon\in\op{T}_{\upsilon}\!\Upsilon$. 
It follows from (\ref{TPhitau}), with $\Phi$ and 
$\upsilon$ replaced by $\Theta$ and $\delta\upsilon 
:=(0,\,\delta\lambda ,\,\delta z)$, 
respectively, where $\delta\theta =0$ and 
$(\delta\lambda ,\,\delta\rho )=A^{-1}\xi$, that 
$\op{T}\!\Theta$ maps the vector field 
$(0,\,\delta\lambda ,\,\delta z)$ to the sum 
of $(0,\,\delta\lambda ,\,\delta z)$ and 
\[
((\tau '_v\,\delta v)_{\got{k}_0},\, 0,\, 
\iota ((\tau '_v\,\delta v)_{\got{h}_0} ).
\]
Because $\widetilde{\Phi}\circ\widetilde{\Psi}^{-1}\circ\Xi 
=\Phi\circ\Psi ^{-1}\circ\Theta$, 
$L_{\xi}^{\widetilde{\Phi}\circ\widetilde{\Psi}^{-1}\circ\Xi }
=L_{\xi}^{\widetilde{\Phi}}$, and 
$L^{\Phi\circ\Psi ^{-1}}_{\xi}=L^{\Phi}_{\xi}$, the conclusion of 
the lemma follows with $\alpha (\widetilde{x})=\tau '_\upsilon
\,\delta\upsilon$ if 
$\widetilde{x}=\Phi\circ\Psi ^{-1}(\upsilon )$. 
\end{proof}
\begin{definition}
We use the atlas of local models as in Lemma 
\ref{modellem}, with a fixed linear projection 
$X\mapsto X_{\got{l}}$ from $\got{t}$ onto 
$\got{l}$. For every $\xi\in\got{l}^*$, an {\em admissible lift  
of $\xi$} is a smooth $T$\--invariant vector field 
$L_{\xi}$ on $M_{\scriptop{reg}}$ such that for each local model 
as in Lemma \ref{modellem} there is a smooth $T$\--invariant 
$\got{l}$\--valued function $\alpha _{\xi}$ on $U$, such that 
$L_{\xi}(x)=L^{\Phi}_{\xi}(x)+\alpha _{\xi}(x)_M(x)$ for every $x\in U$. 
Here $L^{\Phi}_{\xi}$ is the local model lift introduced in 
Definition \ref{localmodelliftdef}. 

If we are at an orbit type $\Sigma$ with stabilizer group $H$, 
and $\xi$ is equal to zero on the Lie algebra $\got{h}\subset\got{l}$ 
of $H$, then $L_{\xi}$ has a unique smooth $T$\--invariant 
extension to an open 
neighborhood of $\Sigma$ in $M$, which will also be denoted by $L_{\xi}$. 
In particular, if $\zeta\in N:=(\got{l}/\got{t}_{\scriptop{h}})^*$, 
the space of linear forms on $\got{l}$ which vanish 
on $\got{t}_{\scriptop{h}}$, then $L_{\zeta}$ is a smooth 
$T$\--invariant vector field 
on the whole manifold $M$. 

An {\em admissible connection} for the principal $T$\--bundle 
$\pi :M_{\scriptop{reg}}\to (M/T)_{\scriptop{reg}}$ is a 
linear mapping 
$\xi\mapsto L_{\xi}$ from $\got{l}^*$ to the space of 
smooth vector fields on $M_{\scriptop{reg}}$, 
such that, for each $\xi\in\got{l}^*$, $L_{\xi}$ 
is an admissible lift of $\xi$. Because we work with a
fixed action of the torus $T$, we will just write 
``admissible connection'' in the sequel. 

In the literature, the term ``admissible connection'' has been 
used in various different frameworks and with correspondingly 
different meanings. Our usage of the term ``admissible connection'' 
continues this.  
\label{liftdef}
\end{definition}
\begin{lemma}
There exist admissible connections $\xi\mapsto L_{\xi}$. 
For each admissible connection $\xi\mapsto L_{\xi}$, we have 
\begin{equation}
\sigma (L_{\xi},\, X_M)=\xi (X_{\got{l}}),
\quad \xi\in\got{l}^*,\quad X\in\got{t}. 
\label{stronglift}
\end{equation}  
\label{liftexistencelem}
\end{lemma}
\begin{proof}
If we piece the local model lifts $L^{\Phi}_{\xi}$, introduced in 
Definition \ref{localmodelliftdef}, together by means of 
a partition of unity consisting of smooth 
$T$\--invariant functions with supports in 
the local model neighborhoods $U$, then it follows 
from Lemma \ref{LPhilem} that the resulting connection 
is admissible. 

In the local model of Lemma \ref{modellem}, 
where we use the $(\theta ,\,\rho )$\--coordinates 
in $\C ^m$ as in (\ref{rhodef}), (\ref{Aiota}), we 
have $L^{\Phi}_{\xi}=(0,\,\delta\lambda _{\xi},\, 
(0,\,\delta\rho _{\xi}))$ 
with $(\delta\lambda _{\xi},\,\delta\rho _{\xi})=A^{-1}(\xi )$. 
Furthermore $X\cdot (k,\,\lambda ,\, z)=
(X_{\got{k}},\, 0,\, (\iota (X_{\got{h}})/\op{i},\, 0))$. 
It follows that 
\[
\sigma (L^{\Phi}_{\xi},\, X_M)=
\delta\lambda _{\xi}((X_{\got{k}})_{\got{l}})+\sum_{j=1}^m\, 
\iota (X_{\got{h}})^j\, 
(\delta\rho _{\xi})_j/\op{i}=\xi (X_{\got{l}}).
\]
Here we have used (\ref{gotkgotl}) with $Y$ replaced by $X$. 
On the other hand 
\[
\sigma _x(\alpha _{\xi}(x)_M,\, X_M(x))
=\sigma ^{\got{t}}(\alpha _{\xi}(x),\, X)=0
\] 
for every $X\in\got{t}$ if $\alpha_{\xi}(x)\in\got{l}
:=\op{ker}\sigma ^{\got{t}}$, 
and (\ref{stronglift}) now follows from Definition \ref{liftdef}. 
\end{proof}
The equation (\ref{stronglift}) improves upon 
(\ref{lift}) if $\got{l}$ is a proper linear subspace 
of $\got{t}$, that is, if $\sigma ^{\got{t}}\neq 0$. If 
the principal orbits are Lagrangian 
submanifolds of $M$, then $\got{l}=\got{t}$ and 
(\ref{stronglift}) is the same as (\ref{lift}).  

\subsection{Special admissible connections}
\label{nicess}
Recall the Hamiltonian torus $T_{\scriptop{h}}$, the unique maximal stabilizer 
subgroup $T_{\scriptop{h}}$ of $T$ as in Remark \ref{T1rem} and 
Lemma \ref{stabilizerlem}, with Lie algebra $\got{t}_{\scriptop{h}}\subset 
\got{l}$. In our quest for nice admissible lifts, we will use 
a decomposition of $T$ into the subtorus $T_{\scriptop{h}}$ and a 
complementary subtorus $T_{\scriptop{f}}$, as in Lemma \ref{tbaselem} 
with $U=T_{\scriptop{h}}$. Note that the torus $T_{\scriptop{f}}$ acts freely 
on $M$, because if $x\in M$, then $T_x\subset T_{\scriptop{h}}$, hence 
$T_x\cap T_{\scriptop{f}}\subset T_{\scriptop{h}}
\cap T_{\scriptop{f}}=\{ 1\}$. This explains our choice of the 
subscript f in $T_{\scriptop{f}}$. 
Note also that the choice of a complementary subtorus $T_{\scriptop{f}}$ 
to $T_{\scriptop{h}}$ is far from unique if $\{ 1\}\neq T_{\scriptop{h}}
\neq T$. See Remark \ref{tbaserem}. We will refer to 
$T_{\scriptop{f}}$ as {\em a freely acting complementary 
torus to the Hamiltonian torus} $T_{\scriptop{h}}$.  

If $\got{t}_{\scriptop{f}}$ denotes the Lie algebra of 
$T_{\scriptop{f}}$, then 
we have a corresponding direct sum decomposition 
$\got{t}=\got{t}_{\scriptop{h}}\oplus\got{t}_{\scriptop{f}}$ 
of Lie algebras. 
Each linear form on ${\got{t}_{\scriptop{h}}}^*$ has a unique extension 
to a linear form on $\got{l}$ which is equal to zero 
on $\got{t}_{\scriptop{f}}$. This leads to 
an isomorphism of ${\got{t}_{\scriptop{h}}}^*$ 
with the linear subspace 
$(\got{l}/\got{l}\cap\got{t}_{\scriptop{f}})^*$ of 
$\got{l}^*$. This isomorphism depends on the choice of the complementary 
freely acting torus $T_{\scriptop{f}}$ to the Hamiltonian 
torus in $T$. 
Note that the direct sum decomposition 
$\got{l}=\got{t}_{\scriptop{h}}\oplus 
(\got{l}\cap\got{t}_{\scriptop{f}})$ implies 
the direct sum decomposition 
\begin{equation}
\got{l}^*=(\got{l}/\got{l}\cap\got{t}_{\scriptop{f}})^*\oplus 
(\got{l}/\got{t}_{\scriptop{h}})^*.
\label{ldec}
\end{equation}
Let
\begin{equation}
\mu :M\to \Delta\subset 
(\got{l}/\got{l}\cap\got{t}_{\scriptop{f}})^*
\simeq {\got{t}_{\scriptop{h}}}^*
\label{pi^1}
\end{equation}
denote the projection 
$\pi :M\to M/T$, followed by the projection from 
$M/T\simeq\Delta\times (N/P)$ onto the first factor. 
Here we use the isomorphism $\Phi _p:\Delta\times (N/P)\to M/T$ 
of Proposition \ref{M/Tprop}, with $N=(\got{l}/\got{t}_{\scriptop{h}})^*$ 
and $C=(\got{l}/\got{l}\cap\got{t}_{\scriptop{f}})^*
\simeq {\got{t}_{\scriptop{h}}}^*$. Note that $\mu :M\to 
{\got{t}_{\scriptop{h}}}^*$ is a momentum mapping for the 
Hamiltonian $T_{\scriptop{h}}$\--action on $M$ as in  
Corollary \ref{Tstabcor}. 

With these notations, Proposition \ref{liftprop} below 
yields the existence of an admissible connection for which 
both the Lie brackets and the symplectic products of 
the $L_{\xi}$ take an extremely simple form. 
(We are tempted to call such a connection a 
``minimal admissible connection'', see also 
Subsection \ref{minimalcouplingss}, but we 
do not have a proposal for a functional which 
is minimized exactly by the connections in 
Proposition \ref{liftprop}.) 
In Remark \ref{curvrem} and Remark \ref{integralcrem} 
we discuss the topological meaning of the antisymmetric 
bilinear form $c:N\times N\to\got{l}$. From these remarks 
it follows that $c$ is unique. The freedom in the choice of the  
admissible connection in Proposition \ref{liftprop}
will be described in Lemma \ref{liftfreedomlem}.  

\begin{proposition}
There exists an admissible connection 
$\got{l}^*\ni\xi\mapsto L_{\xi}$ as in Definition 
\ref{liftdef}, and an antisymmetric bilinear mapping 
$c:N\times N\to\got{l}$, with the following properties. 
\begin{itemize}
\item[i)] 
$\left[ L_{\eta},\, L_{\eta'}\right] =0$ 
~for all $\eta ,\,\eta '\in C$, 
\item[ii)] 
$\left[L_{\eta},\, L_{\zeta}\right]=0$ ~for all
$\eta\in C$ and 
$\zeta\in N$, 
\item[iii)] 
$\left[ L_{\zeta},\, L_{\zeta'}\right] =c(\zeta ,\,\zeta ')_M$  
~for all $\zeta ,\,\zeta '\in N$, 
\item[iv)] 
$\sigma( L_{\eta},\, L_{\eta '})=0$ 
~for all $\eta ,\,\eta '\in C$, 
\item[v)] 
$\sigma(L_{\eta},\, L_{\zeta})=0$ ~for all
$\eta\in C$ and 
$\zeta\in N$, and finally 
\item[vi)] 
$\sigma _x(L_{\zeta}(x),\, L_{\zeta '}(x))
=-\mu (x)(c_{\scriptop{h}}(\zeta ,\zeta '))$ 
~for all $\zeta ,\,\zeta '\in N$ and 
$x\in M$. 
Here $c_{\scriptop{h}}(\zeta ,\zeta ')$ denotes
the $\got{t}_{\scriptop{h}}$\--component of $c(\zeta ,\,\zeta ')$ 
in the direct sum decomposition $\got{l}=\got{t}_{\scriptop{h}}
\oplus (\got{l}\cap\got{t}_{\scriptop{f}})$. 
\end{itemize}
The antisymmetric bilinear mapping $c:N\times N\to\got{l}$ 
in part iii) satisfies the relation 
\begin{equation}
\zeta (c(\zeta ',\,\zeta ''))
+\zeta '(c(\zeta '',\,\zeta ))
+\zeta ''(c(\zeta ,\,\zeta '))=0
\label{ft}
\end{equation}
for every $\zeta ,\,\zeta ',\,\zeta ''\in 
N$. Note that 
$\zeta$ is a linear 
form on $\got{l}$ which vanishes on $\got{t}_{\scriptop{h}}$, 
and therefore 
$\zeta (c(\zeta ',\,\zeta ''))$ is a real number  
which only depends on the projection 
of $c(\zeta ,\,\zeta ')$ to 
$\got{l}/\got{t}_{\scriptop{h}}$. 
\label{liftprop}
\end{proposition}
\begin{proof}
We start with an arbitrary admissible connection 
$\got{l}^*\ni\xi\mapsto L_{\xi}$, which exists 
according to Lemma \ref{liftexistencelem}, and first 
simplify the Lie brackets. 

We use the isomorphism 
$\Phi _p:\Delta\times (N/P)\to M/T$ 
of Proposition \ref{M/Tprop}, in order to identify 
$M/T$ with  $\Delta\times (N/P)=(\Delta\times N)/P\subset 
\got{l}^*/P$. In view of Lemma \ref{basicformlem}, 
the smooth basic $k$\--forms on $M$ satisfy 
$\omega =\pi^*\nu$ on $M_{\scriptop{reg}}$ for 
uniquely determined 
smooth $k$\--forms $\nu$ on $(\Delta _{\scriptop{reg}}\times N)/P$, 
such that $\nu$ extends to a smooth $k$\--form on $\got{l}^*/P$. 
This leads to an identification of the space of 
all smooth basic $k$\--forms 
on $M$ with the space of all restrictions to $(\Delta\times N)/P$ 
of smooth $k$\--forms on $\got{l}^*/P$. If we view 
$\xi\in\got{l}^*$ as a constant vector field on $\got{l}^*/P$, 
then the fact that $\op{T}\!\pi$ maps $L_{\xi}$ to $\xi$ 
implies that 
$
(\pi ^*\nu )(L_{\xi ^1},\,\ldots ,\, L_{\xi ^k})=
\nu (\xi ^1,\,\ldots ,\, \xi ^k). 
$
Because $\pi$ intertwines the flow of $L_{\xi}$ 
with the flow of the constant vector field $\xi$, 
we also have the identity $\op{L}_{L_{\xi}}(\pi ^*\nu )
=\pi ^*(\op{L}_{\xi}\nu )$ for the Lie derivatives. 
In particular the differentiation of $T$\--invariant smooth functions 
on $M$ in the direction of the vector field $L_{\xi}$ 
corresponds to the differentiation $\partial _{\xi}$ 
of smooth functions on $M/T$ in the direction of the 
constant vector field $\xi$. 

Let $X_i$, $1\leq i\leq d_{\got{l}}:=\op{dim}\got{l}$, 
be a basis of $\got{l}$. We will write 
$\alpha =\alpha ^i\, X_i$, in which 
the real numbers $\alpha ^i$ are the coordinates 
of $\alpha \in\got{l}$ with respect to 
this basis, and we use Einstein's summation 
convention when summing over indices which run from 
$1$ up to $d_{\got{l}}$. 

The local model lifts 
$L^{\Phi}_{\xi}$ and $L^{\Phi}_{\xi '}$, 
introduced in Definition \ref{localmodelliftdef}, commute 
because they are constant vector fields on $K\times E_0$ in the 
local model of Lemma \ref{modellem}, where we use 
the $(\theta ,\,\rho )$\--coordinates 
for $z$. If $L_{\xi}=L^{\Phi}_{\xi}+\alpha ^i_{\xi}\, (X_i)_M$ as in 
Definition \ref{liftdef}, then the fact that the vector fields 
$L^{\Phi}_{\xi}$ and $L^{\Phi}_{\xi '}$ commute as well 
as the vector fields $\alpha ^i_{\xi}\, (X_i)_M$ and 
$\alpha ^j_{\xi '}\, (X_j)_M$, implies that 
$[L_{\xi},\, L_{\xi '}]=\beta ^i_{\xi ,\,\xi '}\, (X_i)_M$, 
in which the uniquely determined $T$\--invariant functions 
$\beta ^i_{\xi ,\,\xi '}$ on $U=\Phi (K\times E_0)$ are given by 
\[
\beta ^i_{\xi ,\,\xi '}=L_{\xi}\,\alpha ^i_{\xi '}
-L_{\xi '}\,\alpha ^i_{\xi}
=\partial _{\xi}\,\alpha ^i_{\xi '}
-\partial _{\xi '}\,\alpha ^i_{\xi}.
\]
That is, $\beta =\op{d}\!\alpha$ on $U$, if we define 
the basic $\got{l}$\--valued one\--form $\alpha$ and 
two\--form $\beta$ by $\alpha (\xi )=\alpha ^i_{\xi}\, X_i$ 
and $\beta (\xi ,\,\xi ')=\beta ^i_{\xi ,\,\xi '}\, X_i$, 
respectively. Because $\beta$ is locally exact, it follows that 
the globally defined smooth basic two\--form $\beta$ is closed. 

Any other connection $\got{l}^*\ni\xi\to \widetilde{L}_{\xi}$ is 
admissible, if and only if 
\begin{equation}
\widetilde{L}_{\xi}=L_{\xi}+\widetilde{\alpha} ^i_{\xi}\, (X_i)_M, 
\label{tildeU}
\end{equation} 
in which $\widetilde{\alpha}(\xi ):=\widetilde{\alpha} ^i_{\xi}\, X_i$ defines 
a smooth basic $\got{l}^*$\--valued one\--form 
$\widetilde{\alpha}$ on $M$. 
With the same reasoning as above we obtain that 
$[\widetilde{L}_{\xi},\,\widetilde{L}_{\xi '}]
=\widetilde{\beta}^i_{\xi ,\,\xi '}\, (X_i)_M$, 
in which 
$\widetilde{\beta}(\xi ,\,\xi ')
:=\widetilde{\beta}^i_{\xi ,\,\xi '}\, X_i$ 
defines a smooth basic $\got{l}^*$\--valued two\--form 
$\widetilde{\beta}$, such that $\widetilde{\beta}=\beta 
+\op{d}\!\widetilde{\alpha}$. 

According to Corollary \ref{cohomcor}, the de Rham cohomology class of 
$\beta$ contains a unique $c=\widetilde{\beta}$, 
such that $c(\xi ,\,\xi ')=0$ when 
$\xi\in C$ or $\xi '\in C$, 
and for $\xi ,\,\xi '\in N$ the $\got{l}$\--valued 
function $c(\xi ,\,\xi ')$ 
is a constant, equal to 
the average of $\beta (\xi,\,\xi ')$ over any $N/P$\--orbit 
in $M/T$. This leads to the desired properties of the Lie brackets,   
where the uniqueness of $c$ follows from the injectivity of 
the mapping $\Lambda ^2N^*\to\op{H}^2(N/P,\,\R )$ 
in Corollary \ref{cohomcor}.  

We now turn to the symplectic inner products. 
Let $\xi,\,\xi '\in\got{l}^*$. 
It follows from Definition \ref{localmodelliftdef} 
that $\sigma (L^{\Phi}_{\xi},\, L^{\Phi}_{\xi'})=0$. 
In view of Definition \ref{liftdef} and formula 
(\ref{stronglift}) we conclude that 
$\sigma (L_{\xi},\, L_{\xi '})$ 
is a smooth function on $M$, which moreover is $T$\--invariant. 
Therefore $(\xi,\,\xi ')\mapsto \sigma (L_{\xi},\, L_{\xi '})$ 
defines a smooth basic two\--form $s$ on $M$. If in 
(\ref{tildeU}) we take $\widetilde{\alpha} ^i_{\xi}
=\partial _{\xi}\,\varphi ^i$ for 
$\varphi ^i\in\op{C}^{\infty}(M/T)$, that is, 
$\widetilde{\alpha} =\op{d}\!\varphi$, 
then the Lie brackets do not change, but 
\[
\widetilde{s}(\xi ,\,\xi ') 
:=\sigma (\widetilde{L}_{\xi},\,\widetilde{L}_{\xi '})
=s(\xi ,\,\xi ')+\partial _{\xi '}\, (\varphi (\xi ))
-\partial _{\xi}\, (\varphi (\xi ')), 
\]
where we have used that $\sigma ((X_i)_M,\, (X_j)_M)=
\sigma ^{\got{t}}(X_i,\, X_j)=0$, because $X_i,\, X_j\in\got{l}$. 
This means that $\widetilde{s}=s-\op{d}\!\varphi$. 

In order investigate the exterior derivative of $s$, 
we recall the identity
\begin{eqnarray}
(\op{d}\!\omega)(u,\, v,\, w)&=&
\partial _u(\omega(v,\, w))
+\partial _v(\omega(w,\, u))
+\partial _w(\omega(u,\, v))\nonumber\\
&&+\omega(u,\, [v,\, w])+\omega(v,\, [w,\, u])+\omega(w,\, [u,\, v]),
\label{domega}
\end{eqnarray}
which holds for any smooth two\--form $\omega$ and 
smooth vector fields $u,\, v,\, w$. It follows that 
\begin{eqnarray*}
(\op{d}\! s)(\xi ,\,\xi ',\,\xi '')&=&
\partial _{\xi}\, s(\xi ',\,\xi '')
+\partial _{\xi '}\, s(\xi '',\,\xi )
+\partial _{\xi ''}\, s(\xi ,\,\xi ')\\
&=&L_{\xi}(\sigma(L_{\xi '},\, L_{\xi ''}))
+L_{\xi '}(\sigma(L_{\xi ''},\, L_{\xi}))
+L_{\xi ''}(\sigma(L_{\xi},\, L_{\xi '}))\\
&=&-\sigma(L_{\xi},\, [L_{\xi '},\, L_{\xi ''}])
-\sigma(L_{\xi '},\, [L_{\xi ''},\, L_{\xi}])
-\sigma(L_{\xi ''},\, [L_{\xi},\, L_{\xi '}])\\
&=&-\xi (c(\xi ',\,\xi ''))-\xi '(c(\xi '',\,\xi ))-\xi ''(c(\xi ,\,\xi ')),
\end{eqnarray*}
where in the third identity we have used that $\op{d}\!\sigma = 0$,  
and in the last identity we have inserted 
$[L_{\xi},\, L_{\xi '}]=c(\xi ,\,\xi ')^i\, (X_i)_M$ 
and (\ref{lift}). 

This shows that $\op{d}\! s$ is constant. In the notation of 
Corollary \ref{cohomcor}, we have that 
$\op{d} (\iota_p^*s)=\iota_p^*(\op{d}\! s)$ is constant, 
and cohomologically equal 
to zero, which in view of the first part of the last statement 
in Corollary \ref{cohomcor} implies that 
$\iota_p^*(\op{d}\! s)=0$. That is, 
$(\op{d}\! s)(\xi ,\,\xi ',\,\xi '')=0$ 
when $\xi ,\,\xi ',\,\xi ''\in N$, which in turn is equivalent to 
(\ref{ft}). On the other hand it follows from 
the already proved statements about the Lie brackets that 
$c(\xi ,\,\xi ')=0$ 
if $\xi\in C$ or $\xi '\in C$, and hence 
$(\op{d}\! s)(\xi ,\,\xi ',\,\xi '')=0$ unless 
one of the vectors $\xi ,\,\xi ',\,\xi ''$ 
belongs to $C$ 
and the other two belong to $N$. 
Moreover, if $\xi\in C$ and $\xi ',\,\xi ''\in N$, 
then we obtain that 
$(\op{d}\! s)(\xi ,\,\xi ',\,\xi '')=\, -\xi (c(\xi ',\,\xi ''))$. 

In other words, the smooth basic two\--form 
$S:=s+\mu\, c_{\scriptop{h}}$ is closed. Here $\mu$ is viewed as a 
${\got{t}_{\scriptop{h}}}^*$\--valued $T$\--invariant function on $M$, 
and the pairing with the $\got{t}_{\scriptop{h}}$\--valued 
antisymmetric bilinear form $c_{\scriptop{h}}$ 
yields a smooth basic two\--form $\mu\, c_{\scriptop{h}}$ 
on $M$.  
According to Corollary \ref{cohomcor}, the smooth basic 
one\--form $\varphi$ can be now chosen such that if 
$\widetilde{S}=S-\op{d}\!\varphi$, then 
$\widetilde{S}(\xi ,\,\xi ')=0$ when $\xi\in C$ or $\xi '\in C$, 
and for $\xi ,\,\xi '\in N$ the function 
$\widetilde{S}(\xi ,\,\xi ')$ is constant. 

We finally observe that if $\alpha :\xi\mapsto\alpha _{\xi}$ 
is a linear mapping from $\got{l}^*$ to $\got{l}$, which 
is viewed as a constant, hence closed 
one\--form on the $\got{l}^*$\--parallel space $M/T$, 
then the Lie brackets of the $L_{\xi}$'s do not change if 
we replace $L_{\xi}$ by $L_{\xi}+(\alpha _{\xi})_M$. 
However, $\sigma (L_{\xi},\, L_{\xi '})$ then gets replaced 
by $\sigma(L_{\xi},L_{\xi'})+\xi(\alpha_{\xi'})-\xi'(\alpha_{\xi})$. 
Because any antisymmetric bilinear form on $\got{l}^*$ 
is of the form 
$(\xi,\,\xi')\mapsto \xi(\alpha_{\xi'})-\xi'(\alpha_{\xi})$, 
for a suitable linear mapping $\alpha:\got{l}^*\to\got{l}$, 
we can arrange that $\widetilde{S}=0$, which leads to 
vi) in Proposition \ref{liftprop}. 
\end{proof}

\begin{remark}
Because the left hand side of (\ref{ft}) is antisymmetric in 
$\zeta ,\,\zeta ',\,\zeta ''$, it is automatically equal to zero when 
$\op{dim}N=\op{dim}\got{l}-\op{dim}\got{t}_{\scriptop{h}}\leq 2$. 
However, when $\op{dim}N\geq 3$, then the equations (\ref{ft}) 
impose nontrivial conditions on the 
$\got{l}$\--valued two\--form $c$ on $N$. 
\label{ftrem}
\end{remark}

\begin{remark}
For every $x\in M_{\scriptop{reg}}$, let $H_x$ denote the linear 
span in $\op{T}_x\! M$ of the vectors $L_{\xi}(x)$, $\xi\in\got{l}^*$. 
Then the $H_x$, $x\in M_{\scriptop{reg}}$, define a $T$\--invariant 
infinitesimal connection of the principal $T$\--bundle 
$M_{\scriptop{reg}}$ over $(M/T)_{\scriptop{reg}}\simeq 
\Delta ^{\scriptop{int}}\times (N/P)$. Here 
$\Delta ^{\scriptop{int}}$ denotes the interior 
of the Delzant polytope $\Delta$. 
Any connection 
of this principal $T$\--bundle has 
a {\em curvature form} which is a smooth $\got{t}$\--valued 
two\--form on $M_{\scriptop{reg}}/T$. The cohomology class of 
the curvature form is an element of 
$\op{H}^2(M_{\scriptop{reg}}/T,\,\got{t})$, which is independent of the 
choice of the connection. The action of 
$N$ on $M/T$ leaves $M_{\scriptop{reg}}/T\simeq (M/T)_{\scriptop{reg}}$ 
invariant, with orbits isomorphic to the torus $N/P$, and the 
pull\--back to the $N$\--orbits defines an isomorphism 
from $\op{H}^2(M_{\scriptop{reg}}/T,\,\got{t})$ onto 
$\op{H}^2(N/P,\,\got{t})$, which in turn is identified 
with $(\Lambda ^2N^*)\otimes\got{t}$ as in Corollary 
\ref{cohomcor}. 

The proof of Proposition \ref{liftprop} 
shows that the element 
$c\in (\Lambda ^2N^*)\otimes\got{l}\subset(\Lambda ^2N^*)\otimes\got{t}$ 
is equal to the negative of the pull\--back to an $N$\--orbit of the 
cohomology class of the curvature form. This proves in particular that 
the antisymmetric bilinear mapping $c:N\times N\to\got{l}$ 
in Proposition \ref{liftprop} is independent of the choice 
of the freely acting complementary torus $T_{\scriptop{f}}$ to 
the Hamiltonian torus in $T$. More precisely, if 
$(\widetilde{M},\,\widetilde{\sigma},\, T)$ is another 
symplectic manifold with an effective symplectic $T$\--action with 
coisotropic principal orbits, then Proposition \ref{liftprop} 
with $(M,\,\sigma ,\, T)$ replaced by 
$(\widetilde{M},\widetilde{\sigma}, T)$ yields 
an antisymmetric bilinear mapping $\widetilde{c}$ 
instead of $c$. If there exists 
a $T$\--equivariant symplectomorphism $\Phi$ from 
$(M,\,\sigma ,\, T)$ onto $(\widetilde{M},\widetilde{\sigma}, T)$
then $\widetilde{c}=c$. 
\label{curvrem}
\end{remark}

\begin{remark}
The retrivializations of the principal $T$\--bundle 
$\pi :M_{\scriptop{reg}}\to M_{\scriptop{reg}}/T$ define a one\--cocycle 
of smooth $\got{t}$\--valued functions on $M_{\scriptop{reg}}/T$, 
of which the
sheaf (= \v{C}ech) cohomology class $\tau$ in 
$\op{H}^1(M_{\scriptop{reg}}/T,\, 
\op{C}^{\infty}(\cdot ,\, T))$ 
classifies the principal $T$\--bundle 
$\pi :M_{\scriptop{reg}}\to M_{\scriptop{reg}}/T$. 
Because the sheaf $\op{C}^{\infty}(\cdot,\,\got{t})$ is fine, 
the short exact sequence 
\[
0\to T_{\Z}\to\op{C}^{\infty}(\cdot ,\,\got{t})
\stackrel{\scriptop{exp}}{\to}
\op{C}^{\infty}(\cdot ,\, T)\to 1
\] 
induces an isomorphism $\delta :\op{H}^1(M_{\scriptop{reg}}/T,\, 
\op{C}^{\infty}(\cdot ,\, T))\to\op{H}^2(M_{\scriptop{reg}}/T,\, T_{\Z})$. 
Here exp denotes the exponential mapping $\got{t}\to T$.  
The cohomology class 
$\delta (\tau )\in\op{H}^2(M_{\scriptop{reg}}/T,\, T_{\Z})$ 
is called the {\em Chern class} of the principal $T$\--bundle 
$\pi :M_{\scriptop{reg}}\to M_{\scriptop{reg}}/T$.  
It is a general fact, see for instance the 
arguments in \cite[Sec. 15.3]{heatkernel}, that the image 
of $\delta (\gamma)$ in $\op{H}^2(M_{\scriptop{reg}}/T,\,\got{t})$ 
under the coefficient homomorphism 
$\op{H}^2(M_{\scriptop{reg}}/T,\, T_{\Z})
\to\op{H}^2(M_{\scriptop{reg}}/T,\,\got{t})$ 
is equal to the negative of the cohomology 
class of the curvature form of any connection in the principal 
$T$\--bundle. In view of Remark \ref{curvrem}, we therefore 
conclude that $c$ represents the Chern class of the  
principal $T$\--bundle 
$\pi :M_{\scriptop{reg}}\to M_{\scriptop{reg}}/T$. 

In view of the canonical isomorphism between sheaf cohomology 
and singular cohomology, this implies that the integral of 
$c$ over every two\--cycle in $(M/T)_{\scriptop{reg}}$ 
belongs to $T_{\Z}$. If $\zeta ,\,\zeta '\in P$, then 
for every $p\in (M/T)_{\scriptop{reg}}$ the mapping 
\[
\iota _{\zeta ,\,\zeta'}:(t,\, t')\mapsto p+(t\,\zeta +t'\,\zeta '):
\R ^2/\Z ^2\to (M/T)_{\scriptop{reg}}
\]
defines a two\--cycle in $(M/T)_{\scriptop{reg}}$, and 
\[
c(\zeta ,\,\zeta ')
=\int _{\R ^2/\Z ^2}\, (\iota _{\zeta ,\,\zeta '})^*c.
\]
It follows that $c(\zeta ,\,\zeta ')\in T_{\Z}$ for 
every $\zeta ,\,\zeta '\in P$.  

In Lemma \ref{TPlem}, this conclusion will be proved by means of 
a group theoretical consideration.  
Other topological interpretations of 
$c$ will be given in 
Proposition \ref{pi1prop} and Proposition \ref{chernprop}. 
\label{integralcrem}
\end{remark}

Let $\op{Lin}(E, F)$ denote the space of all linear mappings 
from a vector space $E$ to a vector space $F$. 
In the following lemma we use that a smooth $T$\--invariant mapping 
$\alpha :M\to\op{Lin}(\got{l}^*,\,\got{l})$ 
corresponds to 
a unique smooth basic $\got{l}$\--valued one\--form on $M$, 
which we also denote by $\alpha$, such that 
$\alpha (x)(\xi )=\alpha _x(v)$ for every 
$v\in\op{T}_x\! M$ such that $\widehat{\sigma}_x(v)=\xi$. 
In view of Lemma \ref{basicformlem}, $\alpha$ can also be viewed 
as the restriction to $M/T$ of a smooth $\got{l}$\--valued 
one\--form on $\got{l}^*/P$, if we identify the 
$\got{l}^*$\--parallel space $M/T$ with a subset 
of $\got{l}^*/P$ as in Proposition \ref{M/Tprop}.  
\begin{lemma}
Let $\got{l}^*\ni\xi\mapsto L_{\xi}$ be a connection 
as in Proposition \ref{liftprop}. Then 
$\got{l}^*\ni\xi\mapsto\widetilde{L}_{\xi}$ 
is an admissible connection, 
if and only if there exists a smooth $T$\--invariant 
mapping 
$\alpha :x\mapsto (\xi\mapsto\alpha _{\xi}(x))$ from 
$M$ to $\op{Lin}(\got{l}^*,\,\got{l})$, such that 
$\widetilde{L}_{\xi}(x)=
L_{\xi}(x)+\alpha_{\xi}(x)_M(x)$ for every 
$x\in M$ and $\xi\in\got{l}^*$. 
Proposition \ref{liftprop} holds with 
$L$ replaced by $\widetilde{L}$, if and only if 
$\alpha$ is closed when considered as a smooth 
basic $\got{l}$\--valued 
one\--form on $M$, and moreover $\alpha$ is 
symmetric in the sense that 
\begin{equation}
\xi(\alpha_{\xi'}(x))-\xi'(\alpha_{\xi}(x)) = 0
\label{alphasym}
\end{equation}
for all $\xi,\,\xi '\in\got{l}^*$ and all $x\in M$. 
\label{liftfreedomlem}
\end{lemma}
\begin{proof}
The first statement follows from Definition 
\ref{liftdef}, the definition of admissible connections. 
It follows from the proof of Proposition \ref{liftprop} 
that $[\widetilde{L}_{\xi },\,\widetilde{L}_{\xi'}] 
\equiv [L_{\xi},\, L_{\xi'}]$ if and only if 
$\alpha$ is closed. In view of the uniqueness of 
$c$, see Remark \ref{curvrem}, we have iv), v), vi) in 
Proposition \ref{liftprop} with $L$ replaced by 
$\widetilde{L}$, if and only if  
$\sigma (\widetilde{L}_{\xi},\,\widetilde{L}_{\xi '})
\equiv\sigma (L_{\xi},\, L_{\xi '})$, which is equivalent to   
(\ref{alphasym}). 
\end{proof}

\section{Delzant submanifolds}
\label{delzantsec}
Recall the especially nice admissible connection 
introduced in Proposition \ref{liftprop}, the construction 
of which is based on the identification in Proposition \ref{M/Tprop} 
of the orbit space $M/T$ 
with the $\got{l}^*$\--parallel space $\Delta\times (N/P)$. 
Proposition \ref{embedprop} below implies that the vector fields 
$Y_M$, $Y\in\got{t}_{\scriptop{h}}$, and 
$L_{\eta}$, $\eta\in C$, are tangent to the 
fibers of a fibration of $M$ by Delzant submanifolds. 
From this section on, the word {\em fibration} is short for a 
locally trivial smooth fiber bundle. The remainder of 
this section is devoted to the proof and further precision of 
Proposition \ref{embedprop}. For any subset $Y$ of a set $X$, 
the {\em inclusion mapping} $\iota _Y$ is the identity on $Y$, 
viewed as a mapping from $Y$ to $X$.  
\begin{proposition}
Let $\got{l}^*\ni\xi\mapsto L_{\xi}$ be an admissible connection 
as in Proposition \ref{liftprop}. 
Then there is a unique smooth $T$\--invariant 
distribution $D$ on $M$ such that,  
for every $x\in M_{\scriptop{reg}}$, $D_x$ is equal to the linear span 
in $\op{T}_x\! M$ of the vectors $Y_M(x)$ with 
$Y\in\got{t}_{\scriptop{h}}$  
and $L_{\eta}(x)$, $\eta\in C:=(\got{l}/\got{l}\cap 
\got{t}_{\scriptop{f}})^*\simeq {\got{t}_{\scriptop{h}}}^*$. 

The distribution $D$ is integrable. Each integral manifold 
manifold $I$ of $D$ is invariant under the action of the 
Hamiltonian torus $T_{\scriptop{h}}$, and 
$(I,\, {\iota _I}^*\,\sigma ,\, T_{\scriptop{h}})$ 
is a Delzant manifold with the Delzant polytope $\Delta$ 
introduced in Proposition \ref{M/Tprop}. 
Here $\iota _I:I\to M$ is the inclusion mapping from $I$ 
into $M$. 
The integral manifolds of $D$ form 
a smooth fibration of $M$ into 
Delzant submanifolds with Delzant polytope $\Delta$. 
\label{embedprop}
\end{proposition}
\begin{proof}
This follows from Lemma \ref{embedlem} below, which in 
turn uses Lemma \ref{pi2lem}. 
\end{proof}

\begin{lemma}
Let $\pi _{N/P}:M\to N/P$ be the mapping which is equal to 
$\pi :M\to M/T$, followed by the inverse $M/T\to\Delta\times (N/P)$ 
of the isomorphism $\Phi _p$ in Proposition \ref{M/Tprop}, iii), 
followed by the projection $\Delta\times (N/P)\to N/P$ onto the 
second factor. 

Then $\pi _{N/P}:M\to N/P$ defines a smooth 
fibration of $M$ 
over the torus $N/P$. Each fiber $F$ of $\pi _{N/P}:M\to N/P$ 
is a connected compact $T$\--invariant 
smooth submanifold of $M$. For each fiber $F$ of 
$\pi _{N/P}:M\to N/P$, 
$F\cap M_{\scriptop{reg}}$ is dense in $F$. 
\label{pi2lem}
\end{lemma}
\begin{proof}
Because $\pi :M_{\scriptop{reg}}\to (M/T)_{\scriptop{reg}}$ and 
the projection from  $(M/T)_{\scriptop{reg}}\simeq 
\Delta _{\scriptop{reg}}\times (N/P)$ onto 
$N/P$ are smooth fibrations with connected 
fibers, it follows that the restriction to $M_{\scriptop{reg}}$ 
of $\pi _{N/P}$ is a smooth fibration with connected fibers.  

In the local model of Lemma \ref{modellem}, 
the mapping $\pi _{N/P}$ corresponds to the 
mapping 
\[
(k,\,\lambda ,\, z)\mapsto\lambda _{(\got{l}/\got{t}_{\scriptop{h}})^*}+P
\in N/P,
\]
where we have used the direct sum decomposition (\ref{ldec}). 
This shows that $\pi _{N/P}$ is a smooth submersion. 
Moreover, for each fiber $F$ of $\pi _{N/P}:M\to N/P$, 
$F\cap M_{\scriptop{reg}}$ is dense in $F$, because 
the point $(k,\,\lambda ,\, z)$ is regular if and only if 
$z^j\neq 0$ for every $j$. 
Because the fiber $F\cap M_{\scriptop{reg}}$ of the restriction to 
$M_{\scriptop{reg}}$ of $\pi _{N/P}$ is connected, it follows that 
$F$ is connected. Because $M$ is compact, the submersion 
$\pi _{N/P}$ is proper, and because every proper submersion is 
a fibration, it follows that $\pi _{N/P}$ 
is a fibration. 
\end{proof}
As observed in the beginning 
of Subsection \ref{nicess}, the action on $M$ of the complementary 
torus $T_{\scriptop{f}}$ to $T_{\scriptop{h}}$ is free. 
This exhibits each fiber $F$ of $\pi _{N/P}$ as a 
principal $T_{\scriptop{f}}$\--bundle 
$\pi _{F/T_{\scriptop{f}}}:F\to F/T_{\scriptop{f}}$, 
in which the $T_{\scriptop{f}}$\--orbit space $F/T_{\scriptop{f}}$ 
is a compact, connected smooth manifold, on which we still have the 
action of the Hamiltonian torus $T_{\scriptop{h}}$. 
The following lemma says that there is a symplectic form 
$\sigma _{F/T_{\scriptop{f}}}$ on $F/T_{\scriptop{f}}$ 
such that 
\begin{equation}
(F/T_{\scriptop{f}},\,\sigma _{F/T_{\scriptop{f}}},
\, T_{\scriptop{h}})
\label{F/Tdelzant}
\end{equation}
is a Delzant manifold defined by the Delzant polytope 
$\Delta$, and that the 
fibration $\pi _{F/T_{\scriptop{f}}}:
F\to F/T_{\scriptop{f}}$ is 
trivial, exhibiting $F$ as the Cartesian product 
of the Delzant manifold $F/T_{\scriptop{f}}$ with 
$T_{\scriptop{f}}$. 

\begin{lemma}
There is a unique smooth distribution $D$ on $M$ such that,  
for every $x\in M_{\scriptop{reg}}$, $D_x$ is equal to the linear span 
in $\op{T}_x\! M$ of the vectors $Y_M(x)$, 
$Y\in\got{t}_{\scriptop{h}}$, 
and $L_{\eta}(x)$, $\eta\in C$. 
The distribution $D$ is integrable and $T$\--invariant. 

For every fiber $F$ of the fibration $\pi _{N/P}$ in 
Lemma \ref{pi2lem}, we have $D|_F\subset\op{T}\! F$, 
which implies that $I\subset F$ or $I\cap F=\emptyset$ 
for every integral manifold $I$ of $D$. 
Let $f_0\in F$ and let $I_0$ be the integral manifold 
of $D$ such that $f_0\in I_0$. For each $y\in F/T_{\scriptop{f}}$ 
there is a unique $i(y)\in I_0$ 
such that 
$\pi _{F/T_{\scriptop{f}}}(i(y))=y$. The mapping 
\[
(y,\, t_{\scriptop{f}})\mapsto t_{\scriptop{f}}\cdot i(y): 
(F/T_{\scriptop{f}})\times T_{\scriptop{f}}\to F
\]
is the inverse of a trivialization $\tau$ of the principal 
$T_{\scriptop{f}}$\--fibration 
$\pi _{F/T_{\scriptop{f}}}:F\to F/T_{\scriptop{f}}$. 
The trivialization $\tau$ is $T$\--equivariant, where  
$t\in T$ acts on $(F/T_{\scriptop{f}})\times 
T_{\scriptop{f}}$ by sending 
$(\pi _{F/T_{\scriptop{f}}}(f),\,\widetilde{t}_{\scriptop{f}})$ to 
$(\pi _{F/T_{\scriptop{f}}}(t_{\scriptop{h}}\cdot f),\, t_{\scriptop{f}}\, 
\widetilde{t}_{\scriptop{f}})$,  
if $t=t_{\scriptop{h}}\, t_{\scriptop{f}}$, 
with $t_{\scriptop{h}}\in T_{\scriptop{h}}$ and 
$t_{\scriptop{f}}\in T_{\scriptop{f}}$. 

Finally, there is a unique symplectic form  
$\sigma _{F/T_{\scriptop{f}}}$ 
on $F/T_{\scriptop{f}}$ such that, for any 
integral manifold $I$ of $D$ in $F$,  
\begin{equation}
(\pi _{F/T_{\scriptop{f}}}\circ\iota _I)^*
\,\sigma _{F/T_{\scriptop{f}}}
={\iota _I}^*\sigma ,
\label{piiotasigma}
\end{equation}
if $\iota _I:I\to F$ denotes the inclusion mapping from $I$ into $F$. 
With this symplectic form, (\ref{F/Tdelzant}) is a Delzant manifold 
with Delzant polytope $\Delta$. For each integral manifold $I$ of $D$ 
in $F$, $(I,\, {\iota _I}^*\sigma ,\, T_{\scriptop{h}})$ 
is a Delzant manifold with Delzant polytope $\Delta$, 
and  $\pi _{F/T_{\scriptop{f}}}\circ\iota _I$ is a $T_{\scriptop{h}}$\--equivariant 
symplectomorphism from $(I,\, {\iota _I}^*\sigma ,\, T_{\scriptop{h}})$ 
onto the Delzant manifold (\ref{F/Tdelzant}). 
\label{embedlem}
\end{lemma}
\begin{proof}
In order to investigate the $D_x$ with $x\in M_{\scriptop{reg}}$ 
near a singular point $x_0$, we use a local model as in Lemma \ref{modellem}, 
with the $(\theta ,\,\rho )$\--coordinates in $\C ^m$ 
as in (\ref{rhodef}). Here $H=T_{x_0}$ is a subtorus of the 
Hamiltonian torus $T_{\scriptop{h}}$. Let $K_0$ be a complementary 
subtorus to $H$ in $T_{\scriptop{h}}$. We will take 
$K=K_0\, T_{\scriptop{f}}$ as the complementary subtorus 
to $H$ in $T$. For the Lie algebras we have the corresponding 
direct sum decompositions $\got{t}_{\scriptop{h}}=\got{t}_x\oplus 
\got{k}_0$ and 
$\got{k}=\got{k}_0\oplus\got{t}_{\scriptop{f}}$.  
The span of the infinitesimal actions 
of the elements $Y\in\got{h}$ is equal to 
the span of the vector fields $(0,\, 0,\, \partial/\partial\theta _j)$, 
$1\leq j\leq m$, and the vector fields $(Y,\, 0,\, 0)$ 
with $\got{k}_0$. The $L^{\Phi}_{\eta}$, 
$\eta\in C$, are the linear combinations 
of the vector fields $(0,\, 0,\,\partial/\partial\rho _j)$, 
$1\leq j\leq m$, and the vector fields 
$(0,\,\delta\lambda ,\, 0)$, with constant 
$\delta\lambda\in (\got{l}/(\got{h}
\oplus (\got{l}\cap\got{t}_{\scriptop{f}}))^*$. 
According to Definition \ref{liftdef}, the definition of admissible lifts, 
$L_{\eta}=L^{\Phi}_{\eta}+v_{\eta}$ in which the vector 
field $v_{\eta}$ is smooth on $M$, of the form 
$v_{\eta}(x)=\alpha _{\eta}(x)_M(x)$ for a smooth $T$\--invariant 
$\got{l}$\--valued function $\alpha$ on $M$. We write 
$v_j$ instead of $v_{\eta}$ if $\eta$ is such that 
$L_{\eta}^{\Phi}=(0,\, 0,\,\partial/\partial\rho _j)$. 
The problem is that the vector fields $\partial/\partial\theta ^j$ 
and $\partial/\partial\rho _j$ have a zero and a pole 
at $z^j=p^j+\op{i}\, q^j=0$. 

Now 
\[
\frac{\partial}{\partial\rho _j}=(2\rho _j)^{-1}\, 
(p^j\,\frac{\partial}{\partial p^j}+q^j\,\frac{\partial}{\partial q^j})
\quad\mbox{\rm and}\quad
\frac{\partial}{\partial\theta ^j}=-q^j\,\frac{\partial}{\partial p^j}
+p^j\,\frac{\partial}{\partial q^j}
\]
imply that 
\[
p^j\, L_{\eta ^j}-\frac{q^j}{2\rho _j}\, (Y_j)_M 
=\frac{\partial}{\partial p^j}+p^j\, v_j
\quad\mbox{\rm and}\quad
q^j\, L_{\eta ^j}+\frac{p^k}{2\rho _k}\, (Y_j)_M
=\frac{\partial}{\partial q^j}+q^j\, v_j.
\]
These two vector fields are smooth and converge to 
$\partial/\partial p^j$ and $\partial/\partial q^j$, 
respectively, as $z^j\to 0$. This proves the first statement 
in  the lemma. 
We also obtain for every $x\in M$ that 
$\op{T}_xM=D_x\oplus E_x$, if $E_x$ denotes the linear 
span of the $Z_M(x)$, $Z\in\got{t}_{\scriptop{f}}$, 
and the $L_{\eta}(x)$, $\zeta\in N$. 

In view of (\ref{lift}), conclusion i) in Proposition 
\ref{liftprop}, and the commutativity of the 
infinitesimal action of $\got{t}_{\scriptop{h}}$ on $M$, 
the vector fields $Y_M$ and $L_{\eta}$ 
all commute with each other. This implies that 
on $M_{\scriptop{reg}}$ the distribution $D$ satisfies 
the Frobenius integrability condition. Because 
$M_{\scriptop{reg}}$ is dense in $M$, it follows by 
continuity that $D$ is integrable on $M$. 
Because the vector fields $Y_M$, and $L_{\eta}$ 
are $T$\--invariant, the restriction to $M_{\scriptop{reg}}$ of $D$ 
is $T$\--invariant, and it follows by continuity that 
$D$ is $T$\--invariant. 

For each $x\in M_{\scriptop{reg}}$, the vectors 
$X_M(x)$, $X\in\got{t}$, and $L_{\eta}(x)$, 
$\eta\in C:=(\got{l}/\got{l}\cap 
\got{t}_{\scriptop{f}})^*$, together 
span $\op{T}_x\! F=\op{ker}\op{T}_x\!\pi _{N/P}$, 
hence $D_x\subset\op{T}_x\! F$. Because $M_{\scriptop{reg}}\cap F$ is dense 
in $F$, see Lemma \ref{pi2lem}, it follows by continuity that 
$D|_F\subset\op{T}\! F$. This implies in turn that 
if $I$ is an integral manifold of $D$ in $M$ and $I\cap F\neq\emptyset$, 
then $I\subset F$ and $I$ is an integral manifold of $D|_F$.  

Because for every $x\in F$ the linear subspaces 
$D_x$ and $((\got{t}_{\scriptop{f}})_M)_x$ of 
$\op{T}_x\! F$ have zero intersection 
and their dimensions add up to the dimension of $F$, 
we have that $D_x$ is a complementary linear subspace 
to $((\got{t}_{\scriptop{f}})_M)_x$ in $\op{T}_x\! F$, 
and it follows that 
$D|_F$ defines a $T_{\scriptop{f}}$\--invariant infinitesimal connection 
for the principal $T_{\scriptop{f}}$\--bundle 
$\pi _{F/T_{\scriptop{f}}}:F\to F/T_{\scriptop{f}}$. 

It follows from (\ref{stronglift}) and the conclusion v) 
in Proposition \ref{liftprop}, that for every 
$x\in M_{\scriptop{reg}}$ the complementary linear subspaces 
$D_x$ and $E_x$ of 
$\op{T}_x\! M$ are $\sigma _x$\--orthogonal, 
and by continuity the same conclusion follows for every $x\in M$. 
This implies that, for every $x\in M$, $D_x$ is a symplectic 
vector subspace of $\op{T}_x\! M$, and therefore every 
integral manifold $I$ of $D$ is a symplectic submanifold 
of $(M,\,\sigma )$. 

If $I$ is an integral manifold of $D$, 
then the restriction to $I$ of $\pi _{F/T_{\scriptop{f}}}$ is a covering 
from $I$ onto $F/T_{\scriptop{f}}$. 
Because $\sigma$ is invariant under the action of 
$T_{\scriptop{f}}$, there is a unique two\--form 
$\sigma _{F/T_{\scriptop{f}}}$ 
on $F/T_{\scriptop{f}}$ such that 
(\ref{piiotasigma}) holds, and because 
$\pi _{F/T_{\scriptop{f}}}\circ\iota _I$ is a covering, it follows 
that $\sigma _{F/T_{\scriptop{f}}}$ is a smooth 
symplectic form on $F/T_{\scriptop{f}}$. 

The mapping from $F/T_{\scriptop{f}}$ to $\Delta$ induced 
by (\ref{pi^1}), which we also denote by $\mu$,  
is a momentum mapping for the $T_{\scriptop{h}}$\--action 
on the symplectic manifold $(F/T_{\scriptop{f}},\, 
\sigma _{F/T_{\scriptop{f}}})$. Because for any 
$q\in N/P$ the pre\--image of $\{ q\}$ under the 
projection from $M/T\simeq\Delta\times (N/P)$ onto the 
second factor is equal to $\Delta\times\{ q\}$, and 
$\mu$ forgets the second factor, we have that 
$\mu (F)=\Delta$, and therefore $\mu (F/T_{\scriptop{f}})=\Delta$. 
Because $F$ is compact and connected, see Lemma \ref{pi2lem}, 
the image $F/T_{\scriptop{f}}$ of $F$ under the continuous projection 
$F\to F/T_{\scriptop{f}}$ is also compact and connected. 
The conclusion is that (\ref{F/Tdelzant}) is a Delzant 
manifold defined by the Delzant polytope $\Delta$. 

Because $F/T_{\scriptop{f}}$ is simply connected 
in view of  
Lemma \ref{pi1delzantlem}, $(\pi _{F/T_{\scriptop{f}}})|_I:I
\to F/T_{\scriptop{f}}$ 
is a diffeomorphism. The other statements in the lemma 
now readily follow. 
\end{proof}
  
\begin{lemma}
Every Delzant manifold is simply connected. 
\label{pi1delzantlem}
\end{lemma}
\begin{proof}
Every Delzant manifold can be provided with the 
structure of a toric variety defined by a complete fan, 
cf. Delzant \cite{delzant} and Guillemin \cite[App. 1]{guillemin}, 
and Danilov \cite[Th. 9.1]{danilov} observed that such a 
toric variety is simply connected. The argument is that 
the toric variety has an open cell which is isomorphic to 
$\C ^n$, of which the complement is a complex subvariety 
of complex codimension one. Therefore any loop can be deformed 
into the cell and contracted within the cell to a point. 
\end{proof}
\begin{remark}
The pull\--back to each $T_{\scriptop{f}}$\--orbit of 
the symplectic form $\sigma$ on $M$ 
is given by 
\[
\sigma _x(X_M(x),\, Y_M(x))=\sigma ^{\got{t}}(X,\, Y)
\quad\mbox{\rm  for all}\quad X,\, Y\in\got{t}_{\scriptop{f}}.
\]
Because 
$\got{t}=\got{t}_{\scriptop{h}}\oplus \got{t}_{\scriptop{f}}$ 
and $\got{t}_{\scriptop{h}}\subset\got{l}:=\op{ker}\sigma ^{\got{t}}$, 
we have that this pull\--back is equal to zero if and only 
$\sigma ^{\fop{t}}=0$, that is, the principal $T$\--orbits 
are Lagrangian. In this case the tangent spaces of the 
$T_{\scriptop{f}}$\--orbits in $F$ are the kernels of the 
pull\--back to $F$ of $\sigma$, and the symplectic form 
$\sigma _{F/T_{\scriptop{f}}}$ 
on $F/T_{\scriptop{f}}$ is the reduced form of 
the pull\--back to $F$ of $\sigma$. In other words, 
$(F/T_{\scriptop{f}},\,\sigma _{F/T_{\scriptop{f}}})$ 
is a reduced phase space for the ``momentum mapping'' 
$\pi _{N/P}:M\to N/P$ for the $T_{\scriptop{f}}$\--action, 
where the word momentum mapping is put between parentheses because 
the free $T_{\scriptop{f}}$\--action is not Hamiltonian. 

The $T$\--invariant projection 
$\pi _{N/P}:M\to N/P$ induces a $T_{\scriptop{h}}$\--invariant 
projection $\pi _{N/P}:M/T_{\scriptop{f}}\to N/P$, of which 
the fibers are canonically identified with the $F/T_{\scriptop{f}}$, 
where the $F$ are the fibers of $\pi _{N/P}:M\to N/P$.  
If $\sigma ^{\got{t}}=0$, then the symplectic leaves 
in $M/T_{\scriptop{f}}$ of the Poisson structure on 
$\op{C}^{\infty}(M/T_{\scriptop{f}}) 
=\op{C}^{\infty}(M)^{T_{\scriptop{f}}}$ 
are equal to the fibers $F/T_{\scriptop{f}}$ \
of $\pi _{N/P}:M/T_{\scriptop{f}}\to N/P$, provided with 
the symplectic forms $\sigma _{F/T_{\scriptop{f}}}$. 
It is quite remarkable that the symplectic leaves form a 
fibration, because in general the symplectic leaves of 
a Poisson structure are only immersed submanifolds, 
not necessarily closed. 
\label{Fredrem}
\end{remark}
\section{A global model}
\label{normalformsec}
Let $(M,\,\sigma )$ be our compact connected symplectic 
manifold, together with an effective action of 
the torus $T$ by means of symplectomorphisms 
of $(M,\,\sigma )$, such that some (all) principal orbits 
of the $T$\--action are coisotropic submanifolds of $(M,\,\sigma )$.  

In Subsection \ref{Gss} we will show that the $T$\--action together with 
the infinitesimal action of the vector fields 
$L_{\zeta}$, $\zeta\in N$, introduced in Proposition \ref{liftprop}, 
lead to an action on $M$ of a two\--step nilpotent Lie group $G$, 
where $G$ is explicitly defined in terms of the 
antisymmetric bilinear mapping $c:N\times N\to\got{l}$ 
introduced in Proposition \ref{liftprop}. 
Subsection \ref{Gss} is a sequence of definitions, 
together with some of their immediate consequences.

Recall the fibration of $M$ into Delzant submanifolds  
introduced in Proposition \ref{embedprop}. The action of $G$ on $M$
will be used to 
exhibit this fibration   
as a $G$\--homogeneous bundle over the homogeneous space $G/H$ 
with fiber equal to a Delzant manifold 
defined by the Delzant polytope $\Delta$. 
Here $H$ is a closed Lie subgroup of 
$G$ which is explicitly defined in terms 
of $c$ and the period group $P$ in $N$, 
defined in Lemma \ref{perlem} with $Q=M/T$, 
$V=\got{l}^*$, and $N=(\got{l}/\got{t}_{\scriptop{h}})^*$. 
See Proposition \ref{Aprop}. 

The symplectic form on this bundle of Delzant manifolds 
is given explicitly by means of the formula 
(\ref{A*sigma}), in terms of the antisymmetric bilinear form  
$\sigma ^{\got{t}}$ on $\got{t}$ introduced in Lemma \ref{constlem}, 
the antisymmetric bilinear mapping $c:N\times N\to\got{l}$ 
introduced in Proposition \ref{liftprop}, and the 
symplectic form $\sigma _{\scriptop{h}}$ on the Delzant manifold 
$M_{\scriptop{h}}$.  
In this way we obtain an explicit global model  
for our symplectic manifold $(M,\,\sigma )$ with symplectic 
$T$\--action. 

\subsection{An extension $G$ of $T$ by $N$ 
which acts on $M$}
\label{Gss}
In the sequel, ${\cal X}^{\infty}(M)$ denotes 
the Lie algebra of all smooth 
vector fields on $M$, provided with the 
Lie brackets $[u,\, v]$ of $u,\,\ v \in {\cal X}^{\infty}(M)$ 
such that $[u,\, v]\, f=u\, (v\, f)-v\, (u\, f)$ 
for every $f\in\op{C}^{\infty}(M)$. 
We denote the flow after time $t\in\R$ of 
$v\in {\cal X}^{\infty}(M)$ by $\op{e}^{t\, v}$. 
This defines an exponential mapping $v\mapsto\op{e}^v$ 
from ${\cal X}^{\infty}(M)$ to the group $\op{Diff}^{\infty}(M)$ 
of all smooth diffeomorphisms of $M$, which is analogous to 
the exponential mapping $\op{exp}$ from the Lie algebra 
of any Lie group to the Lie group. 

Let $\got{l}^*\ni\xi\mapsto L_{\xi}$ be an admissible 
lift as in Proposition \ref{liftprop}. 
For each $\zeta\in N$, $L_{\zeta}$ is a smooth 
vector field on $M$, see Definition \ref{liftdef}, and   
because $M$ is compact, its flow $\op{e}^{t\, L_{\zeta}}:M\to M$ 
is defined for all $t\in\R$. 

A Lie algebra $\got{g}$ is called {\em two\--step nilpotent} 
if $[[X,\, Y],\, Z]=0$ for all $X,\, Y,\, Z\in\got{g}$. 
Because the vector fields $L_{\zeta}$, $\zeta\in N$ 
commute with the $X_M$, $X\in\got{t}$, and the 
$X_M$, $X\in\got{t}$, commute with each other, it follows 
from iii) in Proposition \ref{liftprop} that 
the linear span of the $X_M$, $X\in\got{t}$, 
and the $L_{\zeta}$, $\zeta\in N$, 
is a two\--step nilpotent Lie subalgebra $\got{g}_M$ 
of ${\cal X}^{\infty}(M)$. 
Moreover, if we provide $\got{g}:=\got{t}\times N$ 
with the structure of a two\--step nilpotent 
Lie algebra defined by 
\begin{equation}
[(X,\,\zeta ),\, (X',\,\zeta ')]=\, -(c(\zeta ,\,\zeta '),\, 0),
\quad (X,\,\zeta ),\, (X',\,\zeta ')\in \got{g}=\got{t}\times N, 
\label{gotg}
\end{equation}
then the mapping $(X,\,\zeta )\mapsto X_M+L_{\zeta}$ is 
an injective anti\--homomorphism of Lie algebras from $\got{g}$ 
to ${\cal X}^{\infty}(M)$, with image equal to $\got{g}_M$.  
 
The vector space $\got{t}\times N$, provided with the product 
\begin{equation}
(X,\,\zeta )\, (X',\,\zeta ')=
(X+X'-c(\zeta ,\,\zeta ')/2,\,\zeta +\zeta '),
\quad (X,\,\zeta ),\, (X',\,\zeta ')\in \got{t}\times N, 
\label{G}
\end{equation}
is a two\--step nilpotent Lie group with Lie algebra equal to ${\got g}$ 
and the identity as the exponential mapping. It follows that 
the mapping 
\begin{equation}
(X,\,\zeta )\mapsto \op{e}^{X_M+L_{\zeta}}=\op{e}^{X_M}\circ\op{e}^{L_{\zeta}}
\label{Gaction}
\end{equation}
is a (left) action of the group $\got{t}\times N$ on $M$, 
that is a homomorphism from 
the group $\got{t}\times N$ 
to the group $\op{Diff}^{\infty}(M)$, with infinitesimal 
action given by $(X,\,\zeta )\mapsto X_M+L_{\zeta}$. 
It follows that 
\begin{equation}
\op{e}^{X_M+L_{\zeta}}\circ\op{e}^{X'_M+L_{\zeta '}}
=\op{e}^{(X+X'-c(\zeta ,\,\zeta ')/2)_M+L_{\zeta +\zeta '}}.  
\label{expproduct}
\end{equation}

The kernel of the homomorphism (\ref{Gaction}) is equal to 
the discrete normal subgroup $T_{\Z}\times\{ 0\}$ of $
\got{t}\times N$, 
in which $T_{\Z}=\op{ker}\op{exp}$ is the integral 
lattice in the Lie algebra $\got{t}$ of $T$. 
It follows that the connected Lie group $G=T\times N
\simeq (\got{t}/T_{\Z})\times N$ acts smoothly on 
$M$, where in $T\times N$ we have the product 
\begin{equation}
(t,\,\zeta )\, (t',\,\zeta ')=
(t\, t'\,\op{e}^{-c(\zeta ,\,\zeta ')/2},\,\zeta +\zeta ')
\label{Hproduct}
\end{equation}
and the action is given by 
\begin{equation}
(t,\,\zeta )\mapsto t_M\circ\op{e}^{L_{\zeta}}.
\label{Haction}
\end{equation}

Note that the Lie algebra of $G$ is equal to 
the previously introduced two\--step nilpotent 
Lie algebra $\got{g}=\got{t}\times N$. 
Also note that the $T$\--orbit map $\pi :M\to M/T$ 
intertwines the action of $G$ on $M$ with the 
translational action of $N$ on $M/T$, 
in the sense that $\pi ((t,\,\zeta )\cdot x)=
\pi(x)+\zeta$ for every $(t,\,\zeta )\in G=T\times N$.  

\subsection{The holonomy of the connection}
\label{holss}
Let $\got{l}^*\ni\xi\mapsto L_{\xi}$ be an admissible 
connection as in Proposition \ref{liftprop}.  
For each $\zeta\in P$ and $p\in M/T$, the curve 
$\gamma_{\zeta}(t):=p+t\,\zeta$, 
$0\leq t\leq 1$, is a loop in $M/T$. If 
$x\in M$ and $p=\pi (x)$, then the curve 
$\delta (t)=\op{e}^{t\, L_{\zeta}}(x)$, $0\leq t\leq 1$, 
is called the {\em horizontal lift in $M$ of the loop $\gamma _{\zeta}$ 
which starts at $x$}, because $\delta (0)=x$, 
$\delta '(t)=L_{\zeta}(\delta (t))$ is a horizontal 
tangent vector which is mapped by $\op{T}_{\delta (t)}\!\pi$ 
to the constant vector $\zeta$, which implies that  
$\pi (\delta (t))=\gamma _{\zeta}(t)$, $0\leq t\leq 1$.  
The element of $T$ which maps 
the initial point $\delta (0)=x$ to the end point 
$\delta (1)$ is called the {\em holonomy 
$\tau _{\zeta}(x)$ of the 
loop $\gamma _{\zeta}$ and the intial point $x$} 
with respect to the given 
connection. Because $\delta (1)
=\op{e}^{L_{\zeta}}(x)$, we have 
$\tau _{\zeta}(x)\cdot x=\op{e}^{L_{\zeta}}(x)$. 
In Lemma \ref{TPlem} below we investigate the dependence 
of the holonomy element $\tau _{\zeta}(x)\in T$ 
on the point $x\in M$ and the period $\zeta\in P$. 

\begin{lemma}
Let $\zeta\in N$. Then the following conditions 
are equivalent.
\begin{itemize}
\item[i)] There exists an $x\in M$ and a $t\in T$ such that~ 
$\op{e}^{L_{\zeta}}(x)=t\cdot x$. 
\item[ii)] $\zeta\in P$, where $P$ is the period 
group in $N$ for the translational action of 
$N$ on $M/T$, as defined in Lemma 
\ref{perlem} with $Q=M/T$, $V=\got{l}^*$, and 
$N=(\got{l}/\got{t}_{\scriptop{h}})^*$. 
\item[iii)] The diffeomorphism $\op{e}^{L_{\zeta}}$ leaves all 
$T$\--orbits in $M$ invariant.
\end{itemize}

For each $\zeta\in P$ there is a unique $T$\--invariant smooth mapping 
$\tau _{\zeta}:M\to T$ such that $\op{e}^{L_{\zeta}}(x)=
\tau _{\zeta}(x)\cdot x$ for every $x\in M$. We have 
\begin{equation}
\tau _{\zeta}(t\cdot\op{e}^{L_{\zeta '}}(x))
=\op{e}^{c(\zeta ,\,\zeta ')}\,\tau _{\zeta}(x)
\label{tauzeta}
\end{equation}
for every $(t,\,\zeta ')\in T\times N$. 

We have $c(\zeta ,\,\zeta ')\in T_{\Z}$ 
whenever $\zeta ,\,\zeta '\in P$, and $T\times P$ is a 
commutative subgroup of $G$.  

Finally, the mapping $\tau _{\zeta}:M\to T$ is constant on 
every fiber of the fibration of $M$ into 
Delzant submanifolds introduced in Proposition 
\ref{embedprop}, and 
satisfies 
\begin{equation}
\tau _{\zeta '}(x)\,\tau _{\zeta}(x)=\tau _{\zeta +\zeta '}(x)
\,\op{e}^{c(\zeta ',\,\zeta )/2},\quad x\in M,\quad \zeta,\,\zeta '\in P. 
\label{tauplus}
\end{equation} 
\label{TPlem}
\end{lemma}
\begin{proof}
Because the action of $\op{e}^{L_{\zeta}}$ on the $T$\--orbits 
is equal to the transformation $p\mapsto p+\zeta$ in $M/T$, 
the equivalence between i), ii), iii) follows from Lemma \ref{perlem} 
with $Q=M/T$, $V=\got{l}^*$, and $N=(\got{l}/\got{t}_{\scriptop{h}})^*$. 

If $\zeta\in P$, then $\op{e}^{L_{\zeta }}$ leaves 
each $T$\--orbit invariant. Because, for every 
$\zeta\in N$, $\op{e}^{L_{\zeta }}$ 
commutes with the $T$\--action, this implies the 
existence of the smooth mapping $\tau _{\zeta}$ in view 
of Lemma \ref{hslem}.  

In order to show that (\ref{tauzeta}) holds, we observe that 
\begin{eqnarray*}
\op{e}^{L_{\zeta '}}(\tau _{\zeta}(\op{e}^{L_{\zeta '}}(x))\cdot x)
&=&\tau _{\zeta}(\op{e}^{L_{\zeta '}}(x))\cdot\op{e}^{L_{\zeta '}}(x)
=\op{e}^{L_{\zeta}}(\op{e}^{L_{\zeta '}}(x))
=(\op{e}^{L_{\zeta}}\circ\op{e}^{L_{\zeta '}}\circ\op{e}^{-L_{\zeta}})
(\op{e}^{L_{\zeta}}(x))\\
&=&\op{e}^{L_{\zeta '}+[L_{\zeta},\, L_{\zeta '}]}
(\tau _{\zeta}(x)\cdot x)
=\op{e}^{L_{\zeta '}}(\op{e}^{c(\zeta ,\,\zeta ')}\cdot
(\tau _{\zeta}(x)\cdot x)),
\end{eqnarray*}
which implies that  
$
\tau _{\zeta}(\op{e}^{L_{\zeta '}}(x))
=\op{e}^{c(\zeta ,\,\zeta ')}\,\tau _{\zeta}(x).
$
In combination with the $T$\--invariance of $\tau _{\zeta}$ 
this yields (\ref{tauzeta}). 

If $\zeta '\in P$, then we have for every $x\in M$ 
that $\op{e}^{L_{\zeta '}}(x)\in T\cdot x$, 
hence $\tau _{\zeta}(\op{e}^{L_{\zeta '}}(x))=\tau _{\zeta}(x)$, 
which in view of (\ref{tauzeta}) implies that 
$\op{e}^{c(\zeta ,\,\zeta ')}=1$, hence $c(\zeta ,\,\zeta ')\in T_{\Z}$. 
The fact that $c(\zeta ,\,\zeta ')\in T_{\Z}$ for all 
$\zeta ,\,\zeta '\in P$ implies in view of (\ref{Hproduct}) that 
$T\times P$ is a commutative 
subgroup of $T\times N$. 

Because $L_{\zeta}$ commutes with all $L_{\eta}$, 
$\eta\in C:=(\got{l}/\got{l}\cap\got{t}_{\scriptop{f}})^*
\simeq {\got{t}_{\scriptop{h}}}^*$, see 
ii) in Proposition \ref{liftprop}, we have 
\[
\op{e}^{L_{\eta}}(\tau _{\zeta}(\op{e}^{L_{\eta}}(x))\cdot x)
=\tau _{\zeta}(\op{e}^{L_{\eta}}(x))\cdot \op{e}^{L_{\eta}}(x)
=\op{e}^{L_{\zeta}}(\op{e}^{L_{\eta}}(x))
=\op{e}^{L_{\eta}}(\op{e}^{L_{\zeta}}(x))
=\op{e}^{L_{\eta}}(\tau _{\zeta}(x)\cdot x), 
\]
which for regular $x$ implies that 
$\tau _{\zeta}(\op{e}^{L_{\eta}}(x))=\tau _{\zeta}(x)$. 
By continuity this identity extends to all $x\in M$. 
Because also $\tau _{\zeta}(t\cdot x)=\tau _{\zeta}(x)$ 
for all $t\in T_{\scriptop{h}}$, it follows from the definition 
in Proposition 
\ref{embedprop} of the fibration of $M$ into 
Delzant submanifolds, that $\tau _{\zeta}$ is 
constant on its fibers. 

If $\zeta ,\,\zeta '\in P$, then we obtain, using 
(\ref{expproduct}), that 
\begin{eqnarray*}
\tau_{\zeta '}(x)\cdot (\tau_{\zeta}(x)\cdot x)
&=&\tau_{\zeta '}(x)\cdot\op{e}^{L_{\zeta }}(x)
=\op{e}^{L_{\zeta}}(\tau_{\zeta '}(x)\cdot x)
=(\op{e}^{L_{\zeta}}\circ\op{e}^{L_{\zeta '}})(x)\\
&=&\op{e}^{-c(\zeta ,\,\zeta ')_M/2}\cdot\op{e}^{L_{\zeta +\zeta '}}(x) 
=\op{e}^{c(\zeta ',\,\zeta )_M/2}\cdot
(\tau _{\zeta +\zeta '}(x)\cdot x), 
\end{eqnarray*}
which implies (\ref{tauplus}). 
\end{proof}

Let $\varepsilon ^l$, $1\leq l\leq d_N:=\op{dim}N$,  
be a $\Z$\--basis of $P$. For any $\zeta\in P$ 
we have $\zeta =\sum_l\,\zeta _l\,\varepsilon ^l$ 
for unique integral coordinates $\zeta _l\in\Z$. 
With the notation $c^{l\, l'}
:=c(\varepsilon ^l,\,\varepsilon ^{l'})\in\got{l}\cap T_{\Z}$, 
the formula (\ref{tauplus}) leads to the formula 
\begin{equation}
\tau _{\zeta}(x)=
\op{e}^{\sum_{l<l'}\,\zeta _l\,\zeta _{l'}\, c^{l\, l'}/2}
\,\prod_{l=1}^{d_N}\,\tau _{\varepsilon ^l}(x)^{\zeta _l}
\label{tauzetaform}
\end{equation}
for $\tau _{\zeta}(x)$ in terms of the 
elements $\tau _{\varepsilon ^l}(x)\in T$. In other words, 
all holonomies at a given point $x\in M$ can be expressed 
in terms of the holonomies of the basic loops 
$\gamma _{\varepsilon ^l}$, $1\leq l\leq d_N$, 
by means of the formula (\ref{tauzetaform}). 

\subsection{$M$ as a $G$\--homogeneous bundle with the 
Delzant manifold as fiber}
In this subsection, let 
$(M_{\scriptop{h}},\,\sigma _{\scriptop{h}},\, T_{\scriptop{h}})$ 
be one of the Delzant submanifolds of $(M,\,\sigma ,\, T)$ 
in Proposition \ref{embedprop}. That is, $M_{\scriptop{h}}$ is an integral 
manifold $I$ of the distribution $H$, and $\sigma _{\scriptop{h}}=
{\iota _I}^*\sigma$, if $\iota _I$ denotes the inclusion mapping 
from $I$ to $M$. 
Recall that all Delzant manifolds with the same Delzant polytope 
are $T_{\scriptop{h}}$\--equivariantly symplectomorphic, 
which means that one may identify $(M_{\scriptop{h}},\,\sigma _{\scriptop{h}},\, 
T_{\scriptop{h}})$ 
with any favourite explicit model of a Delzant 
manifold with Delzant polytope $\Delta$. 
We will construct a model for our symplectic 
$T$\--manifold $(M,\,\sigma ,\, T)$ by means of the 
mapping $A:G\times M_{\scriptop{h}}\to M$ which is defined by 
\begin{equation}
A((t,\,\zeta ),\, x)=t\cdot\op{e}^{L_{\zeta}}(x),\quad 
t\in T,\; \zeta\in N,\; x\in M_{\scriptop{h}}. 
\label{A}
\end{equation}
Write $\tau _{\zeta}$ for the common value of the 
$\tau _{\zeta}(x)$ for all $x\in M_{\scriptop{h}}$, 
see Lemma \ref{TPlem}. 
Define 
\begin{equation}
H:=\left\{ (t,\,\zeta )\in G\mid \zeta\in P\;
\mbox{\rm and}\; t\,\tau _{\zeta}\in T_{\scriptop{h}}\right\} .
\label{Hdef}
\end{equation}
Then $H$ is a closed Lie subgroup of $G$, commutative because 
$T\times P$ is commutative, see Lemma \ref{TPlem}. 
Furthermore, 
\begin{equation}
((t,\,\zeta ),\, x)\mapsto (t\,\tau _{\zeta})\cdot x:
H\times M_{\scriptop{h}}\to M_{\scriptop{h}}
\label{Hact}
\end{equation}
defines a smooth action of $H$ on the Delzant manifold 
$M_{\scriptop{h}}$. 
\begin{proposition}
The mapping {\em (\ref{A})} induces a 
diffeomorphism $\alpha$ from $G\times _HM_{\scriptop{h}}$ onto $M$, 
where $h\in H$ acts on $G\times M_{\scriptop{h}}$ by sending 
$(g,\, x)$ to $(g\, h^{-1},\, h\cdot x)$. 

The diffeomorphism 
$\alpha$ intertwines the action of $G$ on $G\times _HM_{\scriptop{h}}$, 
which is induced by the action 
$(g,\, (g', x))\mapsto (g\, g', x)$ 
of $G$ on $G\times M_{\scriptop{h}}$, with the action of 
$G$ on $M$, and therefore also the action of the normal 
subgroup $T\times\{ 0\}\simeq T$ of $G$ with the 
$T$\--action on $M$. The projection 
$(g,\ x)\mapsto g:G\times M_{\scriptop{h}}\to G$ 
induces a $G$\--equivariant smooth fibration 
$\psi :G\times _HM_{\scriptop{h}}\to G/H$, and  
$\delta =\psi\circ\alpha ^{-1}:M\to G/H$ is a $G$\--equivariant 
smooth fibration of which the Delzant submanifolds of $M$ 
introduced in Proposition \ref{embedprop} are the fibers. 
\label{Aprop}
\end{proposition}
\begin{proof}
Let $x_0\in M_{\scriptop{h}}$ and $y\in M$. 
For each $\zeta\in N$, the projection 
$\pi _{N/P}:M\to N/P$ 
defined in Lemma \ref{pi2lem} intertwines the diffeomorphism 
$\op{e}^{L_{\zeta}}$ in $M$ with the translation in 
$N/P$ over the vector $\zeta$. Because these translations 
act transitively on $N/P$, there exists a $\zeta\in N$  
such that $\pi _{N/P}(y)=\pi _{N/P}(x_0)+\zeta$, 
which implies that $\op{e}^{-L_{\zeta}}(y)$ belongs to 
the same fiber $F$ of $\pi _{N/T}$ as $x_0$. With such 
a choice of $\zeta$, it follows from 
Lemma \ref{embedlem} that there exists  
$t_{\scriptop{f}}\in T_{\scriptop{f}}$ 
and $x\in M_{\scriptop{h}}$ such that 
$\op{e}^{-L_{\zeta}}(y)=t_{\scriptop{f}}\cdot x$, or equivalently 
$y=t_{\scriptop{f}}\cdot\op{e}^{L_{\zeta}}(x)$. This shows that 
already the restriction of $A$ to 
$(T_{\scriptop{f}}\times N)\times M_{\scriptop{h}}$ is surjective. 

Let $g,\, g'\in G$, $x,\, x'\in M_{\scriptop{h}}$ and 
$g\cdot x=g'\cdot x'$. Then $x'=h\cdot x$ 
in which $h:=(g')^{-1}\, g$. Write $h=(t,\,\zeta )$ 
with $t\in T$ and $\zeta\in N$. 
Then 
\[
\pi _{N/P}(x)=\pi _{N/P}(x')=\pi _{N/P}(t\cdot\op{e}^{L_{\zeta}}(x))
=\pi _{N/P}(\op{e}^{L_{\zeta}}(x))=\pi _{N/P}(x)+\zeta
\]
implies that $\zeta\in P$, and it follows from 
Lemma \ref{TPlem} that 
\[
x'=t\cdot\op{e}^{L_{\zeta}}(x)=t\cdot \tau _{\zeta}\cdot x
=(t\,\tau _{\zeta})_{\scriptop{f}}
\cdot ((t\,\tau _{\zeta})_{\scriptop{h}}\cdot x.
\]
Because $x'$ and $(t\,\tau _{\zeta})_{\scriptop{h}}\cdot x$ belong to the 
same integral manifold $I$ of $H$, it follows from 
Lemma \ref{embedlem} that the element 
$(t\,\tau _{\zeta})_{\scriptop{f}}$ of $T_{\scriptop{f}}$ 
is equal to the identity element, hence $t\,\tau _{\zeta}\in T_{\scriptop{h}}$ 
and $x'=(t\,\tau _{\zeta})\cdot x$. In other words, 
$h\in H$, $g' =g\, h^{-1}$ and $x'=(t\,\tau _{\zeta})\cdot x$. 
This proves that the mapping $A$ induces a bijective 
mapping $\alpha$ from $G\times _HM_{\scriptop{h}}$ onto $M$. 

The closedness of $H$ in $G$ implies that 
the right action of $H$ on $G$ is proper and free, hence  
the action of $H$ on $G\times M_{\scriptop{h}}$ 
is proper and free, and the orbit space $G\times _HM_{\scriptop{h}}$ 
has a unique smooth structure for which the projection 
$G\times M_{\scriptop{h}}\to G\times _HM_{\scriptop{h}}$ 
is a principal $H$\--bundle. 
With respect to this smooth structure on 
$G\times _HM_{\scriptop{h}}$, the mapping 
$\alpha :G\times _HM_{\scriptop{h}}\to M$ 
is smooth. 
The transversality to $\op{T}\! M_{\scriptop{h}}$ of the span of $Z_M$, 
$Z\in\got{t}_{\scriptop{f}}$ and 
the  $L_{\zeta}$, $\zeta\in N$, implies that 
at every point the tangent mapping of $A$ is surjective. 
Hence $\alpha$ is a submersion, and because $\alpha$ is 
bijective, it follows from the inverse mapping theorem that 
$\alpha$ is a diffeomorphism.  
The other statements in 
the proposition are general facts about induced 
fiber bundles $G\times _HM_{\scriptop{h}}$ over $G/H$ 
with fiber $M_{\scriptop{h}}$, 
see for instance \cite[Sec. 2.4]{dk}.  
\end{proof}

\begin{remark}
On $G/H$ we still have the free action of the torus $T/T_{\scriptop{h}}$, 
which exhibits $G/H$ as a principal $T/T_{\scriptop{h}}$\--bundle 
over the torus $(G/H)/T\simeq N/P$. Palais and Stewart \cite{ps} 
showed that every principal torus bundle over a torus is diffeomorphic to 
a nilmanifold for a two\--step nilpotent Lie group. In this remark 
we will give an explicit nilmanifold description of $G/H$. 

The Hamiltonian torus $T_{\scriptop{h}}$, or rather 
the identity component $H^o=T_{\scriptop{h}}\times\{ 0\}$ 
of $H$, is a closed normal Lie subgroup 
of both $G=T\times N$ and $H$, and the mapping 
$(g\, H^o)\, (H/H^o)\mapsto g\, H$ 
is a $G$\--equivariant 
diffeomorphism from $(G/H^o)/(H/H^o)$ onto $G/H$. 
The group structure in $G/H^o=(T/T_{\scriptop{h}})\times N$ is defined by 
\begin{equation}
(t,\,\zeta)\, (t',\,\zeta ')
=(t\, t'\,\op{e}^{-c_{\got{l}/\got{t}_{\tinyop{h}}}
(\zeta ,\,\zeta ')/2},\,\zeta +\zeta '), 
\quad t,\, t'\in T/T_{\scriptop{h}},\quad\zeta ,\,\zeta '\in N,
\label{G2product}
\end{equation}
and $c_{\got{l}/\got{t}_{\scriptop{h}}}
:N\times N\to\got{l}/\got{t}_{\scriptop{h}}$ 
is equal to $c:N\times N\to\got{l}$, followed by the 
projection $\got{l}\to\got{l}/\got{t}_{\scriptop{h}}$. 
This exhibits $G/H^o$ as a two\--step nilpotent Lie group 
with universal covering equal to $(\got{t}/\got{t}_{\scriptop{h}})\times N$ 
and covering group $(T/T_{\scriptop{h}})_{\Z}\simeq 
T_{\Z}/(T_{\scriptop{h}})_{\Z}$. Also note that 
$\iota :\zeta\mapsto ({\tau _{\zeta}}^{-1},\,\zeta )\, H^o$ 
is an isomorphism from the period group $P$ onto $H/H^o$. 

In view of (\ref{Hdef}), we conclude that 
the compact homogeneous $G$\--space $G/H$ is isomorphic 
to the quotient of the simply connected two\--step 
nilpotent Lie group $(\got{t}/\got{t}_{\scriptop{h}})\times N$ 
by the discrete subgroup of $(\got{t}/\got{t}_{\scriptop{h}})\times N$ 
which consists of all $(Z,\,\zeta )\in 
(\got{t}/\got{t}_{\scriptop{h}})\times P$ such that 
$\op{e}^Z\,\tau _{\zeta}\in T_{\scriptop{h}}$. 
\label{Horem}
\end{remark}

\subsection{The symplectic form on the global model}
In Proposition \ref{Aprop} we have described the global model  
$M_{\scriptop{model}}:=G\times _HM_{\scriptop{h}}$ 
for the $T$\--manifold $M$, where the multiplication 
in the Lie group $G=T\times N$ 
is defined by (\ref{Hproduct}). We now describe the symplectic form on 
$M_{\scriptop{model}}$. 

\begin{proposition}
Let $\omega$ be the pull\--back of $\sigma$ to $G\times M_{\scriptop{h}}
=(T\times N)\times M_{\scriptop{h}}$ 
by means of the mapping $A$ in (\ref{A}). 
Let $\delta a=((\delta t,\,\delta\zeta ),\,\delta x)$ 
and $\delta' a=((\delta 't,\,\delta '\zeta ),\,\delta ' x)$ 
be tangent vectors to $G\times M_{\scriptop{h}}$ at 
$a=((t,\,\zeta ),\, x)$, 
where we identify each tangent space of the torus $T$ with 
$\got{t}$. Write $X=\delta t+c(\delta\zeta ,\,\zeta )/2$ 
and $X'=\delta 't+c(\delta '\zeta ,\,\zeta )/2$. Then
\begin{eqnarray}
\omega _a(\delta a,\,\delta 'a)
&=&\sigma ^{\got{t}}(\delta t,\,\delta 't)
+\delta\zeta ({X'}_{\got{l}})
-\delta '\zeta (X_{\got{l}})
-\mu (x)(c_{\scriptop{h}}(\delta\zeta ,\,\delta '\zeta ))
\nonumber\\
&&+\, (\sigma _{\scriptop{h}})_x(\delta x,\, 
({X'}_{\scriptop{h}})_{M_{\scriptop{h}}}(x))
-(\sigma _{\scriptop{h}})_x(\delta 'x,\, 
(X_{\scriptop{h}})_{M_{\scriptop{h}}}(x))
\nonumber\\
&&+\, (\sigma _{\scriptop{h}})_x(\delta x,\,\delta 'x).
\label{A*sigma}
\end{eqnarray}
Here $X_{\scriptop{h}}$ denotes the 
$\got{t}_{\scriptop{h}}$\--component of $X\in\got{t}$ with respect 
to the direct sum decomposition $\got{t}_{\scriptop{h}}\oplus 
\got{t}_{\scriptop{f}}$. 

If $\pi _{M_{\scriptop{model}}}$ denotes the canonical projection from 
$G\times M_{\scriptop{h}}$ onto 
$M_{\scriptop{model}}:=G\times _HM_{\scriptop{h}}$, 
then the $T$\--invariant symplectic form 
$\sigma _{\scriptop{model}}:=\alpha ^*\,\sigma$ on $M_{\scriptop{model}}$ 
is the unique two\--form $\beta$ on 
$M_{\scriptop{model}}$ such that 
$\omega ={\pi _{M_{\scriptop{model}}}}^*\,\beta$. 
\label{omegalem}
\end{proposition}
\begin{proof}
It follows from (\ref{expproduct}) that 
\[
\op{e}^{L_{\zeta '+\zeta}}=\op{e}^{c(\zeta ',\,\zeta )_M/2}
\circ\op{e}^{L_{\zeta '}}\circ\op{e}^{L_{\zeta}}.
\]
Therefore, if we substitute $\zeta '=\epsilon\,\delta\zeta$ 
and differentiate with respect to $\epsilon$ at $\epsilon =0$, we 
get the vector $c(\delta\zeta ,\,\zeta )_M/2+L_{\delta\zeta}$ 
at the image point under the mapping $\op{e}^{L_{\zeta}}$. 
Because $\op{e}^{L_{\zeta}}$ commutes with the $T$\--action, 
it follows, with the notations 
$y=A(a)$ and $B=t_M\circ\op{e}^{L_{\zeta}}$, 
that 
\[
\delta y=(\op{T}_a\! A)(\delta a)=(X_{M}
+L_{\delta\zeta})(y)
+(\op{T}_x\! B)(\delta x), 
\]
in which $X=\delta t+c(\delta\zeta ,\,\zeta )/2$. 

If $x$ is a regular point in $M_{\scriptop{h}}$, then we can write 
$\delta x=(Y_M+L_{\eta})(p)$ 
for uniquely determined $Y\in\got{t}_{\scriptop{h}}$ and 
$\eta\in C=(\got{l}/\got{l}\cap\got{t}_{\scriptop{f}})^*
\simeq (\got{t}_{\scriptop{h}})^*$. 
The vector fields $Y_M$, 
$Y\in\got{t}_{\scriptop{h}}$, and $L_{\eta}$, $\eta\in C$, 
commute with the vector fields 
$X_M$, $X\in\got{t}$, and $L_{\zeta}$, $\zeta\in N$, 
because of ii) in Proposition \ref{liftprop} and the fact that 
all the vector fields are $T$\--invariant. 
Therefore 
$(\op{T}_x\! B)(\delta x)=(Y_M+L_{\eta})(y)$,
and we obtain that 
$(\op{T}_a\! A)(\delta a)$ is equal to the 
value at $y=A(a)$ of the vector field 
$(X+Y)_M+L_{\delta\zeta +\eta}$. 

In view of (\ref{stronglift}) and 
iv), v) in Proposition \ref{liftprop}, 
the symplectic product of this vector with the one in which 
$\delta t$, $\delta\zeta$, $Y$, $\eta$ are 
replaced by $\delta 't$, $\delta '\zeta$, $Y'$, $\eta '$, respectively, 
is equal to 
\[
(\delta\zeta +\eta)((X'+Y')_{\got{l}})-
(\delta '\zeta +\eta ')((X+Y)_{\got{l}})
+\sigma _y(L_{\delta\zeta}(y),\, L_{\delta '\zeta }(y)), 
\]
in which $X=\delta t+c(\delta\zeta,\, \zeta )/2+Y$ 
and $X'=\delta 't+c(\delta '\zeta ,\,\zeta )/2+Y'$.  
Collecting terms and using the equations   
$\eta ({X'}_{\got{l}})=\eta ({X'}_{\scriptop{h}})
=(\sigma _{\scriptop{h}})_x(\delta x,\, (X'_{\scriptop{h}})
_{M_{\scriptop{h}}}(x))$, 
$\eta '(X_{\got{l}})=\eta '(X_{\scriptop{h}})
=(\sigma _{\scriptop{h}})_x(\delta 'x,\, (X_{\scriptop{h}})
_{M_{\scriptop{h}}}(x))$, 
$\eta (Y')-\eta '(Y)=(\sigma _{\scriptop{h}})_x(\delta x,\,\delta 'x)$, 
and vi) in Proposition \ref{liftprop}, we arrive at 
(\ref{A*sigma}). 

Because $A=\alpha\circ\pi _{M_{\scriptop{model}}}$, we have 
$\omega =A^*\,\sigma ={\pi _{M_{\scriptop{model}}}}^* (\alpha ^*\,\sigma )
={\pi _{M_{\scriptop{model}}}}^*\,\sigma _{\scriptop{model}}$. 
The uniqueness 
in the last statement follows because $\pi _{M_{\scriptop{model}}}$ is 
a submersion. 
\end{proof}

\begin{lemma}
Let $T_{\scriptop{f}}$ be a complementary torus to the 
Hamiltonian torus $T_{\scriptop{h}}$ in $T$. Then 
the following conditions are equivalent. 
\begin{itemize}
\item[a)] 
$(M,\,\sigma ,\, T)$ is $T$\--equivariantly symplectomorphic 
to the Cartesian product of a symplectic $T_{\scriptop{f}}$\--
space $(M_{\scriptop{f}},\,\sigma _{\scriptop{f}},\, T_{\scriptop{f}})$ 
on which the $T_{\scriptop{f}}$\--action is free and a 
Delzant manifold $(M_{\scriptop{h}},\,\sigma _{\scriptop{h}}, 
\, T_{\scriptop{h}})$. Here $t\in T$ acts on $M_{\scriptop{f}}\times 
M_{\scriptop{h}}$ by sending $(x_{\scriptop{f}},\, x_{\scriptop{h}})$ 
to $(t_{\scriptop{f}}\cdot x_{\scriptop{f}}
,\, t_{\scriptop{h}}\cdot x_{\scriptop{h}})$, 
if $t=t_{\scriptop{f}}\, t_{\scriptop{h}}$ with 
$t_{\scriptop{f}}\in T_{\scriptop{f}}$  
and $t_{\scriptop{h}}\in T_{\scriptop{h}}$. 
\item[b)] $c(P\times P)\subset
\got{t}_{\scriptop{f}}$. 
\item[c)] $c(N\times N)\subset
\got{t}_{\scriptop{f}}$. 
\item[d)] 
The $\got{t}_{\scriptop{h}}$\--component 
$c_{\scriptop{h}}$ of $c$ in the direct sum decomposition 
$\got{l}=\got{t}_{\scriptop{h}}\oplus 
(\got{l}\cap\got{t}_{\scriptop{f}})$ is equal to zero. 
\end{itemize}
\label{e0rem}
\end{lemma}
\begin{proof}
The equivalence between b) and c) follows from the fact that 
$P$ has a $\Z$\--basis which is an $\R$\--basis of $N$. 
The equivalence between c) and d) is obvious. 

If $(M,\,\sigma ,\, T)$ is equal to the Cartesian 
product of a Delzant manifold 
$(M_{\scriptop{h}},\,\sigma _{\scriptop{h}},\, T_{\scriptop{h}})$ 
and a symplectic $T_{\scriptop{f}}$\--space 
$(M_{\scriptop{f}},\,\sigma _{\scriptop{f}},\, T_{\scriptop{f}})$ 
for which the 
$T_{\scriptop{f}}$\--action on $M_{\scriptop{f}}$ is free, then we can choose 
the $L_{\zeta}$ in the direction of the second component 
$M_{\scriptop{f}}$. In this case we have for every 
$\zeta ,\,\zeta '\in N$ that $[L_{\zeta},\, L_{\zeta '}]\in 
\got{t}_{\scriptop{f}}$, which means that the antisymmetric 
bilinear mapping $c:N\times N\to\got{l}$ 
in Proposition \ref{liftprop} has the property that 
$c(N\times N)\subset\got{t}_{\scriptop{f}}$, 
or equivalently $c_{\scriptop{h}}=0$. 

For the converse, assume that $c(N\times N)\subset\got{t}_{\scriptop{f}}$, 
which implies that $c_{\scriptop{h}}=0$.  
Then the Lie group 
$G$ is equal to the Cartesian product 
$T_{\scriptop{h}}\times G_{\scriptop{f}}$, 
in which $G_{\scriptop{f}}=T_{\scriptop{f}}\times N$, where the product in 
$G_{\scriptop{f}}$ is defined as in (\ref{Hproduct}) with $T$ 
replaced by $T_{\scriptop{f}}$. 
According to Subsection \ref{holinvariantss} we can multiply the elements 
$\tau _{\varepsilon ^l}(x)$, for $x\in M_{\scriptop{h}}$, 
by any element of $\op{exp}(\got{l})$. Because 
$\got{t}_{\scriptop{h}}\subset\got{l}$, it follows 
that we can arrange that $\tau _{\varepsilon ^l}(x)\in 
T_{\scriptop{f}}$ for every $1\leq l\leq d_N$, 
and then it follows from 
(\ref{tauzetaform}) that 
$\tau _{\zeta}\in T_{\scriptop{f}}$ for every $\zeta\in P$. 
The mapping $\iota :\zeta\mapsto ({\tau_{\zeta}}^{-1},\,\zeta )$ 
is a homomorphism from $P$ onto a discrete cocompact subgroup 
of $G_{\scriptop{f}}$. Write $M_{\scriptop{f}}:=G_{\scriptop{f}}/\iota (P)$.  
It follows that the mapping 
\begin{equation}
A_{\scriptop{f}}:((t_{\scriptop{f}},\,\zeta ), x)\mapsto 
t_{\scriptop{f}}\cdot\op{e}^{L_{\zeta}}(x)
:G_{\scriptop{f}}\times M_{\scriptop{h}}\to M 
\label{A2}
\end{equation}
induces a diffeomorphism $\alpha _{\scriptop{f}}$ from 
$M_{\scriptop{f}}\times M_{\scriptop{h}}$ onto 
$M$. Moreover, it follows from (\ref{A*sigma}) that the symplectic form 
${\alpha _{\scriptop{f}}}^*\,\sigma$ 
on $M_{\scriptop{f}}\times M_{\scriptop{h}}$ is equal to 
${\pi _{\scriptop{f}}}^*\,\sigma _{\scriptop{f}}+
{\pi  _{\scriptop{h}}}^*\,\sigma _{\scriptop{h}}$, if $\pi _{\scriptop{f}}$ 
and $\pi _{\scriptop{h}}$ is the projection form 
$M_{\scriptop{f}}\times M_{\scriptop{h}}$ onto the first and the 
second factor, respectively, and the symplectic form 
$\sigma _{\scriptop{f}}$ on $M_{\scriptop{f}}$ is given by 
\begin{equation}
(\sigma _{\scriptop{f}})_{b}(\delta b,\,\delta 'b)=
\sigma ^{\got{t}}(\delta t,\,\delta 't)
+\delta\zeta ({Z'}_{\got{l}})-\delta '\zeta (Z_{\got{l}}).
\label{sigma2}
\end{equation}
Here $b=(t,\,\zeta )\,\iota (P)\in (G_{\scriptop{f}}/\iota (P))$, the tangent vectors 
$\delta b=(\delta t,\,\delta\zeta )$ 
and $\delta 'b=(\delta 't,\,\delta '\zeta )$ 
are elements of $\got{t}_{\scriptop{f}}\times N$, the vectors 
$X:=\delta t+c(\delta\zeta ,\,\zeta )/2$ and 
$X':=\delta 't+c(\delta '\zeta ,\,\zeta )/2$ 
are elements of $\got{t}_{\scriptop{f}}$, and finally 
$\sigma (L_{\delta\zeta},\, 
L_{\delta '\zeta})=0$ 
because in vi) in Proposition \ref{liftprop} 
we have $c_{\scriptop{h}}=0$. 
It follows that  $(M,\,\sigma ,\, T)$ is $T$\--equivariantly 
symplectomorphic to 
$(M_{\scriptop{f}},\,\sigma _{\scriptop{f}},\, T_{\scriptop{f}})
\times (M_{\scriptop{h}},\,\sigma _{\scriptop{h}},\, 
T_{\scriptop{h}})$, in which  
$(M_{\scriptop{f}},\,\sigma _{\scriptop{f}},\, T_{\scriptop{f}})$ is a 
compact connected symplectic manifold with a free 
symplectic action $T_{\scriptop{f}}$\--action. 
\end{proof}

\begin{remark}
In the proof of Lemma \ref{e0rem} we have also given 
a global model for 
$(M,\,\sigma ,\, T)$ in the case that the action 
of $T$ is free. 
Note that the mapping 
$g\mapsto g\, H^o$ defines an isomorphism from 
$G_{\scriptop{f}}$ onto the group $G/H^o$ in 
Remark \ref{Horem}, and an isomorphism from 
$\iota (P)$ onto $H/H^o$, which leads to an identification 
of $M_{\scriptop{f}}=G_{\scriptop{f}}/\iota (P)$ with 
the manifold $G/H$. In Remark \ref{Horem}, 
$G/H$ has been described as a principal 
$T_{\scriptop{f}}$\--bundle over the torus $N/P$, 
and as a nilmanifold for a two\--step nilpotent 
Lie group. 

The first example of Thurston \cite{thurston} is 
equal to $G_{\scriptop{f}}/\iota (P)$ with $T_{\scriptop{f}}=\R ^2/\Z ^2$, 
$N=\R ^2$, $P=\Z ^2$, $\sigma ^{\got{t}}=0$, $s=0$, 
and $c(e_1,\, e_2)=e_1$ if $e_1$, $e_2$ denotes the standard basis in $\R ^2$. 
For more examples, see McDuff and Salamon \cite[Ex. 3.8 on p.88]{MS} 
and the references therein.  
\label{freerem}
\end{remark}

\begin{remark}
If $\{ 1\}\neq T_{\scriptop{h}}\neq T$, 
then the choice of a complementary torus $T_{\scriptop{f}}$ 
to $T_{\scriptop{h}}$ in 
$T$ is far from unique, see Remark \ref{tbaserem}. It can happen 
that for some choice of $T_{\scriptop{f}}$ we have 
$c(N\times N)\subset\got{t}_{\scriptop{f}}$, whereas for 
another choice we have not. 

However, if 
$c(N\times N)\subset\got{t}_{\scriptop{h}}$ and $c\neq 0$, 
then there is no choice of a complementary torus 
$T_{\scriptop{f}}$ to $T_{\scriptop{h}}$ such that 
$c(N\times N)\subset\got{t}_{\scriptop{f}}$, and therefore 
$(M,\,\sigma ,\, T)$ is in no way $T$\--equivariantly 
symplectomorphic to a Cartesian product of a 
symplectic manifold with a free torus action 
and a Delzant manifold. 
\label{makee0rem}
\end{remark}

\begin{remark}
If $\op{dim}N\leq 1$ then $c(\zeta ,\,\zeta ')\equiv 0$ because 
every antisymmetric bilinear form on a one\--dimensional space 
is equal to zero, and we conclude that 
$(M,\,\sigma ,\, T)$ is a Cartesian product of 
a Delzant manifold with a two\--dimensional homogeneous symplectic torus. 

If $M$ is four\--dimensional, that is, $n=2$, then we have 
only few possibilities. 
\begin{itemize}
\item[a)] The homogeneous symplectic torus, where 
$T$ acts freely and transitively on $M$. 
Cases where $T$ is replaced by a subtorus with coisotropic 
orbits are treated as subcases. 
\item[b)] $(M,\,\sigma ,\, T)$ is 
$T$\--equivariantly symplectomorphic with the Cartesian product 
of a two\--di\-men\-si\-onal homogeneous symplectic torus 
and a sphere, provided with a rotationally invariant area form. 
\item[c)] $(M,\,\sigma ,T)$ is a four\--dimensional 
Delzant manifold. 
\item[d)] The action of $T$ is free 
with Lagrangian orbits, but not in case a). 
See the proof of Lemma \ref{e0rem} for a more 
detailed description of this case. The example of 
Thurston mentioned in Remark \ref{freerem} was the first 
one in the literature. 
\end{itemize}
\label{fourdimrem}
\end{remark}

\subsection{The holonomy invariant} 
\label{holinvariantss}
In view of the surjectivity of the mapping $A$ in 
Proposition \ref{Aprop}, Lemma \ref{TPlem} contains 
the description of the dependence of the 
holonomy $\tau (x):\zeta\mapsto\tau _{\zeta}(x):P\to T$ 
on all points $x\in M$. The only change which occurs is 
that if $x$ is replaced by $x'=\op{e}^{L_{\zeta '}}(x)$, 
$\zeta '\in N$, then $\tau _{\zeta}(x)$ is replaced by 
$\tau _{\zeta}(x)\,\op{e}^{c(\zeta ,\,\zeta ')}$, see 
(\ref{tauzeta}). 
We now investigate  
the dependence of the holonomy 
on the choice of the admissible connection 
as in Proposition \ref{liftprop}. 

It follows from Lemma 
\ref{liftfreedomlem} that $\got{l}^*\ni\xi\mapsto\widetilde{L}_{\xi}$ 
is another connection as in Proposition \ref{liftprop}, 
if and only if there exists a smooth $T$\--invariant 
mapping $\alpha :x\mapsto (\xi\mapsto\alpha _{\xi}(x))$ 
from $M$ to $\op{Lin}(\got{l}^*,\,\got{l})$, 
which is {\em closed} when viewed as an $\got{l}$\--valued 
one\--form on $M/T$, and {\em symmetric}  
in the sense of (\ref{alphasym}), such that 
$\widetilde{L}_{\xi}(x)=
L_{\xi}(x)+\alpha_{\xi}(x)_M(x)$ for every 
$x\in M$ and $\xi\in\got{l}^*$. 
The change from $L$ to $\widetilde{L}$ leads to a 
change from $\tau_{\zeta}(x)$ to 
\begin{equation}
\widetilde{\tau}_{\zeta}(x)
=\tau_{\zeta}(x)\,  \op{e}^{\int_{\gamma_{\zeta}}\,\alpha}.
\label{tildetau}
\end{equation} 
Here $\gamma_{\zeta}(t):=\pi (x)+t\,\zeta$, 
$0\leq t\leq 1$, is a loop in $M/T$ because $\zeta\in P$. 
Because $\alpha$ is closed, the integral 
$\int_{\gamma_{\zeta}}\,\alpha$ only depends on the de Rham 
cohomology class of $\alpha$, which means that for the effect on 
the $\tau_{\zeta}(x)$'s we can restrict ourselves to constant 
$\got{l}$-valued one-forms on $\got{l}^*$, that is, 
linear mappings from $\got{l}^*$ to $\got{l}$, or equivalently, 
bilinear forms on $\got{l}^*$. Therefore, at a given 
point $x\in M$, the allowed changes in $\tau_{\zeta}(x)$, 
$\zeta\in P$,  
consist of the multiplications with $\op{e}^{\alpha_{\zeta}}$, 
where $\alpha$ ranges over the space of 
linear mappings $\xi\mapsto\alpha _{\xi}$ 
from $\got{l}^*$ to $\got{l}$, which are 
symmetric in the sense of (\ref{alphasym}). 

\begin{definition}
Let $\op{Hom}_c(P,\, T)$ denote the space of mappings 
$\tau :\zeta\mapsto \tau_{\zeta}:P\to T$ such that 
\begin{equation}
\tau _{\zeta '}\,\tau _{\zeta}=\tau _{\zeta +\zeta '}
\,\op{e}^{c(\zeta ',\,\zeta )/2},\quad \zeta,\,\zeta '\in P. 
\label{tauplushom}
\end{equation} 
Because $c(P\times P)\subset T_{\Z}$, 
the factor $\op{e}^{c(\zeta ',\zeta )/2}$ in (\ref{tauplushom}) 
is an element of order two in $T$. Therefore 
the elements of $\op{Hom}_c(P,\, T)$ are quite close 
to being homomorphisms from $P$ to $T$. They 
are homomorphisms from $P$ to $T$ if 
$c(P\times P)\subset 2T_{\Z}$. 
\label{Homcdef}
\end{definition}

If $h:\zeta\mapsto h_{\zeta}$ is a homomorphism from $P$ to $T$, 
then $h\cdot\tau :\zeta\mapsto \tau _{\zeta}\, h_{\zeta} 
\in\op{Hom}_c(P,\, T)$ for every $\tau\in\op{Hom}_c(P,\, T)$, 
and $(h,\,\tau )\mapsto h\cdot\tau$ defines a 
free, proper, and transitive action of $\op{Hom}(P,\, T)$ 
on $\op{Hom}_c(P,\, T)$. If $\varepsilon ^l$, 
$1\leq l\leq d_N:=\op{dim}N$ is a $\Z$\--basis of $P$, 
then the mapping $h\mapsto (h_{\varepsilon ^l})_{1\leq l\leq d_N}$ 
is an isomorphism from $\op{Hom}(P,\, T)$ onto the torus 
$T^{d_N}$. Therefore $\op{Hom}_c(P,\, T)$ is diffeomorphic 
to a torus of dimension $\op{dim}N\,\op{dim}T$. 

For each $\zeta '\in N$, $c(\cdot ,\,\zeta '):\zeta\mapsto 
c(\zeta ,\,\zeta ')$ is a homomorphism from $P$ to 
$\got{t}$, actually $\got{l}$\--valued. Write 
$c(\cdot ,\, N)$ for the set of all $c(\cdot ,\,\zeta ')
\in\op{Hom}(P,\,\got{t})$ such that $\zeta '\in N$. 
$c(\cdot ,\, N)$ is a linear subspace of the Lie 
algebra $\op{Hom}(P,\,\got{t})$ of $\op{Hom}(P,\, T)$. 

\begin{definition}
Let $\op{Sym}$ denote the space of all linear mappings 
$\alpha :\got{l}^*\to\got{l}$ which are symmetric in the 
sense of (\ref{alphasym}). For each $\alpha\in\op{Sym}$, 
the restriction $\alpha |_P$ of $\alpha$ to $P$ is a 
homomorphism from $P$ to $\got{l}\subset\got{t}$. 
In this way the set $\op{Sym}|_P$ of all $\alpha |_P$ 
such that $\alpha\in\op{Sym}$ is another linear 
subspace of $\op{Hom}(P,\,\got{t})$. Write 
\begin{equation}
{\cal T}:=\op{Hom}_c(P,\, T)/\op{exp}{\cal A},
\quad {\cal A}:=c(\cdot ,\, N)+\op{Sym}|_P
\label{calT}
\end{equation} 
for the orbit space of the action of  
the Lie subgroup $\op{exp}{\cal A}$ 
of $\op{Hom}(P,\, T)$ on $\op{Hom}_c(P,\, T)$. 
Because $\op{exp}{\cal A}$ 
need not be a closed subgroup of $\op{Hom}(P,\, T)$, 
the quotient topology of ${\cal T}$ need not be Hausdorff. 

It follows from (\ref{tauplus}) and (\ref{tauplushom}), 
that for every choice of a connection as in Proposition 
\ref{liftprop} and every $x\in M$, the mapping 
$\tau (x):\zeta\mapsto\tau _{\zeta}(x)$ is an element 
of $\op{Hom}_c(P,\, T)$. It is the point of 
(\ref{calT}), that the right hand side in  
\begin{equation}
\overline{\tau}:=
(\op{exp}{\cal A})
\cdot\tau (x)\in {\cal T}
\label{overlinetau}
\end{equation}
defines an invariant $\overline{\tau}$ of our symplectic $T$\--space 
$(M,\,\sigma ,\, T)$, in the sense that it neither depends 
on the choice of the point $x\in M$, nor on the 
choice of the connection as in Proposition \ref{liftprop}. 
\label{overlinetaudef}
\end{definition} 
\begin{remark}
In order to obtain some more insight in the vector 
space ${\cal A}$ in (\ref{calT}), we use the direct sum 
decomposition $\got{t}=\got{t}_{\scriptop{h}}\oplus 
\got{t}_{\scriptop{f}}$, where 
$\got{t}_{\scriptop{h}}\subset\got{l}$ is the 
Lie algebra of the Hamiltonian torus $T_{\scriptop{h}}$ 
and $\got{t}_{\scriptop{f}}$ is the Lie algebra of 
a complementary torus $T_{\scriptop{f}}$ to $T_{\scriptop{h}}$ 
in $T$. This leads to an identification of 
$N=(\got{l}/\got{t}_{\scriptop{h}})^*$ with 
$(\got{l}\cap\got{t}_{\scriptop{f}})^*$ and 
of its linear complement $C$ in $\got{l}^*$ 
with ${\got{t}_{\scriptop{h}}}^*$. 

Let $(\op{Sym}_{\scriptop{f}})|_P$ denote the 
space of all linear mappings 
$\alpha :(\got{l}\cap\got{t}_{\scriptop{f}})^*\to 
\got{l}\cap\got{t}_{\scriptop{f}}$, which satisfy 
the symmetry condition (\ref{alphasym}) with 
$\got{l}$ replaced by $\got{l}\cap\got{t}_{\scriptop{f}}$. 
The space $\op{Sym}|_P$ of all restrictions to 
$P\subset N=(\got{l}\cap\got{t}_{\scriptop{f}})^*$ 
of linear mappings 
$\alpha :\got{l}^*\to\got{l}$ which satisfy the 
the symmetry condition (\ref{alphasym}) 
is equal to the direct sum of the space 
$\op{Hom}(P,\,\got{t}_{\scriptop{h}})$ of {\em all} 
homomorphisms from $P$ to $\got{t}_{\scriptop{h}}$, 
and the space $(\op{Sym}_{\scriptop{f}})|_P$. 
This means that in the space ${\cal T}$ in 
(\ref{calT}) we dispose of the $T_{\scriptop{h}}$\--
components, and in the computation of 
${\cal A}$ we can replace $c$ by its 
$\got{l}\cap\got{t}_{\scriptop{f}}$\--component $c_{\scriptop{f}}$. 

Now suppose that $\zeta '\in N$ and 
$c_{\scriptop{f}}(\cdot ,\,\zeta ')
\in (\op{Sym}_{\scriptop{f}})|_P$. 
This is equivalent to the condition that 
\[
-\zeta '(c(\zeta ,\,\zeta ''))
=\zeta ''(c(\zeta ,\,\zeta '))+\zeta (c(\zeta ',\zeta ''))
=\zeta ''(c_{\scriptop{f}}(\zeta ,\,\zeta '))
-\zeta (c_{\scriptop{f}}(\zeta '',\zeta '))=0
\]
for all $\zeta ,\,\zeta ''\in P\subset N
=(\got{l}\cap\got{t}_{\scriptop{f}})^*$. Here we have 
used (\ref{ft}) in the 
first equality. In the second equality 
we have used the antisymmetry of $c$ and the fact 
that the elements $\zeta '',\,\zeta 
\in N=(\got{l}/\got{t}_{\scriptop{h}})^*$ are equal to 
zero on $\got{t}_{\scriptop{h}}$. 
In other words, $c_{\scriptop{f}}(\cdot ,\,\zeta ')
\in (\op{Sym}_{\scriptop{f}})|_P$ if and only 
if $\zeta '=0$ on $c(P\times P)$, or equivalently 
$\zeta '=0$ on the linear subspace of $\got{l}$ which 
is spanned by $c(N\times N)$. Let $c^0$ denote the 
space of all $\zeta '\in\got{l}^*$ which are 
equal to zero on the linear span of $c(N\times N)$. 
In view of (\ref{ft}) we have $\op{ker}c\subset c^0$, 
and it follows that the dimension of $c_{\scriptop{f}}(\cdot,\, N)
\cap (\op{Sym}_{\scriptop{f}})|_P$ is equal to 
$\op{dim}c^0-\op{dim}\op{ker}c$, whereas the dimension of 
$c_{\scriptop{f}}(\cdot,\, N)$ is equal to 
$\op{dim}N-\op{dim}\op{ker}c_{\scriptop{f}}$. 
It follows that the dimension of 
$c_{\scriptop{f}}(\cdot,\, N)+(\op{Sym}_{\scriptop{f}})|_P$
is equal to 
$d_N(d_N+1)/2-(\op{dim}c^0-\op{dim}\,\op{ker}c)$, 
in which $d_N=\op{dim}N=\op{dim}(\got{l}\cap\got{t}_{\scriptop{f}})$. 
Therefore the codimension of 
$c_{\scriptop{f}}(\cdot,\, N)+(\op{Sym}_{\scriptop{f}})|_P$ 
in the ${d_N}^2$\--dimensional space $\op{Hom}(P,\,
\got{l}\cap\got{t}_{\scriptop{f}})$ is equal to 
$d_N(d_N-3)/2
+\op{dim}\,\op{ker}c_{\scriptop{f}}
-\op{dim}\,\op{ker}c+\op{dim}c^0$.
Because all elements of ${\cal A}$ map to $\got{l}$, it follows that 
\begin{equation}
\op{dim}{\cal T}=d_N(\op{dim}T-d_N)+d_N(d_N-3)/2
+\op{dim}\,\op{ker}c_{\scriptop{f}}
-\op{dim}\,\op{ker}c+\op{dim}c^0.
\label{dimcalT}
\end{equation}
\label{dimholinvariantrem}
\end{remark}

\section{Applications of the global model}
\label{applsec}
In this section, which is not 
needed for the classification in Section \ref{invariantsec}, 
we give some applications of Proposition \ref{Aprop} 
to minimal coupling, the reduced phase spaces, 
the topology of the torus action, and 
to the universal covering of our symplectic $T$\--space $M$. 

\subsection{Minimal coupling}
\label{minimalcouplingss}
The fibration of $M$ by Delzant manifolds is a fibration by symplectic
submanifolds with a structure group $H$ which acts on the fiber 
by means of symplectomorphisms. See Proposition \ref{Aprop}. 
Moreover, the distribution  
spanned by the $Z_M$, $Z\in\got{t}_{\scriptop{f}}$, and 
the $L_{\zeta}$, $\zeta\in N$, which we used in the construction 
of the model, is the symplectic orthogonal complement of the 
fibers. This follows from the fact that at the regular points 
the tangent space to the Delzant submanifold is 
spanned by the $Y_M$, $Y\in\got{t}_{\scriptop{h}}$, 
and the $L_{\eta}$, $\eta\in C$, combined with 
Lemma \ref{constlem} and $\got{t}_{\scriptop{h}}\subset
\got{l}:=\op{ker}\sigma ^{\got{t}}$, 
the equation (\ref{stronglift}), and v) in Proposition \ref{liftprop}. 
Because the $Z_M$, $Z\in\got{t}_{\scriptop{f}}$, 
commute with each other and with the $L_{\zeta}$, 
$\zeta\in N$, the only nonzero Lie brackets of 
horizontal vector fields are the $[L_{\zeta},\, L_{\zeta}]=
c(\zeta ,\,\zeta ')_M$, $\zeta ,\,\zeta '\in N$, see 
iii) in Proposition \ref{liftprop}. The vertical part 
of $[L_{\zeta},\, L_{\zeta}]$ is equal to 
$c_{\scriptop{h}}(\zeta ,\,\zeta ')_M$. 
Because, for every $Y\in\got{t}_{\scriptop{h}}$,  
$Y_M$ is the Hamiltonian vector field defined by 
the function $x\mapsto \mu (x)(Y)$, 
the vertical part of $[L_{\zeta},\, L_{\zeta}]$ is 
the Hamiltonian vector field defined by the 
function $x\mapsto \mu (x)(c_{\scriptop{h}}(\zeta ,\,\zeta '))$. 
It follows from vi) in Proposition \ref{liftprop} that 
the derivative of this function is equal to the negative of  
the derivative of $\sigma (L_{\zeta},\, L_{\zeta '})$, 
and therefore the vertical part of $[L_{\zeta},\, L_{\zeta '}]$ 
is equal to $-\op{Ham}_{\sigma (L_{\zeta},\, L_{\zeta '})}$. 
This equation, which holds in great generality 
for the curvature of the symplectically 
orthogonal connection in a fibration by symplectic manifolds,  
is known as {\em minimal coupling}, see 
Guillemin, Lerman and Sternberg \cite[Sec. 1.3]{gls}. 
In this way equation vi) in Proposition \ref{liftprop} 
represents the minimal coupling term in the 
symplectic form on $M$. 
This observation was 
suggested to us by Yael Karshon. 

Recall that the fibration of $M$ by Delzant manifolds 
was not a priori given. It has been constructed using 
the special admissible connection introduced in 
Proposition \ref{liftprop}, and it is not unique 
if $\{ 1\}\neq T_{\scriptop{h}}\neq T$.    

\subsection{The reduced phase spaces}
\label{redG/Hrem}
On the symplectic manifold $(M,\,\sigma )$ we have the Hamiltonian 
action of the torus $T_{\scriptop{h}}$, with momentum mapping 
$\mu :M\to {\got{t}_{\scriptop{h}}}^*$, where $\mu (M)\simeq\Delta$. 
Let $q\in\mu (M)$. Then, restricting the discussion 
to the orbit type stratum which contains $\mu ^{-1}(\{ q\})$, 
we obtain that $\mu ^{-1}(\{ q\})$ is a compact and connected 
smooth submanifold of $M$, on which $T_{\scriptop{h}}/H$ acts 
freely, where $H$ denotes the common stabilizer subgroup 
of the elements in $\mu ^{-1}(\{ q\})$. It follows that 
the orbit space $M^q:=\mu ^{-1}(\{ q\})/T_{\scriptop{h}}$ has a unique 
structure of a compact connected smooth manifold, such that 
the projection $\pi ^q:\mu ^{-1}(\{ q\})\to M^q$ 
is a principal $T_{\scriptop{h}}/H$\--fibration. 

At each point of $\mu ^{-1}(\{ q\})$, the kernel of the 
pull\--back to $\mu ^{-1}(\{ q\})$ of $\sigma$ is equal to 
the tangent space of the $T_{\scriptop{h}}$\--orbit through that point, 
and it follows that there is a unique symplectic form 
$\sigma ^q$ on $M^q$ such that 
$(\pi ^q)^*\,\sigma ^q=(\iota ^q)^*\,\sigma$, if 
$\iota ^q$ denotes the inclusion mapping from 
$\mu ^{-1}(\{ q\})$ to $M$. The symplectic manifold 
$(M^q,\,\sigma ^q)$ is called the {\em 
reduced phase space at the $\mu$\--value $q$} for the 
Hamiltonian action of $T_{\scriptop{h}}$ on $(M,\,\sigma )$.  

On $M^q$ we still have the action of the torus $T/T_{\scriptop{h}}$, 
which is free, leaves the symplectic form $\sigma ^q$ 
invariant, and has coisotropic orbits. The vector fields 
$L_{\zeta}$, $\zeta\in N$, are tangent to $\mu ^{-1}(\{ q\})$, 
and are intertwined by $\pi ^q$ with unique smooth 
vector fields $L_{\zeta}^q$ on $M^q$. In combination 
with  $(\pi ^q)^*\,\sigma ^q=(\iota ^q)^*\,\sigma$, this implies that 
$(\pi ^q)^* (\sigma ^q(L_{\zeta}^q,\, L_{\zeta '}^q))
=(\iota ^q)^*(\sigma (L_{\zeta},\, L_{\zeta '}))$, 
as an identity between constant functions on $\mu ^{-1}(\{ q\})$. 
It therefore follows from vi) in Proposition \ref{liftprop} 
that 
\begin{equation}
\sigma ^q(L_{\zeta}^q,\, L_{\zeta '}^q)
=-q(c_{\scriptop{h}}(\zeta ,\,\zeta ')),\quad 
\zeta ,\zeta '\in N.
\label{sigmaq}
\end{equation}

We now show that each of the reduced phase spaces $M^q=
\mu ^{-1}(\{ q\})/T_{\scriptop{h}}$ can be identified with the 
$G$\--homogeneous space $G/H\simeq ((T/T_{\scriptop{h}})\times N)
/\iota (P)$ discussed in 
Remark \ref{Horem}. Moreover, if 
$c(N\times N)\subset\got{t}_{\scriptop{h}}$, 
then (\ref{sigmaq}) corresponds to 
the description of the variation of the cohomology 
class of the symplectic form of the reduced phase spaces 
in Duistermaat and Heckman \cite{dh}. 

Let $x\in\mu ^{-1}(\{ q\})$, 
and write 
\[
H_x=\{ (t,\,\zeta )\in T\times P\mid t\,\tau _{\zeta}\in T_x\} .
\]
Because $T_x$ is a closed Lie subgroup of $T_{\scriptop{h}}$, 
$H_x$ is a closed Lie subgroup of $H$, see 
(\ref{Hdef}), and $G/H_x$ is a 
compact $G$\--homogeneous space. The mapping 
$A_x:(t,\,\zeta )\mapsto t\cdot\op{e}^{L_{\zeta}}(x):
G\to M$ induces an embedding $\alpha _x$ from 
$G/H_x$ into $M$, with image equal to $\mu ^{-1}(\{ q\} )$. 
This exhibits $\mu ^{-1}(\{ q\} )$ as a compact 
and connected smooth submanifold of $M$, and actually 
as a $G$\--homogeneous space. The pull\--back to 
$G/H_x$ of the symplectic form $\sigma$ is given by the formula 
(\ref{A*sigma}), in which $\delta x=\delta 'x =0$. 

Because $T_{\scriptop{h}}/T_x\simeq H/H_x$, the mapping $\alpha _x$ 
induces a $T/T_{\scriptop{h}}$\--equivariant diffeomorphism 
$\beta _x$ from $G/H=(G/H_x)/(H/H_x)$ onto the reduced phase space  
$M^q=\mu ^{-1}(\{ q\} )/(T_{\scriptop{h}}/T_x)$. Because the 
dimension of $\mu ^{-1}(\{ q\} )$ jumps down if 
$q\in\Delta$ moves into a lower\--dimensional orbit type stratum, 
it is quite remarkable that nevertheless the reduced phase 
spaces $M^q$ for all $q\in\Delta$ are isomorphic to the 
same space $G/H$ in a natural way. In this model 
the principal $T_{\scriptop{h}}/T_x$\--fibration $\mu ^{-1}(\{ q\})\to M^q$ 
corresponds to the principal $H/H_x$\--fibration 
$G/H_x\to G/H$, in which $T_{\scriptop{h}}/T_x\simeq H/H_x$ 
is a torus. 

If $c(N\times N)\subset \got{t}_{\scriptop{h}}$, then the Chern class of 
the principal $H/H_x$\--fibration 
$\pi _x:G/H_x\to G/H$, which is an element of 
$\op{H}^2(G/H,\, (H/H_x)_{\Z})$ is equal to 
$\psi ^*\, c$, in which $c\in\op{H}^2(N/P,\, T_{\Z})$ is 
the cohomology class corresponding to the antisymmetric 
bilinear form $c$ introduced in Proposition \ref{liftprop}, 
and $\psi$ is the projection from $G/H$ onto $(G/H)/(T/T_{\scriptop{h}})
\simeq N/P$. Therefore, in the case that 
$c(N\times N)\subset \got{t}_{\scriptop{h}}$, formula 
(\ref{sigmaq}) shows that  
the variation of the cohomology 
class of the symplectic form of the reduced phase spaces 
is equal to the cohomology class $-c$ of the curvature form. 

\subsection{The $T_{\scriptop{h}}$\--fixed point set modulo 
$T_{\scriptop{f}}$}
\label{fixedpointss}
The action of the Hamiltonian torus 
$T_{\scriptop{f}}$ on $M$ has fixed points, which 
are the $x\in M$ such that $\mu (x)$ is equal to a 
vertex $v$ of the Delzant polytope $\Delta$, 
if $\mu :M\to\Delta\subset {\got{t}_{\scriptop{h}}}^*$ 
denotes the momentum mapping of the Hamiltonian 
$T_{\scriptop{h}}$\--action as in (\ref{pi^1}). 
Let $v$ be a vertex of $\Delta$. 
Because $T_x=T_{\scriptop{h}}$ for every 
$x\in\mu ^{-1}(\{ v\})$, the 
reduced phase space $\mu ^{-1}(\{ v\})/T_{\scriptop{h}}$ 
at the level $v$, introduced in Subsection \ref{redG/Hrem}, 
is equal to $\mu ^{-1}(\{ v\})$. Because the reduced phase 
spaces are connected, the $\mu ^{-1}(\{ v\})$, where 
$v$ ranges over the vertices of $\Delta$, are the 
connected components ${\cal F}$ of the 
fixed point set $M^{T_{\scriptop{h}}}$ of the 
$T_{\scriptop{h}}$\--action in $M$. Because in this subsection 
we want to find invariants of the $T$\--action, 
disregarding the symplectic structure, we use the 
notation ${\cal F}$ for the connected components 
of $M^{T_{\scriptop{h}}}$, instead of $\mu ^{-1}(\{ v\})$. 

Note that each ${\cal F}$ is a global section 
of the fibration $\delta :M\to G/H\simeq 
((T/T_{\scriptop{h}})\times N)
/\iota (P)$ of $M$ by Delzant 
submanifolds. Using Morse theory with the Hamiltonian 
functions of infinitesimal $T_{\scriptop{h}}$\--actions as Bott\--Morse 
functions, this may lead to useful information 
about the topology of $M$ in terms of the connected components  
${\cal F}$ of $M^{T_{\scriptop{h}}}$. 

Let $T_{\scriptop{f}}$ be a complementary torus 
to $T_{\scriptop{h}}$ in $T$. Because the action 
of $T_{\scriptop{f}}$ is free, we have the principal 
$T_{\scriptop{f}}$\--fibration $M\to M/T_{\scriptop{f}}$, 
and because the actions of $T_{\scriptop{f}}$ and 
$T_{\scriptop{h}}$ commute, we have an induced 
action of $T_{\scriptop{h}}$ on $M/T_{\scriptop{f}}$. 
The manifolds ${\cal F}/T_{\scriptop{h}}$ are the connected 
components of the fixed point set 
$(M/T_{\scriptop{f}})^{T_{\scriptop{h}}}$ of the 
$T_{\scriptop{f}}$\--action in $M/T_{\scriptop{f}}$. 
The fibration $\delta :M\to G/H\simeq 
((T/T_{\scriptop{h}})\times N)
/\iota (P)$ induces a fibration 
\[
\overline{M}:=M/T_{\scriptop{f}}\to (G/H)/T_{\scriptop{f}}\simeq 
(T_{\scriptop{f}}\times N)/\iota (P))/T_{\scriptop{f}}
\simeq N/P
\]
by Delzant manifolds, of which each connected 
component $\overline{\cal F}:={\cal F}/T_{\scriptop{f}}$ 
of the $T_{\scriptop{h}}$\--fixed 
point set is a global section, diffeomorphic to $N/P$.  

Let $x\in {\cal F}$, and write $y=T_{\scriptop{f}}\cdot x
\in\overline{\cal F}$. The tangent action 
of $T_{\scriptop{h}}$ on the normal space 
$N_y:=\op{T}_y\overline{M}/\op{T}_y\! \overline{\cal F}$ 
to $\overline{\cal F}$ can be identified with the tangent action 
of $T_{\scriptop{h}}$ 
on the tangent space of the Delzant manifold through $x$. 
It follows from the local model in Lemma \ref{modellem}, 
with $H=T_x=T_{\scriptop{h}}$,  
$\got{h}=\got{t}_{\scriptop{h}}$, 
and $m=d_{\scriptop{h}}:=\op{dim}T_{\scriptop{h}}$, 
that $N_y$ has a direct sum decomposition into 
$T_{\scriptop{h}}$\--invariant two\--dimensional linear 
subspaces $E_y^j$, $1\leq j\leq d_{\scriptop{h}}$, a complex structure 
on each $E_y^j$, and a corresponding $\Z$\--basis 
$Y_j$, $1\leq j\leq d_{\scriptop{h}}$ of the integral lattice 
$(T_{\scriptop{h}})_{\Z}$ 
in $\got{t}_{\scriptop{h}}$, such that 
the tangent action of $\op{e}^Y$, $Y\in\got{t}_{\scriptop{h}}$ on 
$N_y$ corresponds to 
the multiplication with $\op{e}^{2\pi\fop{i}Y^j}$ in 
$E_y^j$, if $Y=\sum_{j=1}^{n_1}\, Y^j\, Y_j$. 
Although for their existence we referred to the local model in 
Lemma \ref{modellem} for our symplectic $T$\--space, 
all these ingredients are uniquely determined in terms of 
the linearized action of $T_{\scriptop{h}}$ 
on the normal bundle $N$ of $\overline{\cal F}
:={\cal F}/T_{\scriptop{f}}$ in $\overline{M}:=M/T_{\scriptop{f}}$,  
up to a permutation of the indices $j$. That is, 
disregarding the symplectic structure. 

For each $j$, the $E_y^j$, $y\in\overline{\cal F}$, 
form a complex line bundle $E^j$ over 
$\overline{\cal F}\simeq N/P$,  
and the normal bundle $N$ 
of $\overline{\cal F}$ in $\overline{M}$ 
is the direct sum of the complex 
line bundles $E^j$, $1\leq j\leq d_{\scriptop{h}}$. 

Any smooth complex line bundle $L$ over a smooth manifold 
$B$ has a Chern class, 
which is defined as follows. Let $\C ^{\times}$ denote 
the multiplicative group of the nonzero complex 
numbers. The transition functions of local trivializations 
define a 1\--cocycle of germs of smooth $\C ^{\times}$\--valued 
functions, and the bundle $L$ is classified by the 
sheaf (= \v{C}ech) cohomology class $\gamma\in\op{H}^1(B,\,
\op{C}^{\infty}(\cdot ,\,\C ^{\times}))$ of the 
1-cocycle of the transition functions. 
Because the sheaf $\op{C}^{\infty}(\cdot ,\,\C )$ 
is fine, the short 
exact sequence 
\[
0\to\Z\to\op{C}^{\infty}(\cdot ,\,\C )
\stackrel{\fop{e}^{2\pi\scriptop{i}}}{\to}
\op{C}^{\infty}(\cdot ,\,\C ^{\times})\to 1
\]
induces an isomorphism 
$\delta :\op{H}^1(B,\,\op{C}^{\infty}(\cdot ,\,\C ^{\times})
\to\op{H}^2(B,\,\Z )$, and the 
{\em Chern class of the complex line bundle  
$L$ over $B$} is defined as the cohomology class 
$\op{c}(L):=\delta (\gamma )\in\op{H}^2(B,\,\Z )$. 
With these definitions, we have the following conclusions. 
\begin{proposition}
Let the $c_{\scriptop{h}}^j\in\Lambda ^2N^*$, 
$1\leq j\leq d_{\scriptop{h}}$, be defined by 
\[
c_{\scriptop{h}}(\zeta ,\,\zeta ')
=\sum_{j=1}^{d_{\scriptop{h}}}\, 
c_{\scriptop{h}}^j(\zeta ,\,\zeta ')\, Y_j,
\quad \zeta ,\,\zeta '\in N.
\]
Let $1\leq j\leq d_{\scriptop{h}}$. 
Viewing $c_{\scriptop{h}}^j$ as an element of 
$\op{H}^2(N/P,\,\R )\simeq\op{H}^2(\overline{\cal F},\,\R )$ 
as in Corollary \ref{cohomcor}, we have that 
$c_{\scriptop{h}}^j$ is equal to the image 
in ~$\op{H}^2(\overline{\cal F},\,\R )$ 
under the coefficient homomorphism 
$\op{H}^2(\overline{\cal F},\,\Z)
\to \op{H}^2(\overline{\cal F},\,\R )$, of the 
Chern class ~$\op{c}(E^j)$ 
of the com\-plex line bundle $E^j$ over $\overline{\cal F}
\simeq N/P$. 

If $M$ is $T$\--equivariantly diffeomorphic to 
$M_{\scriptop{f}}\times M_{\scriptop{h}}$, 
in which $T_{\scriptop{h}}$ acts only on $M_{\scriptop{h}}$ 
with isolated fixed points, 
and $T_{\scriptop{f}}$ acts only on $M_{\scriptop{f}}$ 
and freely, then 
$c_{\scriptop{h}}=0$, and we have the conclusions 
c) and a) in Lemma \ref{e0rem}.  
\label{chernprop}
\end{proposition}
\begin{proof}
Because the vector fields $L_{\zeta}$, $\zeta\in N$, are 
invariant under the action of $T$, hence under the action 
of $T_{\scriptop{h}}$ and $T_{\scriptop{f}}$, they 
are intertwined by the projection $M\to \overline{M}:=M/T_{\scriptop{f}}$ 
to uniquely determined $T_{\scriptop{h}}$\--invariant 
smooth vector fields on $\overline{M}$, which we 
also denote by $L_{\zeta}$. The identity iii) in 
Proposition \ref{liftprop} leads to the identity 
$[L_{\zeta},\, L_{\zeta '}]=
c_{\scriptop{h}}(\zeta ,\,\zeta ')_{\overline{M}}$ 
for vector fields on $\overline{M}$. 
Because the $L_{\zeta}$ are $T_{\scriptop{h}}$\--invariant, 
their flows leave each connected component $\overline{\cal F}$ 
of the $T_{\scriptop{h}}$\--fixed point set 
$\overline{M}^{T_{\scriptop{h}}}$ invariant, 
and their linearizations define automorphisms 
of the normal bundle $N$ of $\overline{\cal F}$ in 
$\overline{M}$ which commute with the linearized 
action of $T_{\scriptop{h}}$ on $N$. Therefore 
these automorphisms leave each of the complex line 
bundles $E^j$ invariant, and the corresponding  
infinitesimal automorphisms define vector fields 
on $E_j$ which we again denote by $L_{\zeta}$. 
Because the $L_{\zeta}$ are lifts of the constant 
vector fields $\zeta$ on $N/P$, we conclude that 
$N\ni\zeta\mapsto L_{\zeta}$ is a 
$\T$\--invariant connection in $E^j$, where $\T$ is 
the unit circle in $\C$. 
Because the cohomology class in 
$\op{H}^2(N/P,\,\R )\simeq 
\op{H}^2(\overline{\cal F},\,\R )$ of the negative of the curvature 
form is equal to the image of $\op{c}(E^j)$ 
in $\op{H}^2(\overline{\cal F},\,\R )$ under the coefficient 
homomorphism $\op{H}^2(\overline{\cal F},\, \Z )
\to\op{H}^2(\overline{\cal F},\,\R)$, 
the first statement in the proposition follows 
from the combination of the above discussions with the 
identifications in 
Remark \ref{curvrem} and Remark \ref{integralcrem}, 
and the general facts about Chern classes of 
complex line bundles as for instance 
in Bott and Tu \cite[pp. 270, 267, 72, 73]{bt}.  

For the second statement assume that 
$M=M_{\scriptop{f}}\times M_{\scriptop{h}}$, 
in which $T_{\scriptop{h}}$ acts only on $M_{\scriptop{h}}$ 
and has isolated fixed points, 
and $T_{\scriptop{f}}$ acts freely on $M_{\scriptop{f}}$. 
Then   
$
\overline{M}
:=(M_{\scriptop{f}}\times M_{\scriptop{h}})/T_{\scriptop{f}}
=(M_{\scriptop{f}}/T_{\scriptop{f}})\times M_{\scriptop{h}},
$
in which $T_{\scriptop{h}}$ only acts on the second component. 
It follows that the connected components of 
$\overline{M}^{T_{\scriptop{h}}}$ are of the form 
$
\overline{\cal F}=(M_{\scriptop{f}}/T_{\scriptop{f}})\times \{ x\}, 
$
in which $x$ ranges over the isolated fixed 
points of the $T_{\scriptop{h}}$\--action on $M_{\scriptop{h}}$, 
and the normal bundle $N$ of $\overline{\cal F}$ in 
$\overline{M}$ is $T_{\scriptop{h}}$\--equivariantly 
isomorphic to $(M_{\scriptop{f}}/T_{\scriptop{f}})\times 
\op{T}_x\! M_{\scriptop{h}}$. This shows that 
each of the complex line bundles $E^j$ is trivial, 
which implies that $\op{c}(E^j)=0$ and therefore 
$c_{\scriptop{h}}^j=0$ in view of 
the first statement in the proposition. 
Because this holds for every $1\leq j\leq d_{\scriptop{h}}$, 
it follows that $c_{\scriptop{h}}=0$. 
\end{proof}

\subsection{A universal covering of $M$}
\label{coveringss}
In Proposition \ref{pi1prop} below, we describe 
an explicit universal covering 
of the manifold $M$ by a Cartesian product $\widetilde{M}$ of a vector space and the 
Delzant manifold $M_{\scriptop{h}}$, which leads to 
an explicit description of the fundamental group of $M$.  
In Remark \ref{benoist6.16rem} we recover Corollaire 6.16 
of Benoist \cite{benoist}, 
which states that the universal cover of a 
compact connected symplectic $T$\--manifold with 
coisotropic principal orbits is 
$(\got{t}_{\scriptop{f}}\times T_{\scriptop{h}})$\--equivariantly 
symplectomorphic to the Cartesian product of a symplectic vector space 
and a Delzant manifold.  

\medskip\noindent
Let 
$\varepsilon ^l$, $1\leq l\leq 
d_N:=\op{dim}N$, be a $\Z$\--basis of the period group $P$ in $N$. 
If $\zeta,\,\zeta '\in N$ have coordinates 
$\zeta _l,\,\zeta '_l$ with respect to this basis, 
then we write  
\begin{equation}
b(\zeta ,\,\zeta '):=\sum_{l<l'}\, \zeta _l\,\zeta '_{l'}\, 
c^{l\, l'},\quad c^{l\, l'}:=c(\varepsilon ^l,\,\varepsilon ^{l'}). 
\label{b}
\end{equation}
This defines a bilinear mapping $b:N\times N\to\got{l}$ 
such that $c(\zeta ,\,\zeta ')=b(\zeta ,\,\zeta ')-
b(\zeta ',\,\zeta )$. We have $\zeta ,\,\zeta '\in P$ 
if and only if $\zeta _l,\,\zeta '_l\in\Z$ for all $l$. Therefore  
$c(P\times P)\subset T_{\Z}$, see 
Lemma \ref{TPlem}, implies that $b(P\times P)
\subset T_{\Z}$. 

Let $x\in M_{\scriptop{h}}$. For each $1\leq l\leq d_N$ we choose 
$X^l\in\got{t}$ such that $\tau _{\varepsilon ^l}(x)=
\op{e}^{X^l}$. Then (\ref{tauzetaform}) implies that, 
for each $\zeta\in P$, 
\begin{equation}
\tau _{\zeta}:=\tau _{\zeta}(x)
=\op{e}^{b(\zeta ,\,\zeta )+\zeta _l\, X^l},
\label{Xzeta}
\end{equation}
where in the second term in the exponent we use 
Einstein's summation convention. 

Let $T_{\scriptop{f}}$ be the complementary torus 
to the Hamiltonian torus $T_{\scriptop{h}}$ in $T$ which has been used in
Proposition \ref{liftprop}. Finally, let $Z_j$, 
$1\leq j\leq d_{\scriptop{f}}:=\op{dim}T_{\scriptop{f}}$, 
be a $\Z$\--basis of the integral lattice 
$(T_{\scriptop{f}})_{\Z}$ in 
the Lie algebra $\got{t}_{\scriptop{f}}$ of $T_{\scriptop{f}}$. 
For any $X\in\got{t}$ we denote by $X_{\scriptop{h}}$ and 
$X_{\scriptop{f}}$ the $\got{t}_{\scriptop{h}}$\--component 
and the $\got{t}_{\scriptop{f}}$\--component of $X$, respectively. 
With these notations, we have the following conclusions. 
\begin{proposition}
The lattice $\Gamma :=(T_{\scriptop{f}})_{\Z}\times P$ 
is a group with respect to 
the multiplication defined by 
\begin{equation}
(B',\,\beta ')\, (B,\,\beta )=
(B+B'-b_{\scriptop{f}}(\beta ,\beta '),\,\beta +\beta ')
\quad (B,\,\beta ),\, (B',\,\beta ')\in\Gamma .
\label{Gammadef}
\end{equation}
Let $\widetilde{M}:=(\got{t}_{\scriptop{f}}\times N)
\times M_{\scriptop{h}}$. 
Let $(B,\,\beta )\in\Gamma$ act on $\widetilde{M}$ by sending 
$((Z,\,\zeta ),\, x)$ to $((Z',\,\zeta '),\, x')$, where
\begin{eqnarray}
Z'&=&Z+B-\beta _l\, {X^l}_{\scriptop{f}}
+b_{\scriptop{f}}(\beta ,\,\beta )/2
+c_{\scriptop{f}}(\beta ,\,\zeta )/2,
\quad
\zeta '=\zeta +\beta,
\nonumber\\
x'&=&(\op{e}^{c_{\scriptop{h}}(\beta ,\,\zeta )/2}
\, (\tau _{-\beta})_{\scriptop{h}})\cdot x.
\label{Gammaact}
\end{eqnarray}
This defines a proper and free action of $\Gamma$ on 
$\widetilde{M}$, and the mapping 
\begin{equation}
\widetilde{A} :((Z,\,\zeta ),\, x)\mapsto\op{e}^Z\cdot\op{e}^{L_{\zeta}}(x):
\widetilde{M}\to M.
\label{tildeA}
\end{equation}
is a universal covering of $M$ with the action of 
$\Gamma$ on $\widetilde{M}$ as the covering group. 

Let $x\in M_{\scriptop{h}}$ and let 
$\pi _1(M,\, x)$ be the fundamental group 
of $M$ with base point $x$. For any homotopy class 
$[\gamma ]$ of a closed loop $\gamma$ 
based at $x$, let $\iota _x([\gamma ])$ be the element 
of $\Gamma$ of which the action on $\widetilde{M}$ 
is equal to the covering transformation defined by $\gamma$. 
Let $\gamma _j$ be the closed loop 
$\op{e}^{t\, Z_j}\cdot x$, $0\leq t\leq 1$, 
and let $\delta ^l$ be the closed loop based at $x$ which 
consists of $\op{e}^{t\, L_{\varepsilon ^l}}(x)$, $0\leq t\leq 1$, 
followed by $\op{e}^{(1-t)\, X^l}\cdot x$, $0\leq t\leq 1$. 
Then the isomorphism $\iota _x:\pi _1(M,\, x)
\stackrel{\sim}{\to}\Gamma$ is uniquely determined by the 
condition that 
$\iota _x([\gamma _j])=(Z_j,\, 0)$, $1\leq j\leq d_{\scriptop{f}}$, 
and $\iota _x([\delta ^l])=(0,\,\varepsilon ^l)$, 
$1\leq l\leq d_N$.  
\label{pi1prop}
\end{proposition}
\begin{proof}
Let $y\in M$. The surjectivity of the mapping $A$ in (\ref{A}), 
see Proposition \ref{Aprop}, implies that there exist 
$t\in T$, $\zeta\in N$, $x\in M_{\scriptop{h}}$, such that 
$y=t\cdot\op{e}^{L_{\zeta}}(x)$. We have $t=t_{\scriptop{h}}\, 
t_{\scriptop{f}}$ with $t_{\scriptop{h}}\in T_{\scriptop{h}}$ 
and $t_{\scriptop{f}}\in T_{\scriptop{f}}$, and subsequently 
there exists $Z\in\got{t}_{\scriptop{f}}$ such that 
$t_{\scriptop{f}}=\op{e}^Z$. Because $(t_{\scriptop{h}})_M$ 
commutes with $(t_{\scriptop{f}})_M\circ\op{e}^{L_{\zeta}}$ and 
leaves $M_{\scriptop{h}}$ invariant, it follows that 
$y=\widetilde{A}((Z,\,\zeta ),\,t_{\scriptop{h}}\cdot x)$, 
which proves that the mapping $\widetilde{A}$ is surjective. 

Let $\widetilde{A}((Z,\,\zeta ),\, x)
=\widetilde{A}((Z',\,\zeta '),\, x')$. 
The injectivity of the mapping $\alpha$ in Proposition \ref{Aprop} 
implies that there exist $(s,\, -\beta )\in H$ such that 
$(\op{e}^{Z'},\,\zeta ')=(\op{e}^Z,\,\zeta )\, (s,\, -\beta )^{-1}$ 
and $x'=(s\,\tau _{-\beta})\cdot x$. In view of (\ref{Hproduct}) 
and (\ref{Hdef}), this implies that $\beta\in P$, 
$s\,\tau _{-\beta}\in T_{\scriptop{h}}$, 
$\op{e}^{Z'}=\op{e}^Z\, s^{-1}\, 
\op{e}^{c(\beta ,\,\zeta )/2}$,
$\zeta '=\zeta +\beta$, and 
$x'=s\cdot\tau _{-\beta}\cdot x$. 
In view of 
${s_{\scriptop{f}}}^{-1}=(\tau _{-\beta})_{\scriptop{f}}
$ and (\ref{Xzeta}) with $\zeta =-\beta$, 
the $T_{\scriptop{f}}$\--part and the $T_{\scriptop{h}}$\--part 
of the equation for $Z'$ mean that 
\[
Z'\in Z-\beta _l\, {X^l}_{\scriptop{f}}+b_{\scriptop{f}}(\beta ,\,\beta )/2
+c_{\scriptop{f}}(\beta ,\,\zeta )/2+
(T_{\scriptop{f}})_{\Z}
\]
and $s_{\scriptop{h}}=\op{e}^{c_{\scriptop{h}}(\beta ,\,\zeta )/2}$, 
respectively. 

It follows that the fibers of $\widetilde{A}$ are the 
$\Gamma$\--orbits, if we let $(B,\,\beta )\in\Gamma $ 
act on $\widetilde{M}$ as in (\ref{Gammaact}). 
Note that $\zeta '=\zeta$ implies that $\beta =0$, and then 
$Z'=Z$ implies that $B=0$. Therefore the action of 
$\Gamma$ on $\widetilde{M}$ is free, which implies that 
it is effective, in the sense that the mapping from 
$\Gamma$ to the set of diffeomorphisms of 
$\widetilde{M}$ is injective. 

There is a group structure on $\Gamma$ for which 
the action of $\Gamma$ is a group action,   
a homomorphism from $\Gamma$ to the group of 
diffeomorphisms of $\widetilde{M}$, if and only if 
the composition of the actions of two elements of 
$\Gamma$ is an action of an element of $\Gamma$. 
The effectiveness of the action implies that  
if this is the case, then the group structure 
on $\Gamma$ for which this holds 
is unique. 

If we let $(B',\,\beta ')$ 
act on (\ref{Gammaact}), then we arrive at 
$((Z'',\,\zeta ''),\, x'')$, in which 
\begin{eqnarray*}
Z''&=&Z+B-\beta _l\, {X^l}_{\scriptop{f}}
+b_{\scriptop{f}}(\beta ,\,\beta )/2
+c_{\scriptop{f}}(\beta ,\,\zeta )/2\\
&&+B'-\beta '_l\, {X^l}_{\scriptop{f}}
+b_{\scriptop{f}}(\beta ',\,\beta ')/2
+c_{\scriptop{f}}(\beta ',\,\zeta +\beta )/2\\
&=&Z+B+B'-b_{\scriptop{f}}(\beta ,\,\beta ')\\ 
&&-(\beta +\beta ')_l\, {X^l}_{\scriptop{f}}
+b_{\scriptop{f}}(\beta +\beta ',\,\beta +\beta ')/2
+c_{\scriptop{f}}(\beta +\beta ',\,\zeta )/2,
\end{eqnarray*}
$\zeta ''=\zeta +(\beta +\beta ')$, and 
\[
x''=\op{e}^{c_{\scriptop{h}}(\beta ',\,\zeta +\beta )/2}
\cdot (\tau _{-\beta '})_{\scriptop{h}}\cdot 
\op{e}^{c_{\scriptop{h}}(\beta ,\,\zeta )/2}
\cdot (\tau _{-\beta})_{\scriptop{h}}\cdot x
=\op{e}^{c_{\scriptop{h}}(\beta +\beta ',\,\zeta )/2}
\cdot (\tau _{-(\beta +\beta ')})_{\scriptop{h}}\cdot x.
\]
Here we have used that $c(\beta ',\,\beta )
=b(\beta ',\,\beta )-b(\beta ,\,\beta ')$ in the equation for 
$Z''$. Furthermore, in the equation for $x''$ we have used (\ref{tauplus}) 
and the fact that 
if $\beta ,\,\beta '\in P$, then 
$c(\beta ',\,\beta )\in T_{\Z}$, hence 
$c_{\scriptop{h}}(\beta ',\,\beta )\in (T_{\scriptop{h}})_{\Z}$, 
and therefore $\op{e}^{c_{\scriptop{h}}(\beta ',\,\beta )}=1$. 
This proves that $\Gamma$ 
is a group with respect to the multiplication 
defined by (\ref{Gammadef}), and that (\ref{Gammaact}) 
defines a group action of $\Gamma$ on 
$\widetilde{M}$.  

Because the action of $\Gamma$ on $\widetilde{M}$ is obviously proper and 
free, we conclude that $\widetilde{A}$ 
is a covering with covering group equal to the action of 
$\Gamma$. Because $\widetilde{M}$ is simply connected 
as the Cartesian product of a vector space and the 
simply connected Delzant manifold $M_{\scriptop{h}}$, 
see Lemma \ref{pi1delzantlem}, $\widetilde{M}$ is 
a universal covering of $M$. 

It follows from general facts about 
universal coverings, see for instance 
Greenberg \cite[Sec. 5]{greenberg}, that $\iota _x$ is an isomorphism 
from $\pi _1(M,\, x)$ onto $\Gamma$. 
Finally 
$\widetilde{A}$ maps the curve 
$((t\, Z_j,\, 0),\, x)$, $0\leq t\leq 1$, 
which runs from $((0,\, 0),\, x)$ to $((Z_j,\, 0),\, x)$, 
to $\gamma _j$. Furthermore $\widetilde{A}$ maps the curve 
$((0,\, t\,\varepsilon ^l),\, x)$, 
$0\leq t\leq 1$, followed by the curve 
$((-t\, {X^l}_{\scriptop{f}},\, 0),\, 
\op{e}^{-t\, {X^l}_{\scriptop{h}}}\cdot x)$, $0\leq t\leq 1$, 
which runs from $((0,\, 0), x)$ to 
\[
((-{X^l}_{\scriptop{f}},\,\varepsilon ^l),\, 
\op{e}^{-{X^l}_{\scriptop{h}}}\cdot x)
=(0,\,\varepsilon ^l)\cdot ((0,\, 0),\, x),
\]
to $\delta ^l$. This shows that 
$\iota _x([\gamma _j])=(Z_j,\, 0)$ 
and $\iota _x([\delta ^l])=(0,\,\varepsilon ^l)$. 
Because the elements $(Z_j,\, 0)$, $1\leq j\leq d_{\scriptop{f}}$, 
and $(0,\,\varepsilon ^l)$, $1\leq l\leq d_N$, 
together generate $\Gamma$, this proves the last statement 
in the proposition. 
\end{proof}
Viewing $\got{t}_{\scriptop{f}}$ as an additive group, 
the connected commutative Lie group 
$U:=\got{t}_{\scriptop{f}}\times T_{\scriptop{h}}$ 
acts on $\widetilde{M}$, where 
$(Z',\, t)\in U$ 
sends $((Z,\,\zeta ),\, x)$ to 
$((Z+Z',\,\zeta ),\, t\cdot x)$. 
The covering map $\widetilde{A}:\widetilde{M}\to M$ 
intertwines the $U$\--action 
on $\widetilde{M}$ with the $T$\--action on $M$ 
via the covering homomorphism 
$\epsilon :(Z',\, t)\mapsto\op{e}^{Z'}\, t: U\to T$, 
in the sense that $\widetilde{A}(u\cdot p)=
\epsilon (u)\cdot\widetilde{A}(p)$ for every 
$p\in\widetilde{M}$ and $u\in U$. 
\begin{corollary}
The fundamental group of $M$ is commutative if and only if 
$c(P\times P)\subset (T_{\scriptop{h}})_{\Z}$. 
The first homology group $\op{H}_1(M,\,\Z)$ of $M$ with coefficients 
in $\Z$ is isomorphic to 
$((T_{\scriptop{f}})_{\Z}/\Theta )
\times P$, 
in which $\Theta$ 
denotes the additive subgroup of $(T_{\scriptop{f}})_{\Z}$ 
which is generated by the elements 
$c_{\scriptop{f}}(\beta ,\,\beta ')$, such that 
$\beta ,\,\beta '\in P$. 
The first Betti number $\op{dim}\op{H}_1(M,\,\R)$ is equal to 
$\op{dim}M-2\op{dim}T_{\scriptop{h}}-\op{rank}\Theta$. 
\label{pi1cor}
\end{corollary}
\begin{proof}
A straightforward computation shows that 
$(B,\,\beta )^{-1}=(-B-b(\beta ,\,\beta )_{\scriptop{f}}
,\, -\beta )$, that 
\[
(B,\,\beta )^{-1}\, (B',\,\beta ')\, (B,\,\beta ) 
=(B'+c_{\scriptop{f}}(\beta ',\,\beta ),\,\beta '),
\] 
and that the commutator 
$(B',\,\beta ')^{-1}\, (B,\,\beta )^{-1}
\, (B',\,\beta ')\, (B,\,\beta )$ 
is equal to $(c_{\scriptop{f}}(\beta ',\,\beta ),\, 0)$. 
Therefore the subgroup of $\Gamma =(T_{\scriptop{f}})_{\Z}\times P$ 
generated by the commutators 
is equal to $\Theta\times\{ 0\}$, 
and $\Gamma$ 
is commutative if and only if $\Theta =\{ 0\}$ 
if and only if $c(P\times P)\subset 
T_{\Z}\cap\op{t}_{\scriptop{h}}=(T_{\scriptop{h}})_{\Z}$. 

The canonical homomorphism $\pi _1(M,\, x)\to\op{H}_1(M,\,\Z)$ 
is surjective with kernel equal to the subgroup 
of $\pi _1(M,\, x)$ generated by the commutators, 
see Greenberg \cite[Th. 12.1]{greenberg}.  
This induces an isomorphism from 
$((T_{\scriptop{f}})_{\Z}\times P)/(\Theta\times\{ 0\})
=((T_{\scriptop{f}})_{\Z}/\Theta )\times P$
onto $\op{H}_1(M,\,\Z)$. Finally the universal coefficient theorem, 
cf. Greenberg \cite[Th. 29.12]{greenberg} implies that 
for any principal ideal domain $R$, in particular for 
$R=\R$, ~$\op{H}_1(M,\, R)$ is isomorphic to 
$\op{H}_1(M,\,\Z )\otimes _{\Z}\R$. Therefore  
\[
\op{dim}\op{H}_1(M,\,\R )=d_{\scriptop{f}}-\op{rank}\Theta +d_N
=\op{dim}T-\op{dim}T_{\scriptop{h}}-\op{rank}\Theta +\op{dim}\got{l}
-\op{dim}T_{\scriptop{h}}, 
\]
which is equal to ~$\op{dim}M-2\op{dim}T_{\scriptop{h}}
-\op{rank}\Theta$ in view of  
Lemma \ref{coisotropicorbitlem}.
\end{proof}
It follows that the rank of $\Theta$ is a purely topological 
feature of $M$, disregarding both the $T$\--action and the 
symplectic structure on $M$. Note that the generators 
of $\Gamma$ mentioned in Proposition \ref{pi1prop} 
were defined in terms of the action of $T_{\scriptop{f}}$ 
and the $L_{\varepsilon ^l}$, where the latter were defined in 
terms of both the $T$\--action and the symplectic form 
on $M$. 
\begin{remark}
The symplectic form 
$\widetilde{A}^*\sigma$ on the universal covering 
$\widetilde{M}=\got{t}_{\scriptop{f}}\times N\times 
M_{\scriptop{h}}$ is given by (\ref{A*sigma}), 
in which $a$, $\delta a$, $\delta 'a$ 
are replaced by $((Z,\,\zeta ),\, x)$, 
$((\delta Z,\,\delta\zeta ),\,\delta x)$, 
$((\delta 'Z,\,\delta '\zeta ),\,\delta 'x)$, 
respectively, with $Z,\,\delta Z,\,\delta 'Z\in\got{t}_{\scriptop{f}}$. 
We view the linear form 
\[
\mu (x)\, c_{\scriptop{h}}(\cdot ,\,\zeta ):
\delta\zeta\mapsto \mu (x)(c_{\scriptop{h}}(\delta\zeta ,\,\zeta ))
\]
on $N:=(\got{l}/\got{t}_{\scriptop{h}})^*$ as an element of 
$((\got{l}/\got{t}_{\scriptop{h}})^*)^*\simeq
\got{l}/\got{t}_{\scriptop{h}}\simeq\got{l}\cap\got{t}_{\scriptop{f}}$.   
Let $\Psi :\widetilde{N}\to\widetilde{N}$ be defined by 
\[
\Psi:((Z,\,\zeta ),\, x)=
((Z+\mu (x)\, c_{\scriptop{h}}(\cdot ,\,\zeta )/2
,\,\zeta ),\, x),
\quad ((Z,\,\zeta ),\, x)\in\widetilde{M}=
(\got{t}_{\scriptop{f}}\times N)\times M_{\scriptop{h}}.
\]
Then $\Psi$ is a diffeomorphism from $\widetilde{M}$ 
onto $\widetilde{M}$, and the symplectic form 
$\nu :=\Psi ^*(\widetilde{A}^*\sigma )$ on $\widetilde{N}$ 
is given by 
\[
\nu _a(\delta a,\,\delta 'a)
=\sigma ^{\got{t}}(\delta Z,\,\delta 'Z)
+\delta\zeta (\delta 'Z_{\got{l}})
-\delta '\zeta (\delta Z_{\got{l}})
+(\sigma _{\scriptop{h}})_x(\delta x,\,\delta 'x).
\]
That is, $(\widetilde{M},\,\nu )$ is equal to the 
Cartesian product of a symplectic vector space 
$(\got{t}_{\scriptop{f}}\times N,\,\sigma 
^{\got{t}_{\scriptop{f}}\times N})$ and the Delzant manifold 
$(M_{\scriptop{h}},\,\sigma _{\scriptop{h}})$. 
Here 
\[
\sigma ^{\got{t}_{\scriptop{f}}\times N}
((\delta Z,\,\delta\zeta ),\, (\delta 'Z,\,\delta '\zeta ))
=\sigma ^{\got{t}}(\delta Z,\,\delta 'Z)
+\delta\zeta (\delta 'Z_{\got{l}})
-\delta '\zeta (\delta Z_{\got{l}}).
\]
Because $\Psi$ is 
$(\got{t}_{\scriptop{f}}\times T_{\scriptop{h}})$\--equivariant,  
we have recovered Benoist \cite[Cor. 6.16]{benoist}, in which the 
``cocycle $c$'' is equal to our $\sigma ^{\got{t}}$. 
\label{benoist6.16rem}
\end{remark}

\section{The classification}
\label{invariantsec}
\subsection{Invariants}
\label{uniquenesssubsec}
The model in Proposition \ref{Aprop}, of a 
compact connected symplectic manifold $(M,\,\sigma )$ 
with an effective symplectic action of a torus $T$ of which the 
principal orbits are coisotropic submanifolds of $(M,\,\sigma )$,   
has been described in terms of the following ingredients. 
\begin{definition}
Let $T$ be a given torus. 
A {\em list of ingredients for $T$} consists of: 
\begin{itemize}
\item[1)] An antisymmetric bilinear form $\sigma ^{\got{t}}$ 
on the Lie algebra $\got{t}$ of $T$. 
\item[2)] A subtorus $T_{\scriptop{h}}$ of $T$, of which the Lie 
algebra $\got{t}_{\scriptop{h}}$ is contained in 
$\got{l}:=\op{ker}\sigma ^{\got{t}}$.  
\item[3)] A Delzant polytope $\Delta$ in 
${{\got t}_{\scriptop{h}}}^*$ with center of mass at the origin. 

\item[4)] A discrete cocompact subgroup $P$ of the 
additive subgroup $N:=(\got{l}/\got{t}_{\scriptop{h}})^*$ of 
$\got{l}^*$. 
\item[5)] An antisymmetric bilinear mapping 
$c:N\times N\to\got{l}$ 
with the following properties. 
\begin{itemize}
\item[5a)]  
For every $\zeta ,\,\zeta '\in P$, the element 
$c(\zeta ,\,\zeta ')\in\got{l}\subset\got{t}$ belongs to the integral lattice 
$T_{\Z}$ in $\got{t}$, the kernel of the exponential 
mapping $\op{exp}:\got{t}\to T$. 
\item[5b)] For every $\zeta ,\,\zeta ',\,\zeta ''\in N$ 
we have that 
\[
\zeta (c(\zeta ',\,\zeta ''))
+\zeta '(c(\zeta '',\,\zeta ))
+\zeta ''(c(\zeta ,\,\zeta '))=0.
\]
\end{itemize}
\item[6)] An element $\overline{\tau}$ of the space 
${\cal T}$ which has been defined in (\ref{calT}). 
\end{itemize}
\label{ingredientdef}
\end{definition}
\begin{remark}
Regarding the Delzant polytope $\Delta$ in 3) in 
Definition \ref{ingredientdef}, 
we have a corresponding Delzant manifold 
$(M_{\scriptop{h}},\,\sigma _{\scriptop{h}})$, which is a 
$2\op{dim}T_{\scriptop{h}}$\--dimensional 
compact connected symplectic manifold, equipped with an effective 
Hamiltonian $T_{\scriptop{h}}$\--action on $(M_{\scriptop{h}},\,\sigma _{\scriptop{h}})$, 
for which $\Delta$ is equal to the image of the momentum map. 

In 5b), $\zeta\in N$ is viewed as a linear form on 
$\got{l}$ which vanishes on $\got{t}_{\scriptop{h}}$, 
so $\zeta (c(\zeta ',\,\zeta ''))$ is a real number. 

The ingredient 6), the holonomy invariant, has been introduced in 
Subsection \ref{holinvariantss}. As explained there, the 
space ${\cal T}$ to which it belongs can have a non\--Hausdorff 
quotient topology. 
\end{remark}

\begin{definition}
Let $M$ be a compact and connected smooth manifold, $\sigma$ 
a symplectic form on $M$, and $T$ a torus acting 
effectively on $(M,\,\sigma )$ 
by means of symplectomorphisms and with coisotropic principal orbits. 
The {\em list of ingredients of $(M,\,\sigma ,\, T)$}, 
as in Definition \ref{ingredientdef}, consists of:   
\begin{itemize}
\item[i)] $\sigma ^{\got{t}}(M,\,\sigma ,\, T)$ 
is the antisymmetric bilinear form $\sigma ^{\got{t}}$ on 
$\got{t}$ as defined in Lemma \ref{constlem}. 
\item[ii)] $T_{\scriptop{h}}(M,\,\sigma ,\, T)$ 
is the Hamiltonian torus $T_{\scriptop{h}}$, the 
unique maximal stabilizer subgroup $T_{\scriptop{h}}$ for the 
$T$\--action on $M$, see Remark \ref{T1rem} and 
Lemma \ref{stabilizerlem}. 
\item[iii)] $\Delta (M,\,\sigma ,\, T)$ is the image 
$\Delta =\mu (M)$ of the 
momentum map $\mu :M\to{\got{t}_{\scriptop{h}}}^*$ of the 
$T_{\scriptop{h}}$\--action on $(M,\,\sigma )$, which is 
Hamiltonian, cf. Corollary \ref{Tstabcor}, where we eliminated the 
translational ambiguity by putting the center of mass of 
$\Delta$ at the origin. $\Delta$ is a translate of the Delzant 
polytope $\Delta _p$ in Proposition \ref{M/Tprop}. 
\item[iv)] $P(M,\,\sigma ,\, T)$ is the period group $P$ defined in 
Lemma \ref{perlem} with $Q=M/T$, $V=\got{l}^*$, and 
$N=(\got{l}/\got{t}_{\scriptop{h}})^*$, 
which according to Proposition \ref{M/Tprop} is a discrete 
cocompact additive subgroup of $N$. 
\item[v)] $c(M,\,\sigma ,\, T)$ is the 
antisymmetric bilinear mapping $c:N\times N\to\got{l}$ 
defined in Proposition \ref{liftprop}. 
\item[vi)] $\overline{\tau}(M,\,\sigma ,\, T)$
is the holonomy invariant of $(M,\,\sigma ,\, T)$, the element
$\overline{\tau}$ of ${\cal T}$ defined in (\ref{overlinetau}).  
\end{itemize}
Note that all the ingredients in Definition \ref{ingredientdef} 
are defined only in terms of the torus $T$. 
\label{invariantdef}
\end{definition}

\begin{theorem}
Let $T$ be a torus. 
The list of ingredients of $(M,\,\sigma ,\, T)$ 
is a complete set of invariants for the compact 
connected symplectic manifold $(M,\,\sigma )$ with 
effective symplectic $T$\--action with coisotropic principal orbits, 
in the following sense. 
If $(M',\,\sigma ')$ is another 
compact connected symplectic manifold with 
effective symplectic $T$\--action 
with coisotropic principal orbits, then there exists a $T$\--equivariant 
symplectomorphism $\Phi$ from $(M,\,\sigma ,\, T)$ onto 
$(M',\,\sigma ',\, T)$ if and only if the list of ingredients 
of $(M,\,\sigma ,\, T)$ is equal to the 
list of ingredients of $(M',\,\sigma ',\, T)$. 
\label{invariantthm}
\end{theorem}

\begin{proof}
The property 5a) in Definition \ref{ingredientdef} 
follows from Remark \ref{integralcrem}, and also from 
Lemma \ref{TPlem}. Equation   
5b) in Definition \ref{ingredientdef} 
is the equation (\ref{ft}). 

Suppose that $\Phi$ is a $T$\--equivariant 
symplectomorphism from $(M,\,\sigma ,\, T)$ onto 
$(M',\,\sigma ',\, T)$. We will check that the 
ingredients of $(M,\,\sigma ,\, T)$ and 
$(M',\,\sigma ',\, T')$ are the same. In other words, 
the ingredients are invariants of the symplectic $T$\--spaces. 

If $X,\, Y\in\got{t}$, then the $T$\--equivariance 
of $\Phi$ implies that $\Phi ^*\, X_{M'}=X_M$ and  
$\Phi ^*\, Y_{M'}=Y_M$. In combination with 
$\sigma =\Phi ^*\sigma '$, this implies in view of 
Lemma \ref{constlem} that
\begin{eqnarray*}
\sigma ^{\got{t}}(M,\,\sigma ,\, T)(X,\, Y)
&=&\sigma (X_M,\, Y_M)=(\Phi ^*\sigma ')(
\Phi ^*X_{M'},\,\Phi ^*Y_{M'})
=\Phi ^*(\sigma '(X_{M'},\, Y_{M'}))\\
&=&\Phi ^*(\sigma ^{\got{t}}(M',
\,\sigma ',\, T)(X,\, Y))
=\sigma ^{\got{t}}(M',
\,\sigma ',\, T)(X,\, Y),
\end{eqnarray*}
where we have used in the last equation that 
$\sigma ^{\got{t}}(M',\,\sigma ',\, T)(X,\, Y)$ 
is a constant on $M'$. This proves that 
$\sigma ^{\got{t}}(M,\,\sigma ,\, T)
=\sigma ^{\got{t}}(M',\,\sigma ',\, T)$. 
The $T$\--equivariance of $\Phi$  
implies that $T_{\Phi (x)}=T_x$ for every $x\in M$, 
and therefore $T_{\scriptop{h}}(M',\,\sigma ',\, T)
=T_{\scriptop{h}}(M,\,\sigma ,\, T)$. 

In combination with $\Phi ^*\sigma '=\sigma$, the $T$\--equivariance 
of $\Phi$ implies that the $\got{l}^*$\--valued 
closed basic one\--form $\widehat{\sigma}$ defined 
in Lemma \ref{alphalem} is equal to 
$\Phi ^*\widehat{\sigma '}$. It follows that $\Phi$ 
induces an isomorphism of locally convex polyhedral 
$\got{l}^*$\--parallel spaces from $M/T$ onto $M'/T$. 
In view of Proposition \ref{M/Tprop} this implies that 
$P(M',\,\sigma ',\, T)=P(M,\,\sigma ,\, T)$ and 
that $\Delta (M',\,\sigma ',\, T)$ is a translate of 
$\Delta (M,\,\sigma ,\, T)$ in ${\got{t}_{\scriptop{h}}}^*$. 
Because both 
$\Delta (M',\,\sigma ',\, T)$ and
$\Delta (M,\,\sigma ,\, T)$ have their 
center of mass at the origin, it follows that 
$\Delta (M',\,\sigma ',\, T)=\Delta (M,\,\sigma ,\, T)$. 

The $T$\--equivariant symplectomorphism $\Phi$ 
maps an admissible connection as in Proposition \ref{liftprop} 
to an admissible connection as in Proposition \ref{liftprop} 
with $(M,\,\sigma ,\, T)$ replaced by $(M',\,\sigma ',\, T)$. 
it follows that $c(M',\,\sigma ',\ T)=c(M,\,\sigma ,\, T)$ in view of 
Remark \ref{curvrem} and  
$\overline{\tau}(M',\,\sigma ',\ T)
=\overline{\tau}(M,\,\sigma ,\, T)$ 
in view of Subsection \ref{holinvariantss}. 
This proves the ``only if'' part 
of the theorem.  

For the ``if'' part, the completeness of the invariants, 
we observe that the manifold 
$M_{\scriptop{model}}:=G\times _HM_{\scriptop{h}}$ and the 
$T$\--invariant symplectic 
form $\sigma _{\scriptop{model}}$ on $M_{\scriptop{model}}$, 
see Proposition \ref{Aprop} and Proposition \ref{omegalem}, are defined 
in terms of the ingredients 1) -- 6) in Definition \ref{ingredientdef}, 
and the elements $\tau_{\zeta}$, $\zeta\in P$. 
Let $x\in M$ and choose an admissible connection for $(M,\sigma ,\, T)$ 
as in Proposition \ref{liftprop}. Then 
$\overline{\tau}(M',\,\sigma ',\ T)
=\overline{\tau}(M,\,\sigma ,\, T)$ implies in view of  
Subsection \ref{holinvariantss} that 
there exist $x'\in M'$ and a choice 
of an admissible connection for $(M',\,\sigma ',\, T)$ 
as in Proposition \ref{liftprop}, such that 
the holonomy $\tau '_{\zeta}(x')$, $\zeta\in P$, 
defined by this connection and with the initial point $x'$,  
is equal to $\tau _{\zeta}=\tau _{\zeta}(x)$, $\zeta\in P$. 
Therefore the model for $(M',\,\sigma ',\, T)$ in 
Proposition \ref{Aprop}, with $(M,\,\sigma ,\, T)$ 
replaced by $(M',\,\sigma ',\, T)$, can be chosen to be 
equal to the model for 
$(M,\,\sigma ,\, T)$ in Proposition \ref{Aprop}.  
This implies the existence 
of a $T$\--equivariant 
symplectomorphism $\alpha '$ from 
$(M_{\scriptop{model}},\,\sigma _{\scriptop{model}},\, T)$ onto 
$(M',\,\sigma ',\, T)$, and it follows that  
$\Phi :=\alpha '\circ\alpha ^{-1}$ is a $T$\--equivariant 
symplectomorphism from $(M,\,\sigma ,\, T)$ onto 
$(M',\,\sigma ', T)$.
\end{proof}

\begin{remark}
Because $\Delta (M',\,\sigma ',\, T)=\Delta (M,\,\sigma ,\, T)$, 
the point $x'\in M'$ in the last paragraph of the proof of 
Theorem \ref{invariantthm} can be chosen such that 
$\mu '(x ')=\mu (x)$, where $\mu$ and $\mu '$ denote the 
momentum maps of the 
Hamiltonian $T_{\scriptop{h}}$\--actions on $(M,\,\sigma )$ 
and $(M',\sigma ')$, respectively. This implies that 
there is a $T_{\scriptop{h}}$\--equivariant symplectomorphism 
$\Phi _{\scriptop{h}}$ from the Delzant submanifold of 
$(M,\,\sigma )$ through $x$ 
onto the Delzant submanifold of $(M',\,\sigma ')$ through 
$x'$, which maps $x$ to $x'$. Using $\Phi _{\scriptop{h}}$ 
in order to identify both Delzant manifolds with 
$(M_{\scriptop{h}},\,\sigma _{\scriptop{h}},\, T_{\scriptop{h}})$, 
under which identifications $x$ and $x'$ are mapped to the same 
point of $M_{\scriptop{h}}$, we conclude that $\Phi (x)=x'$, if 
$\Phi$ is  the $T$\--equivariant symplectomorphism 
from $(M,\,\sigma ,\, T)$ to $(M',\,\sigma ',\, T)$,   
described in the last paragraph of the proof of 
Theorem \ref{invariantthm}. 

Let $\op{Aut}(M,\,\sigma ,\, T)$ denote the automorphism 
group of $(M,\,\sigma ,\, T)$, the set of all 
$T$\--equivariant symplectomorphisms from 
$(M,\,\sigma ,\, T)$ to $(M,\,\sigma ,\, T)$. 
Each $\Phi\in\op{Aut}(M,\,\sigma ,\, T)$ induces 
a transformation $\Phi _{M/T}$ of $M/T$, which is 
an isomorphism of $\got{l}^*$\--parallel spaces, 
and therefore of the form $p\mapsto p+\nu (\Phi )$ 
for a unique element $\nu ({\Phi})\in N/P$. 
The mapping $\nu :\Phi\mapsto\nu (\Phi )$ is a 
homomorphism from the group $\op{Aut}(M,\,\sigma ,\, T)$ 
to the torus $N/P$. Using the previous paragraph  
with $M'=M$ and $\sigma '=\sigma$ and using 
Subsection \ref{holinvariantss}, it can be proved 
that $\nu (\op{Aut}(M,\,\sigma ,\, T))$ is equal to 
the set of $\zeta '+P\in N/P$, for which there exists 
an $\alpha\in\op{Sym}$ such that 
$\op{e}^{c(\zeta ,\,\zeta ')}=\op{e}^{\alpha _{\zeta}}$ 
for all $\zeta\in P$, where it is sufficient to satisfy these 
equations for all $\zeta$ in a $\Z$\--basis of $P$. 

Using this one can prove that 
$\nu (\op{Aut}(M,\,\sigma ,\, T))$ is a Lie subgroup 
of $N/P$ with Lie algebra equal to $c^0$, the space 
of all elements of $N$ which are equal to zero 
on the span of $c(N\times N)$. 
Actually, $L_{\zeta}$ is an inifinitesimal 
symplectomorphism if and only if $\zeta\in c^0$. In general  
$\nu (\op{Aut}(M,\,\sigma ,\, T))$ need not be a closed 
subgroup of the torus $N/P$, and it neither needs to be 
connected, but it has countably many connected 
components. 

The kernel of the homomorphism $\nu$ from $\op{Aut}(M,\,\sigma ,\, T)$ 
to $N/P$ consists of the group \\ $\op{Aut}_T(M,\,\sigma ,\, T)$ 
of all $T$\--equivariant symplectomorphisms 
$\Phi :(M,\,\sigma ,\, T)\to (M,\,\sigma ,\, T)$ which 
preserve all the $T$\--orbits. This group can be analyzed 
starting from Remark \ref{hsrem}.  
\label{autrem}
\end{remark}

\subsection{Existence}
The following existence theorem completes the classification.
\begin{theorem}
Every list of ingredients as in Definition 
\ref{ingredientdef} 
is equal to the list of 
invariants of a compact connected symplectic manifold $(M,\,\sigma )$ 
with effective symplectic 
$T$\--action with coisotropic principal orbits as in 
Theorem \ref{invariantthm}. 
\label{existencethm}
\end{theorem}
\begin{proof}
A straightforward verification shows that, for any antisymmetric 
bilinear mapping $c:N\times N\to\got{l}$, 
(\ref{gotg}) turns $\got{g}:=\got{t}\times N$ 
into a two\--step nilpotent Lie algebra, and that 
(\ref{Hproduct}) defines a product in $G:=T\times N$ 
for which $G$ is a Lie group with $\got{g}$ as its Lie algebra. 

Choose an element $\tau\in\op{Hom}_c(P,\, T)$ 
such that $\overline{\tau}=(\op{exp}{\cal A})\cdot\tau$, see 
(\ref{calT}). Because the $\tau _{\zeta}$, $\zeta\in P$, 
satisfy (\ref{tauplushom}), 
it follows that (\ref{Hdef}) defines a closed Lie subgroup $H$ of $G$, 
and that (\ref{Hact}) defines a smooth action of $H$ on the Delzant 
manifold $M_{\scriptop{h}}$. Here we have used a choice 
of a complementary torus $T_{\scriptop{f}}$ to $T_{\scriptop{h}}$, 
which will be kept fixed in the remainder of the proof.  

Because $H$ is a closed Lie subgroup of $G$, its right action 
on $G$ is proper and free, and therefore the action of $H$ 
on $G\times M_{\scriptop{h}}$, for which $h\in H$ sends $(g,\, x)$ to 
$(g\, h^{-1},\, h\cdot x)$, is proper and free. The 
orbit space $M:=G\times _HM_{\scriptop{h}}$ 
has a unique structure of a smooth 
manifold such that the canonical projection 
$\pi _M:G\times M_{\scriptop{h}}\to M$ 
is principal $H$\--bundle. Because 
$G$ and $M_{\scriptop{h}}$ are connected, $M$ is connected as the 
image of the connected set $G\times M_{\scriptop{h}}$ under the 
continuous mapping $\pi _M$. The projection 
$(g,\, x)\mapsto g$ induces a $G$\--equivariant smooth 
fibration $\psi :M\to G/H$ with fiber $M_{\scriptop{h}}$, 
the fiber bundle induced from the principal fiber bundle $G\to G/H$ 
by means of the action of $H$ on $M_{\scriptop{h}}$. 
See \cite[Sec. 2.4]{dk}. Because $P$ is cocompact in $N$, 
the base space $G/H$ is compact, and because the fiber 
$M_{\scriptop{h}}$ is a 
compact Delzant manifold, it follows that $M$ is compact. 

On $G\times M_{\scriptop{h}}$ we define the smooth two\--form $\omega$ 
by (\ref{A*sigma}). 
(Note that we cannot use 
the equation $\omega =A^*\,\sigma$ here, 
because we do not have the symplectic form 
$\sigma$ on the manifold $M$ yet, 
we are in the process of defining it.)  
In (\ref{A*sigma}) we have used a choice of a linear projection 
$X\mapsto X_{\got{l}}$ 
from $\got{t}$ onto $\got{l}$, which will be kept 
fixed in the remainder of the proof. 
 
We first verify that $\omega$ is closed. 
The part 
\[
\sigma ^{\got{t}}(\delta t,\,\delta 't)
+\delta\zeta (\delta 't_{\scriptop{l}})
-\delta '\zeta (\delta t_{\scriptop{l}})
\]
of (\ref{A*sigma}) is closed, because it is defined by a constant 
two\--form on $T\times N$.  

For the part 
\[
\varphi _{\zeta}(\delta\zeta ,\,\delta '\zeta )
:=\delta\zeta (c(\delta '\zeta ,\,\zeta )/2)
-\delta '\zeta (c(\delta\zeta ,\,\zeta )/2)
\]
of (\ref{A*sigma}), it follows from 
(\ref{domega}) that 
$(\op{d}\!\varphi )(\delta\zeta ,\,\delta '\zeta ,\,
\delta ''\zeta )$ is equal to the cyclic sum over 
$\delta\zeta ,\,\delta '\zeta ,\,\delta ''\zeta$ of 
$\delta '\zeta (c(\delta ''\zeta ,\,\delta\zeta )/2)
-\delta ''\zeta (c(\delta '\zeta ,\,\delta\zeta )/2)$, 
which is equal to zero because of 5b) in Definition 
\ref{ingredientdef}. 
 
If $A_{\scriptop{h}}:(t,\, x)\mapsto t\cdot x:T_{\scriptop{h}}
\times M_{\scriptop{h}}\to M_{\scriptop{h}}$ 
denotes the action of $T_{\scriptop{h}}$ on 
$M_{\scriptop{h}}$, then the part 
\[
(\sigma _{\scriptop{h}})_x(\delta x,\, (\delta 't_{\scriptop{h}})_M(x)) 
-(\sigma _{\scriptop{h}})_x(\delta 'x,\, (\delta 't_{\scriptop{h}})_M(x)) 
+(\sigma _{\scriptop{h}})_x(\delta x,\,\delta 'x)
\]
of (\ref{A*sigma}) is equal to the pull\--back of 
${A_{\scriptop{h}}}^*\sigma _{\scriptop{h}}$ by means of the 
mapping 
\[
p:((t,\,\zeta ),\, x)\mapsto (t_{\scriptop{h}},\, x)
:(T\times N)\times M_{\scriptop{h}}\to T_{\scriptop{h}}\times 
M_{\scriptop{h}}.
\]
This part of 
(\ref{A*sigma}) is closed, because 
~$\op{d}(p^*({A_{\scriptop{h}}}^*\sigma _{\scriptop{h}}))
=p^*(\op{d}({A_{\scriptop{h}}}^*\sigma _{\scriptop{h}}))
=p^*({A_{\scriptop{h}}}^*(\op{d}\!\sigma _{\scriptop{h}}))=0$. 

The remaining part of (\ref{A*sigma}) is 
\[
-\mu (x)(c_{\scriptop{h}}(\delta\zeta ,\,\delta '\zeta ))
+(\sigma _{\scriptop{h}})_x(\delta x,
\, c_{\scriptop{h}}(\delta '\zeta ,\,\zeta )_{M_{\scriptop{h}}}(x))/2
-(\sigma _{\scriptop{h}})_x(\delta 'x,
\, c_{\scriptop{h}}(\delta \zeta ,\,\zeta )_{M_{\scriptop{h}}}(x))/2.
\]
Because the action of $T_{\scriptop{h}}$ on $M_{\scriptop{h}}$ 
is Hamiltonian with momentum mapping $\mu$, we have for 
every $Y\in\got{t}_{\scriptop{h}}$ that 
$(\sigma _{\scriptop{h}})_x(\delta x,\, Y_{M_{\scriptop{h}}}(x))$ 
is equal to the derivative of $x\mapsto \mu (x)(Y)$ in the 
direction of $\delta x$. If we apply this to 
$Y=c_{\scriptop{h}}(\delta '\zeta ,\,\zeta )/2$, 
then we obtain that the remaining part of (\ref{A*sigma}) 
is equal to $\op{d}\!\gamma$, in which the one\--form 
$\gamma$ is defined by 
\[
\gamma _{((t,\,\zeta ),\, x)}((\delta\gamma ,\,\delta\zeta ),\,\delta x)
=\mu (x)(c_{\scriptop{h}}(\delta\zeta ,\,\zeta ))/2.
\]
Because $\op{d}(\op{d}\!\gamma )=0$, this completes 
the proof that $\op{d}\!\omega =0$.  

The element $(b,\,\beta )\in H$ sends 
$((t,\,\zeta ),\, x)$ to $((\widetilde{t},\,\widetilde{\zeta}),\, 
\widetilde{x})$ with 
$\widetilde{t}=t\, b^{-1}\,\op{e}^{c(\zeta ,\,\beta )/2}$, 
$\widetilde{\zeta}=\zeta -\beta$, and $\widetilde{x}=
(b\,\tau _{\beta})_{\scriptop{h}}\cdot x$. 
Therefore the tangent map of the action of $(b,\,\beta )$ 
sends $((\delta t,\,\delta\zeta ),\, \delta x)$ to 
$((\widetilde{\delta}t,\,\widetilde{\delta}\zeta ),\,\widetilde{\delta}x)$ 
with $\widetilde{\delta}t=\delta t+c(\delta\zeta ,\,\beta )/2$, 
$\widetilde{\delta}\zeta=\delta\zeta$,  and 
$\widetilde{\delta}x
=\op{T}_x((b\,\tau _{\beta})_{\scriptop{h}})_{M_{\scriptop{h}}}\,\delta x$.
Because $\delta t+c(\delta\zeta ,\,\zeta )/2
=\widetilde{\delta}t+c(\widetilde{\delta}\zeta ,\,
\widetilde{\zeta})/2$, and because 
$((b\,\tau _{\beta})_{\scriptop{h}})_{M_{\scriptop{h}}}$ 
is a symplectomorphism on $M_{\scriptop{h}}$ which leaves 
$\mu$ and infinitesimal $T_{\scriptop{h}}$\--actions 
invariant, it follows that the two\--form $\omega$ defined by 
(\ref{A*sigma}) is $H$\--invariant.  
 
The condition that $\omega _a(\delta a,\,\delta 'a)=0$ 
for every $\delta 'a\in\op{T}_a(G\times M_{\scriptop{h}})$ is equivalent 
to $\delta\zeta =0$ (take $\delta '\zeta =0$, 
$\delta 'x=0$, and let $\delta 't$ range over 
$\got{l}\cap\got{t}_{\scriptop{f}}$), 
$\delta t\in\got{l}$ (take $\delta '\zeta =0$, 
$\delta 'x=0$, and let $\delta 't$ range over 
$\got{t}_{\scriptop{f}}$, where we use that we already have 
$\delta\zeta =0$), 
$\delta x+(\delta t_{\scriptop{h}})_{M_{\scriptop{h}}}(x)=0$ 
(take in the remaining equation (\ref{A*sigma}) 
$\delta '\zeta =0$, $\delta 't=0$ and let 
$\delta 'x$ range over $\op{T}_x\! M_{\scriptop{h}}$), 
and finally $\delta t\in\got{t}_{\scriptop{h}}$, 
because the fact that the $T_{\scriptop{h}}$\--orbits in 
$M_{\scriptop{h}}$ are isotropic 
now implies that $-\delta '\zeta (\delta t)=0$ 
for all $\delta '\zeta\in (\got{l}/\got{t}_{\scriptop{h}})^*$. 
It follows that the kernel of $\omega _a$ is equal 
to $\op{T}_a(H\cdot a)=\op{ker}(\op{T}_a\!\pi _M)$. 

The conclusion is that $\omega$ is a basic two\--form 
for the action of $H$ on $G\times M_{\scriptop{h}}$, 
which implies that there is 
a unique smooth two\--form $\sigma$ on $M=G\times _HM_{\scriptop{h}}$ 
such that $\omega ={\pi _M}^*\sigma$. Because 
${\pi _M}^*(\op{d}\!\sigma )=\op{d} ({\pi _M}^*\sigma )=
\op{d}\!\omega =0$ and at every point the tangent mapping 
of $\pi _M$ is surjective, we have that $\op{d}\!\sigma =0$. 
Furthermore $\sigma$ is nondegenerate at every point, because 
the kernel of $\omega$ is equal to the kernel of the tangent mapping 
of $\pi _M$ at every point. Therefore $\sigma$ is 
a symplectic form on $M$. 

On $G\times M_{\scriptop{h}}$ we have the 
action of $s\in T$ which sends $((t,\,\zeta ),\, x)$ 
to $((s\, t,\,\zeta ),\, x)$. This action clearly leaves 
$\omega$ invariant, and it follows that the induced action 
of $T$ on $M:=G\times _HM_{\scriptop{h}}$ leaves $\sigma$ invariant. 
The tangent vectors to the orbits in $G\times M_{\scriptop{h}}$ 
are the $((\delta s ,\, 0),\, 0)$, $\delta s\in\got{t}$, 
and if we substitute these as $\delta 'a$ 
in (\ref{A*sigma}) then we obtain 
\[
\sigma ^{\got{t}}(\delta t,\,\delta s)
+\delta\zeta ((\delta s)_{\got{l}}) 
+(\sigma _{\scriptop{h}})_x(\delta x,\, 
(\delta s_{\scriptop{h}})_{M_{\scriptop{h}}}(x)).
\]
Requiring that this is equal to zero for all 
$\delta s\in\got{t}$ is equivalent to 
$\delta\zeta =0$ (let $\delta s$ range over 
$\got{l}\cap\got{t}_{\scriptop{f}}$), 
$\delta t\in\got{l}$ (let $\delta s$ range over 
$\got{t}_{\scriptop{f}}$ and use that 
$\delta\zeta =0$) and 
$\delta x$ is symplectically orthogonal to 
$\op{T}_x(T_{\scriptop{h}}\cdot x)$. 
If $x\in (M_{\scriptop{h}})_{\scriptop{reg}}$, then 
the last condition implies that $\delta x=Y_{M_{\scriptop{h}}}(x)$ 
for a unique $Y\in\got{t}_{\scriptop{h}}$. 
This shows that the principal orbits of the $T$\--action 
are coisotropic submanifolds of $(M,\,\sigma )$. 

We now verify that the invariants of the compact 
connected symplectic manifold $(M,\,\sigma )$ 
with symplectic $T$\--action with coisotropic principal 
orbits are the ingredients in Definition \ref{ingredientdef} 
we started out with. 

If we substitute $\delta\zeta =\delta '\zeta =0$ 
and $\delta x=\delta 'x=0$ in (\ref{A*sigma}), then we 
get $\sigma ^{\got{t}}(\delta t,\,\delta 't)$, 
which shows that the pull\--back of $\sigma$ to 
the $T$\--orbits is given by $\sigma ^{\got{t}}$. 

If $s\in T$ and $((t,\,\zeta ),\, x)$ are such that 
\[
((s\, t,\,\zeta ),\, x)
=(b,\,\beta )\cdot ((t,\,\zeta ),\, x)
=((t\, b^{-1}\, \,\op{e}^{c(\zeta ,\,\beta )/2},\,\zeta -\beta ),\, 
(b\,\tau _{\beta})_{\scriptop{h}}\cdot x)
\]
for some $(b,\,\beta )\in H$, then $\beta =0$, 
$b=s^{-1}$, and $x={s_{\scriptop{h}}}^{-1}\cdot x$. Because 
$(b,\, 0)\in H$ implies that $(s^{-1})_{\scriptop{f}}
=b_{\scriptop{f}}=1$, 
it follows that $s\in T_{\scriptop{h}}$ and $s\cdot x=x$. This shows 
that $T_{\scriptop{h}}(M,\,\sigma ,\, T)$, the maximal stabilizer 
subgroup of the $T$\--action 
on $M$, see Remark \ref{T1rem}, is equal to $T_{\scriptop{h}}$. 
  
The action of the subtorus $T_{\scriptop{h}}$ of $T$ on 
$M=G\times _HM_{\scriptop{h}}$ is induced by the action 
of $T_{\scriptop{h}}$ on the second factor $M_{\scriptop{h}}$ of 
$G\times M_{\scriptop{h}}$. 
It follows that the action of $T_{\scriptop{h}}$ on $M$ is 
Hamiltonian with image of the momentum mapping equal to a 
translate of $\Delta$. This proves that 
$\Delta (M,\,\sigma ,\, T)$ is equal to a translate of $\Delta$, 
and therefore equal to $\Delta$ if we add a suitable constant 
to the momentum mapping. 

Because $M/T=((T\times N)\times _HM_{\scriptop{h}})/T
\simeq (N/P)\times (M_{\scriptop{h}}/T_{\scriptop{h}})
\simeq (N/P)\times \Delta$, 
we have that $P(M,\,\sigma ,\, T)=P$.  

For each $\zeta \in N$, the infinitesimal action 
of $(0,\,\zeta )\in\got{g}$ on $M$ defines a smooth 
vector field $L_{\zeta}$ on $M$. If the vector fields 
$L_{\scriptop{h},\, \eta}$ on $M_{\scriptop{h}}$
are lifts of $\eta\in C\simeq {\got{t}_{\scriptop{h}}}^*\simeq 
(\got{l}/\got{l}\cap\got{t}_{\scriptop{f}})^*$ 
as in Proposition 
\ref{liftprop} with $(M,\,\sigma ,\, T)$ replaced by 
$(M_{\scriptop{h}},\,\sigma _{\scriptop{h}},\, T_{\scriptop{h}})$, 
then the vector field 
$((0,\, 0),\, L_{\scriptop{h},\,\eta })$ is intertwined by 
$\pi _M$ with a unique vector field $L_{\eta}$ on $M$, 
and the $L_{\eta}$, $L_{\zeta}$ together form a collection 
of lifts of $\eta$, $\zeta$ as in Proposition 
\ref{liftprop}, with $c$ replaced by 
$c(M,\,\sigma ,\, T)$ in Proposition 
\ref{liftprop}, iii). It now follows from (\ref{gotg}) that 
$c(M,\,\sigma ,\, T)=c$.  

Finally, if $\zeta\in P$, then $({\tau _{\zeta}}^{-1},\,\zeta )\in H$, 
and therefore 
\begin{eqnarray*}
(0,\,\zeta )\cdot H\cdot ((0,\, 0),\, x) 
&=&H\cdot ((0,\,\zeta ),\, x)
=H\cdot ((0,\,\zeta )(({\tau _{\zeta}}^{-1},\,\zeta )^{-1},\, x)\\
&=&H\cdot ((\tau _{\zeta},\, 0),\, x)
=\tau _{\zeta}\cdot H\cdot ((0,\, 0,\, x),
\end{eqnarray*}
which implies that $\overline{\tau}(M,\,\sigma ,\, T)
=\overline{\tau}$. 
\end{proof}
According to Theorem \ref{invariantthm}, 
the symplectic manifold $(M,\,\sigma )$ is unique up to 
$T$\--equivariant symplectomorphisms. In particular the 
dimension of $M$ is determined in terms of the 
ingredients in Definition \ref{ingredientdef}.  
Lemma \ref{coisotropicorbitlem} 
implies that $\op{dim}M=\op{dim}T+\op{dim}\got{l}$. 

In the language of categories, see MacLane \cite{maclane}, 
Theorem \ref{invariantthm} and Theorem \ref{existencethm} 
can be summarized as follows. 
\begin{corollary}
Let $T$ be a torus. 
Let ${\cal M}$ denote the category of which 
the objects are the compact connected 
symplectic manifolds 
$(M,\,\sigma )$  together with an effective symplectic 
$T$\--action on $(M,\,\sigma )$ 
with coisotropic principal orbits, and the morphisms are 
the $T$\--equivariant symplectomorphisms. 
Let ${\cal I}$ denote the set of all lists of 
invariants as in Definition \ref{ingredientdef}, 
viewed as a category with only the identities as morphisms. 

Then the assignment $\iota$ in Definition \ref{invariantdef} 
is an equivalence of categories from ${\cal M}$ onto ${\cal I}$. 
In particular it follows that the proper class 
${\cal M}/\!\!\sim$ of isomorphism classes in ${\cal M}$ 
is a set, and the functor $\iota :{\cal M}\to{\cal I}$ 
induces a bijective mapping $\iota/\!\!\sim$ 
from ${\cal M}/\!\!\sim$ onto ${\cal I}$. 
\label{crowncor}
\end{corollary}
\begin{proof}
That $\iota :{\cal M}\to {\cal I}$ is a functor 
and $\iota/\!\!\sim$ is injective 
follows from Theorem \ref{invariantthm}. The surjectivity of $\iota$, 
hence of $\iota /\!\!\sim$,  
follows from Theorem \ref{existencethm}. 
\end{proof}

\begin{remark}
Let $T_{\scriptop{f}}$ be a complementary torus to 
$T_{\scriptop{h}}$ in $T$. 
If $M$ is $T$\--equivariantly diffeomorphic to 
$M_{\scriptop{f}}\times M_{\scriptop{h}}$, 
in which $T_{\scriptop{h}}$ acts only on 
$M_{\scriptop{h}}$ with isolated fixed points, 
and $M_{\scriptop{f}}$ acts freely on 
$M_{\scriptop{f}}$, then $c(N\times N)\subset\got{t}_{\scriptop{f}}$.  
See Proposition \ref{chernprop}. 
Conversely, if this condition is satisfied, then 
Lemma \ref{e0rem} implies the stronger statement that 
$(M,\,\sigma ,\, T)$ is $T$\--equivariantly symplectomorphic 
to the Cartesian product of a symplectic manifold 
$(M_{\scriptop{f}},\,\sigma _{\scriptop{f}},\, T_{\scriptop{f}})$ 
with a free symplectic $T_{\scriptop{f}}$\--action 
and a Delzant manifold 
$(M_{\scriptop{h}},\,\sigma _{\scriptop{h}},\, T_{\scriptop{h}})$.  

Let $c$ in 5) in Definition \ref{ingredientdef} be such that 
$c(\zeta ,\,\zeta ')\notin\got{t}_{\scriptop{f}}$ for some  
$\zeta ,\,\zeta '\in N$. If $(M,\,\sigma ,\, T)$ is 
as in Theorem \ref{existencethm}, then $M$ 
is not $T$\--equivariantly diffeomorphic to 
$M_{\scriptop{f}}\times M_{\scriptop{h}}$, 
in which $T_{\scriptop{f}}$ acts freely on $M_{\scriptop{f}}$
and $T_{\scriptop{h}}$ acts only on $M_{\scriptop{h}}$ and with  
isolated fixed points. Therefore such 
$(M,\,\sigma ,\, T)$ are counterexamples to  
Benoist \cite[Th. 6.6]{benoist}, if in 
\cite[Th. 6.6]{benoist} the word ``isomorphic'' 
implies ``equivariantly diffeomorphic''. 

There exists $c$ in 5) in Definition \ref{ingredientdef},  
such that for every choice of 
a complementary torus $T_{\scriptop{f}}$ to $T_{\scriptop{h}}$ in $T$ 
we have $c(\zeta ,\,\zeta ')\notin\got{t}_{\scriptop{f}}$ for some  
$\zeta ,\,\zeta '\in N$.  
For instance, if $\op{dim}N\geq 2$ and 
$T_{\scriptop{h}}\neq \{ 1\}$,  
then there exists a nonzero antisymmetric 
bilinear mapping $c$ from $N\times N$  
to $\got{t}_{\scriptop{h}}$, which maps $P\times P$ into the integral 
lattice $(T_{\scriptop{h}})_{\Z}$ in the Lie algebra 
$\got{t}_{\scriptop{h}}$ of $T_{\scriptop{h}}$.  
Such a $c$ satisfies 5a) because every $\zeta\in N$ 
is a linear form on $\got{l}$ which vanishes on 
$\got{t}_{\scriptop{h}}$, and it satisfies 5b) by 
assumption. On the other hand $c(\zeta ,\,\zeta ')
\notin\got{t}_{\scriptop{f}}$ as soon as $c(\zeta ,\,\zeta ')\neq 0$. 
\label{productrem}
\end{remark}

\begin{remark}
Let $T_{\scriptop{h}}=\{ 1\}$ and $\sigma ^{\got{t}}=0$, that is, 
the action of $T$ is free, $\op{dim}M=2\op{dim}T$, 
and the orbits are Lagrange submanifolds of $(M,\,\sigma )$. 
In this case the admissible connections are just the smooth $T$\--invariant 
infinitesimal connections for the principal $T$\--bundle 
$\pi :M\to M/T\simeq 
\got{t}^*/P$ over the torus $\got{t}^*/P$, as in 
Remark \ref{curvrem}. The first step in the proof 
of Proposition \ref{liftprop} consists of the construction 
of an infinitesimal connection for the principal $T$\--bundle 
$M$ over the torus $\got{t}^*/P$, of which the curvature form 
is a constant two\--form on the torus $\got{t}^*/P$. 
In this construction the symplectic form did not enter, 
and the principal $T$\--bundle $M$ over $\got{t}^*/P$ 
can be constructed from the ingredients 4), 5) in which 
condition 5a) is kept, but condition 5b) is dropped. 

However, if one has a $T$-invariant 
symplectic form $\sigma$ on $M$ for which the $T$-orbits are 
Lagrangian, then (\ref{ft}), that is 5b), holds. In combination 
with Theorem \ref{existencethm}, we conclude that this principal 
$T$-bundle $M$ over $\got{t}^*/P$ admits a $T$-invariant symplectic form 
for which the $T$-orbits are Lagrangian, if and only if 5b) holds. 
This interpretation of condition 5b) was suggested 
to us by Yael Karshon.  

If $\op{dim}N\geq 3$, then there exist antisymmetric bilinear 
mappings $c:\got{t}^*\times\got{t}^*\to\got{t}$ for which 
5b) does not hold, and it follows that the principal 
$T$\--bundle over $\got{t}^*/P$ defined by $g$ does not admit 
a $T$-invariant symplectic form for which the $T$-orbits are 
Lagrangian. 
\label{cocyclerem}
\end{remark}

\begin{remark}
A slightly different approach to the classification would 
be to allow morphisms $(\Phi ,\,\iota ):
(M,\,\sigma ,\, T)\to (M',\,\sigma ',\, T')$, 
in which $\Phi$ is a symplectomorphism from 
$(M,\, T)$ onto $(M',\,\sigma ')$, 
$\iota$ is an isomorphism of Lie groups from $T$ onto 
$T'$, and $\Phi$ intertwines the $T$\--action on $M$ 
with the $T'$\--action on $M'$ in the sense that 
$\Phi (t\cdot x)=\iota (t)\cdot\Phi (x)$ for 
every $x\in M$ and $t\in T$. 

The isomorphisms between tori are classified by 
the choices of $\Z$\--bases in the integral lattices. 
For instance if we fix a $\Z$\--basis $e_i'$, 
$1\leq i\leq d:=\op{dim}T'=\op{dim}T$, of 
${T'}_{\Z}$, then 
the mapping which assigns to a $\Z$\--basis 
$e_i$, $1\leq i\leq d$, of $T_{\Z}$ the isomorphism 
$\iota :T\to T'$ such that the tangent mapping of 
$\iota$ at the identity element maps $e_i$ to $e_i'$, 
is a bijective mapping from the set of $\Z$\--bases of 
$T_{\Z}$ onto the set of isomorphisms from $T$ onto $T'$. 
In turn the set of $\Z$\--bases in $T_{\Z}$ 
is in bijective correspondence with the group 
$\op{GL}(d,\,\Z )$ of all $d\times d$\--matrices 
with integral coefficients and determinant equal 
to $\pm 1$. 

This can be applied in particular to $T'=\R ^d/\Z ^d$, 
with $e'_i$, $1\leq i\leq d$, equal to the standard basis 
of $\R ^d$. If we also choose a $\Z$\--basis 
$\varepsilon ^l$, $1\leq l\leq\op{dim}N$, in $P$, 
then the ingredients 1) -- 6) in Definition 
\ref{ingredientdef} are determined by their coefficients 
with respect to these bases. Furthermore the 
groups $G$ and $H$ are identified with 
$\R ^d/\Z ^d\times\R ^{d_N}$ and $\R ^{d_{\scriptop{h}}}
/\Z ^{d_{\scriptop{h}}}\times\Z ^{d_N}$, respectively.  
This would lead to a presentation of the model in coordinates, except 
for the Delzant manifold $(M_{\scriptop{h}},\,
\sigma _{\scriptop{h}},\, T_{\scriptop{h}})$. 
Such a model looks even more 
explicit than the one in Proposition \ref{Aprop}. 

The disadvantage of this approach is that the invariants 
are given by coefficient matrices, which are 
determined uniquely only up to the action on these 
matrices of the changes of $\Z$\--bases. 
Also the notations become 
quite a bit heavier if we write out our objects in coordinates.  
\label{baserem}
\end{remark}

\section{$V$\--parallel spaces}
\label{Vsec}
In this section we define the notion of a $V$\--parallel space, and 
prove that every complete, connected and locally convex 
$V$\--parallel space is isomorphic to the Cartesian product of a closed 
convex subset of a vector space and a torus. 
\begin{definition}
Let $V$ be an $n$\--dimensional vector space. 
A {\em $V$\--parallel space} is a Hausdorff 
topological space $Q$, together with an open covering $Q_{\alpha}$, 
$\alpha\in A$, of $Q$ and homeomorphisms $\varphi_{\alpha}$ from 
$Q_{\alpha}$ onto subsets $V_{\alpha}$ of $V$ such that, for every 
$\alpha ,\,\beta\in A$ for which $Q_{\alpha}\cap Q_{\beta}\neq\emptyset$, 
the mapping 
\begin{equation}
x\mapsto \varphi_{\alpha}(x)-\varphi_{\beta}(x)
:Q_{\alpha}\cap Q_{\beta}\to V
\label{pareq}
\end{equation}
is locally constant. 
A subset $U$ of $Q$ is a $V$\--parallel space with the 
$\varphi _{\alpha}$ replaced by their restrictions to 
$U\cap Q_{\alpha}$. If $W$ is another vector space and 
$R$ is another $W$\--parallel space, then 
$Q\times R$ is a $V\times W$\--parallel space in an obvious way. 

The $V$-parallel space $Q$ is called 
{\em locally convex} if the $V_{\alpha}$ are convex subsets of $V$. 
The locally convex $V$\--parallel space $Q$ is called 
{\em locally convex polyhedral} 
if for every $\alpha\in A$ there is a convex open 
subset $V_{\alpha}'$ of $V$ and there are finitely many 
linear forms $v^*_{\alpha,\, i}$ $1\leq i\leq m$, on $V$, 
such that 
\begin{equation}
V_{\alpha}=\{ v\in V_{\alpha}'\mid
v^*_{\alpha,\, j}(v)\geq 0\; 
\mbox{\rm for every}\; 1\leq j\leq m\} .
\label{cornereq}
\end{equation}
\label{cornerdef}
\end{definition}
\begin{remark}
If we have (\ref{cornereq}) with linearly independent 
linear forms $v^*_{\alpha,\, j}$, and (\ref{pareq}) 
is replaced to the weaker condition that 
for every $\alpha,\,\beta\in A$ the mapping 
$\varphi _{\alpha}\circ\varphi_{\beta}^{-1}$ 
is a diffeomorphism, then $Q$ is a ``manifold with corners'' as defined 
for instance in Mather \cite[\S 1]{mather}. If the $V_{\alpha}$ 
are open subsets of $V$  
and (\ref{pareq}) is relaxed to the condition that for every 
$\alpha ,\,\beta\in A$ the mapping $\varphi _{\alpha}\circ\varphi_{\beta}^{-1}$ 
is locally an affine linear transformation, the composition of a 
linear mapping and a translation, then $Q$ is a 
``manifold with affine covering'' as in 
Auslander and Markus \cite[p. 141]{am} or a 
``locally affine manifold'' as in Auslander \cite{auslander}. 

If on the other hand we have (\ref{pareq}) and (\ref{cornereq}) 
with $V_{\alpha}=V$ for all 
$\alpha\in A$, then $Q$ is called an ``affine space modelled over $V$'' in 
geometry. In order to avoid confusion about the 
interpretation of the term ``affine'', we have replaced the  
term ``affine'' by ``parallel'', where the latter word 
reminds of a parallelism, which is 
defined as a global frame in the tangent bundle, or 
equivalently a trivialization of the tangent bundle. 
\label{cornerdefrem}
\end{remark}

\begin{definition}
Let $V$ and $W$ be finite\--dimensional vector spaces, 
$Q$ and $R$ a locally convex $V$\--parallel and $W$\--parallel spaces 
with charts $\varphi _{\alpha}$, $\alpha\in A$ and 
$\psi _{\beta}$, $\beta\in B$, respectively. 
If $L$ is a linear mapping from $W$ to $V$, then 
an {\em $L$\--map from $R$ to $Q$} is a continuous 
map $f:R\to Q$ such that for each  $\beta\in B$ and $\alpha\in A$ 
and each $U$ of $Q_{\alpha}\cap f^{-1}(R_{\beta})$ 
such that $\varphi _{\alpha}(U)$ is convex, 
we have that 
\begin{equation}
\varphi _{\alpha}(f(p))-\varphi _{\alpha}(f(q))
=L\, (\psi _{\beta}(p)-\psi _{\beta}(q)) 
\quad\mbox{\rm for all}\quad p,\, q\in U. 
\label{Lmapeq}
\end{equation}

With such maps as morphisms, the 
locally convex parallel spaces form a 
category. In particular two locally convex parallel spaces are 
called {\em isomorphic} if there exists an $L$\--map 
from the first one to the second one, which is a homeomorphism and 
for which the linear mapping $L$ is bijective. 
\label{Lmapdef}
\end{definition}
\begin{definition}
Let $Q$ be a $V$\--parallel space. For each neighborhood 
$V_0$ of $0$ in $V$, let $P(V_0)$ be the set of 
all $(p,\, q)\in Q\times Q$ such that there exists an 
$\alpha\in A$ for which $p,\, q\in Q_{\alpha}$ and 
$\varphi _{\alpha}(p)-\varphi _{\alpha}(q)\in V_0$. 
The $P(V_0)$ define a {\em uniform structure} on $Q$ in the sense 
of Bourbaki \cite[Ch. II, \S 1, No. 1, D\'ef. 1]{bourbaki}. 
With respect to this uniform structure, we have the notion 
of {\em Cauchy filters} in $Q$ as in Bourbaki 
\cite[Ch. II, \S 3, No. 1, D\'ef. 2]{bourbaki}. 
The $V$\--parallel space is called {\em complete} if 
every Cauchy filter in $Q$ converges to an element of $Q$, 
cf. Bourbaki \cite[Ch. II, \S 3, No. 3, D\'ef. 3]{bourbaki}.
\label{uniformdef} 
\end{definition}
Our next goal is to prove that every complete, connected 
and locally convex  $V$\--parallel space is 
isomorphic to the product of a closed convex subset of 
a vector space and a torus. See Theorem \ref{cornerclassthm}
below for the precise statement. 

Let $Q$ be a locally convex $V$\--parallel 
space, as in Definition \ref{cornerdef}, and let $v\in V$. 
A {\em motion in $Q$ with constant velocity $v$} is a continuous mapping 
$\gamma$ from an interval $I$ in $\R$ to $Q$ such that 
for every $\alpha\in A$ and interval $J$ in $\gamma ^{-1}(Q_{\alpha })$, 
we have 
\[
\varphi _{\alpha}(\gamma (t))=\varphi _{\alpha}(\gamma (s))
+(t-s)\, v\quad\mbox{\rm for all}\quad s,\, t\in J. 
\]
In other words, 
$\gamma :I\to Q$ is an $L$\--map, where $L$ is the 
linear mapping from $\R$ to $V$ which sends $r\in\R$ 
to $r\,v\in V$. 
\begin{lemma}
If $\gamma :I\to Q$ and $\delta :J\to Q$ are motions in $Q$ 
with constant velocity $v$, and $\gamma (s)=\delta (s)$ for some 
$s\in I\cap J$, then $\gamma (t)=\delta (t)$ for all 
$t\in I\cap J$. 
\label{uniquemotionlem}
\end{lemma}
\begin{proof}
Let $K=\{ t\in I\cap J\mid\gamma (t)=\delta (t)\}$. 
Then $K$ is a closed subset of $I\cap J$, because $\gamma$ 
and $\delta$ are continuous. 

Suppose $s\in K$. Because the $Q_{\alpha}$, $\alpha\in A$, 
form a covering of $Q$, there exists an $\alpha\in A$ such that 
$\gamma (s)=\delta (s)\in Q_{\alpha}$. Because $Q_{\alpha}$ 
is open in $Q$ and $\gamma$ and $\delta$ are continuous, 
$H:=\gamma ^{-1}(Q_{\alpha})\cap\delta ^{-1}(Q_{\alpha})$ 
is an open neighborhood of $s$ in $I\cap J$. For every 
$t\in H$ we have 
\[
\varphi _{\alpha}(\gamma (t))
=\varphi _{\alpha}(\gamma (s))+(t-s)\, v
=\varphi _{\alpha}(\delta (s))+(t-s)\, v
=\varphi _{\alpha}(\delta (t)),
\]
hence $\gamma (t)=\delta (t)$ because $\varphi _{\alpha}$ 
is injective. Therefore $H\subset K$, and we have proved that 
$K$ is also an open subset of $I\cap J$. Because 
$I\cap J$ is connected, it follows that $K=\emptyset$ 
or $K=I\cap J$.   
\end{proof}

Let $s\in\R$ and $p\in Q$. Let $\Gamma _{s,\, p}^v$ denote the 
set of all motions $\gamma$ in $Q$ with constant velocity $v$, 
which are defined on an interval $I_{\gamma}$ in $\R$ such that 
$s\in I_{\gamma}$ and $\gamma (s)=p$. Then it follows from 
Lemma \ref{uniquemotionlem} that the $\gamma\in\Gamma _{s,\, p}^v$ 
have a common extension $\gamma _{s,\, p}^v$ to the union 
$I_{s,\, p}^v$ of all the intervals $I_{\gamma}$, 
$\gamma\in\Gamma _{s,\, p}^v$. $I_{s,\, p}^v$ is an interval 
in $\R$ and $\gamma _{s,\, p}^v:I_{s,\, p}^v\to Q$ is a motion 
in $Q$ with constant velocity $v$, the unique maximal motion 
$\gamma :I\to Q$ with constant velocity $v$ such that $s\in I$ 
and $\gamma (s)=p$. 
\begin{definition}
$D$ is the set of all $(v, p) \in V\times Q$ 
such that $1 \in I_{0,\, p}^v$. We write 
$p+v=\gamma _{0,\, p}^v(1)$ when $(v, p)\in D$. 
Note that 
$p+v\in Q$ when $p\in Q$ and $v\in V$, whereas 
$v+w\in V$ when $v,\, w\in V$. 
For every $p\in Q$ we write 
$D_p=\{ v\in V\mid (v,\, p)\in D\}$. 
\label{Ddef}
\end{definition}

\begin{lemma}
Assume that the locally convex 
$V$\--parallel space $Q$ is complete, 
see Definition \ref{cornerdef} and Definition 
\ref{uniformdef}.  Then, insofar 
as defined, the mapping $(v,\, p)\mapsto p + v$ is an action of 
the additive group $(V,\, +)$ on $Q$, in the following sense.  
\begin{itemize}
\item[i)] For every $p\in Q$ we have $(0,\, p)\in D$ and $p+0=p$. 
\item[ii)] If $(v ,\, p)\in D$ then $(-v ,\, p+v )\in D$ and 
$(p+v)+(-v) = p$. 
\item[iii)] If $(v_1,\, p)\in D$ and $(v_2,\, p +v_1)\in D$, then 
$(v_1+ v_2,\, p)\in D$ and $(p + v_1) + v_2 = p + (v_1+ v_2)$. 
\end{itemize}
\label{actionlem}
\end{lemma}
\begin{proof}
The statements i) and ii) follow immediately from the definitions. 
Our proof of iii) is surprisingly long. 

In order to prove iii), assume that 
$(v_1,\, p)\in D$ and $(v_2,\, p+v_1)\in D$. 

Let $L$ be the linear mapping from $\R ^2$ to $V$ which sends 
$e_1=(1,\, 0)$ to $v_1$ and $e_2=(0,\, 1)$ to $v_2$. 
For any $s_1,\, s_2\in\R$, let $C(s_1,\, s_2)$ 
denote the convex hull of $(s_1,\, 0)$, $(1,\, 0)$, 
and $(1,\, s_2)$ in $\R ^2$, which is a solid triangle. 

Let $I$ denote 
the set of all $s_1\in [0,\, 1]$ for which there exists 
$0<s_2\leq 1$ and an $L$\--map $f:C(s_1,\, s_2)\to Q$ 
such that $f(s_1,\, 0)=p+s_1\, v_1$. 
Note that if $s_1\in I$, then $[s_1,\, 1]\subset I$. 
Note also that for every 
$(t_1,\, t_2),\, (t_1',\, t_2')\in C(s_1,\, s_2)$, 
\[
[0,\, 1]\ni u\mapsto f(t_1+u\, (t_1'-t_1),\, 
t_2+u\, (t_2'-t_2))
\]
is a motion in $Q$ with constant velocity 
$L(t_1'-t_1,\, t_2'-t_2)=(t_1'-t_1)\, v_1+(t_2'-t_2)\, v_2$, 
hence 
$((t_1'-t_1)\, v_1+(t_2'-t_2)\, v_2,\, f(t_1,\, t_2))\in D$ 
and 
\[
f(t_1',\, t_2')=
f(t_1,\, t_2)+((t_1'-t_1)\, v_1+(t_2'-t_2)\, v_2).
\]
If we apply this with $(t_1,\, t_2)=(s_1,\, 0)$, then we see 
that $f$ is uniquely determined by the formula 
\begin{equation}
f(t_1,\, t_2)=(p+s_1\, v_1)+((t_1-s_1)\, v_1+t_2\, v_2), 
\label{f+}
\end{equation}
where 
$((t_1-s_1)\, v_1+t_2\, v_2,\, p+s_1\, v)\in D$ for every 
for every $(t_1,\, t_2)\in C(s_1,\, s_2)$, which in turn implies that 
\begin{eqnarray}
(p+s_1\, v_1)+((t_1'-s_1)\, v_1+t_2'\, v_2)
&=&((p+s_1\, v_1)+((t_1-s_1)\, v_1+t_2\, v_2))
\nonumber\\
&&+((t_1'-t_1)\, v_1+(t_2'-t_2)\, v_2), 
\label{+}
\end{eqnarray}
where $((t_1-s_1)\, v_1+t_2\, v_2,\, p+s_1\, v)$, 
$((t_1'-s_1)\, v_1+t_2'\, v_2,\, p+s_1\, v)$ 
and $((t_1'-t_1)\, v_1+(t_2'-t_2)\, v_2,\, 
(p+s_1\, v_1)+((t_1-s_1)\, v_1+t_2\, v_2)$ 
all belong to $D$ 
for every $(t_1,\, t_2),\,(t_1',\, t_2')\in C(s_1,\, s_2)$. 

Let $i_1$ denote the infimum of $I$, which implies that 
$\left] i_1,\, 1\right]\subset I$. We will show that this 
implies that $i_1\in I$, which means that $I=[i_1,\, 1]$. 
Because trivially $1\in I$, we may assume that $0\leq i_1<1$. 
Because $(i_1\, v_1,\ p)\in D$, there exists an $\alpha\in A$ such that 
$p+i_1\, v_1\in Q_{\alpha}$. Because the mapping 
$t\mapsto p+t\, v$ is continuous from $I^{v_1}_{0,\, p}$ to $Q$, 
there exists an $s_1\in\left] i_1,\, 1\right]\subset I$ such that 
$p+s_1\, v_1\in Q_{\alpha}$. Because $s_1\in I$, there exists 
$0<s_2\leq 1$ and an $L$\--map $f:C(s_1,\, s_2)\to Q$ 
such that $f(s_1,\, 0)=p+s_1\, v_1$. Note that 
for each $u\in [0,\, 1]$ we have $c(u):=(s_1,\, 0)
+u\, ((1,\, s_2)-(s_1,\, 0))\in C$. 
Because $f$ is continuous, 
there exists $0<u\leq 1$ such that 
$f(c)\in Q_{\alpha}$, if we write $c:=c(u)$. 
Because $p+i_1\, v_1$, $p+s_1\, v_1$, and 
$f(c)$ all belong to $Q_{\alpha}$, the points 
$\varphi _{\alpha}(p+i_1\, v_1)$, 
$\varphi _{\alpha}(p+s_1\, v_1)$, and 
$\varphi _{\alpha}(f(c))$ all belong to 
the convex subset $V_{\alpha}$ of $V$, which implies that 
their convex hull $B_{\alpha}$ in $V$ is contained in $V_{\alpha}$.   
Let $B$ be the convex hull 
of $(i_1,\, 0)$, $(s_1,\, 0)$, and $c$ in $\R ^2$. 
The $L$\--map from $B$ onto $B_{\alpha}$ 
which sends $(s_1,\, 0)$ to $\varphi _{\alpha}(p+s\, v_1)$, 
followed by ${\varphi_{\alpha}}^{-1}$, defines an $L$\--map $e$ 
from $B$ to $Q$ such that $e(s_1,\, 0)=p+s_1\, v_1$. 

Because of the uniqueness of 
$L$\--maps which map $(s_1,\, 0)$ to $p+s_1\, v_1$, 
see (\ref{f+}), we have that $e=f$ 
on $B\cap C(s_1,\, s_2)$, and therefore $e$ and $f$ have a common extension 
$g:B\cup C(s_1,\, s_2)\to Q$. In order to prove that $g$ is an $L$\--map, 
we observe that the property of being an $L$\--map is local, 
and because $e$ and $f$ are $L$\--maps on the open subsets 
$B\setminus C(s_1,\, s_2)$ and $C(s_1,\, s_2)\setminus B$ of 
$B$ and $C(s_1,\, s_2)$, respectively,  
we have that $g$ is an $L$\--map on 
$(B\setminus C(s_1,\, s_2))\cup (C(s_1,\, s_2)\setminus B)
=(B\cup C(s_1,\, s_2))\setminus (B\cap C(s_1,\, s_2))$. 
On the other hand, if $r\in B\cap C(s_1,\, s_2)$, then there are 
neighborhoods $B_0$ and $C_0$ of $r$ in $B$ and  
$C(s_1,\, s_2)$, respectively, 
such that $\varphi _{\alpha}(e(p))-\varphi _{\alpha}(e(r))=L\, (p-r)$ 
when $p\in B_0\cap e^{-1}(Q_{\alpha})$ 
and $\varphi _{\alpha}(f(q))-\varphi _{\alpha}(f(r))=L\, (p-r)$ 
when $q\in C_0\cap f^{-1}(Q_{\alpha})$. It follows that 
\begin{eqnarray*}
\varphi _{\alpha}(g(p))-\varphi _{\alpha}(g(q))
&=&(\varphi _{\alpha}(g(p))-\varphi _{\alpha}(g(r)))
+(\varphi _{\alpha}(g(r))-\varphi _{\alpha}(g(q)))\\
&=&(\varphi _{\alpha}(e(p))-\varphi _{\alpha}(e(r)))
+(\varphi _{\alpha}(f(r))-\varphi _{\alpha}(f(q)))\\   
&=&L\, (p-r)+L\, (r-q)=L\, (p-q),
\end{eqnarray*}
which implies (\ref{Lmapeq}) with $f$ replaced by $g$. 

Let $d$ be the intersection point of the 
straight line through $(1,\, 0)$ and $(1,\, s_2)$, 
and the straight line through $(i_1,\, 0)$ and 
$c$. A straightforward calculation shows that 
$d=(1,\, s_2')$, with 
\[
s_2':=\frac{u\, (s_1-i_1)+u\, (1-s_1)}{s_1-i_1+u\, (1-s_1)}\, s_2
\leq s_2.  
\] 

\begin{figure}[ht]
\begin{center}
\begin{psfrags}
\epsfxsize=5cm\leavevmode\epsfbox{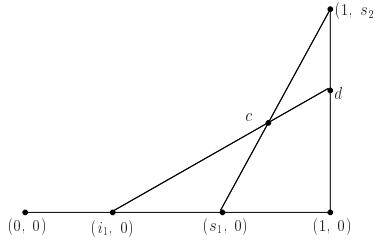}
\end{psfrags}
\end{center}
\caption{The union of $C(s_1,\, s_2)$ and $C(i_1,\, s_2')$.} 
\label{auniffig}
\end{figure}

Because $c$ is lying on the straight line between 
$(i_1,\, 0)$ and $d$, it follows that the convex hull 
$C(i_1,\, s_2')$ of $(i_1,\, 0)$, $(1,\, 0)$ and $d$ 
is equal to to the union of the convex hull 
$B$ of $(i_1,\, 0)$, $(s_1,\, 0)$, $c$, and the 
convex hull $F$ of the points 
$(i_1,\, 0)$, $c$, 
$(1,\, 0)$, and $d$. Because the latter four points 
all lie in $B\cup C(s_1,\, s_2)$, it follows that 
$C(i_1,\, s_2')\subset B\cup C(s_1,\, s_2)$.  
Because the restriction of $g$ to $S(i_1,\, s_2')$ is an 
$L$\--map such that $g(s_1,\, 0)=p+s_1\, v_1$, hence 
\[
g(i_1,\, 0)=(p+s_1\, v_1)+(i_1-s_1)\, v_1=p+i_1\, v_1
\]
in view of (\ref{f+}) with $f$ replaced by $g$, 
it follows that $i_1\in I$.

Now suppose that $i_1>0$. With $\alpha\in A$ such that 
$p+i_1\, v_1\in Q_{\alpha}$, there exists 
$0\leq s_1<i_1$ such that $p+s_1\, v_1\in Q_{\alpha}$. 
The same reasoning as above with $i_1$ and $s_1$ interchanged, 
where we use that $i_1\in I$,  
leads to the conclusion that $s_1\in I$, in contradiction 
with the definition $i_1=\op{inf}I$ of $i_1$.
We conclude that $i_1=0$, or $0\in I$, which 
means that there exists 
$0<s_2\leq 1$ and an $L$\--map $f$ from 
$C(0,\, s_2)$ to $Q$ 
such that $f(0,\, 0)=p$. 

Define $J$ as the set of all $s_2\in [0,\, 1]$ for which there exists 
an $L$\--map $f$ from $C(0,\, s_2)$ to $Q$ 
such that $f(0,\, 0)=p$, and write $s:=\op{sup}J$. 
The uniqueness of the $L$\--maps $f$ from $C(0,\, s_2)$ to $Q$ 
such that $f(0,\, 0)=p$, see (\ref{f+}), where $s_2$ ranges over $J$, 
implies that these $L$\--maps have a common extension to the 
union $C$ over all $s_2\in J$ of the triangles $C(0,\, s_2)$, 
which we also denote by $f$. Note that the closure of 
$C$ in $\R ^2$ is equal to $C(0, s)$. $f:C\to Q$ is an $L$\--map, 
hence uniformly continuous. Because a uniformly continuous 
mapping from a dense subset $A$ of a complete space $X$ to a complete 
space $X'$ has a unique extension to a continuous mapping 
from $X$ to $X'$, see Bourbaki 
\cite[Ch. 2, \S 3, No. 6, Th. 2]{bourbaki}, and because we had 
assumed that $Q$ is complete, it follows that 
$f$ has a unique extension to a continuous mapping 
$\overline{f}:C(0,\, s)\to Q$. Because the restriction of the continuous  
mapping $\overline{f}:C(0,\, s)\to Q$ to the dense subset $C$ of 
$C(0,\, s)$ is an 
$L$\--map, it follows that $\overline{f}$ is an $L$\--map, and 
we conclude that $s\in J$, or equivalently 
$J=[0,\, s]$. 

If $s<1$ then the previous argument leading to $i_1=0$, with 
$v_1$ and $v_2$ replaced by $v_1+s\, v_2$ and $(1-s)\, v_2$, 
respectively, shows that there exists $s<s'\leq 1$ and an 
$L$\--map $f'$ from the convex hull $C'$ of $(0,\, 0)$, $(1,\, s)$, 
and $(1,\, s')$ to $Q$, such that $f'(0,\, 0)=p$. The uniqueness of 
$L$\--maps which send $(0,\, 0)$ to $p$, see (\ref{f+}), 
yields that $f'=f$ on $C'\cap C(0,\, s)$, which implies that 
$f$ and $f'$ have a common extension $f''$ to $C(0,\, s)\cup C'=C$. 
As in the argument leading to $i_1=0$, we have that $f''$ is 
an $L$\--map, and because $f''(0,\, 0)=p$ it follows that 
$s'\in J$, in contradiction with the definition $s=\op{sup}J$ of $s$. 

We arrive at $s=1$, or equivalently $1\in J$, which means that 
there is an $L$\--map $f:C(0,\, 1)\to Q$ such that $f(0,\, 0)=p$. 
Now (\ref{+}) with $(s_1,\, s_2)=(0,\, 1)$, 
$(t_1,\, t_2)=(1,\, 0)$ and 
$(t_1',\, t_2')=(1,\, 1)$ implies that 
$(v_1+v_2,\, p)\in D$ and $p+(v_1+v_2)=(p+v_1)+v_2$. 
This completes the proof of the statement iii) in the lemma. 
\end{proof}

\begin{lemma}
Let $Q$ be a complete locally convex $V$\--parallel space. 
Then $D$ is a closed subset of $V\times Q$  
and $(v,\, p)\mapsto p+v$ is a uniformly continuous $L$\--map 
from $D$ to $Q$. Here $L:V\times V\to V$ is the mapping  
which sends $(v,\, w)\in V\times V$ to $v+w\in V$. 
For each $p\in Q$, the mapping $v\mapsto p+v$ is a local 
homeomorphism from $D_p$ onto an open subset of $Q$. 
\label{contactionlem}
\end{lemma}
\begin{proof}
Let $p\in Q$, and let $D_p$ be defined as in 
Definition \ref{Ddef}. If $\alpha\in A$ and $Q_{\alpha}$ is a 
neighborhood of $p$ in $Q$, then the translate 
$V_{\alpha}-\varphi _{\alpha}(p)$ of $V_{\alpha}$ 
over the vector $-\varphi _{\alpha}(p)$ is 
contained in $D_p$.  
Indeed,  
if $v\in V$, $\varphi _{\alpha}(p)+v\in V_{\alpha}$, then 
the convexity of $V_{\alpha}$ implies that 
$\varphi _{\alpha}(p)+t\, v\in V_{\alpha}$ for every 
$t\in [0,\, 1]$. The continuity of 
${\varphi _{\alpha}}^{-1}$ implies that 
$t\mapsto {\varphi _{\alpha}}^{-1}(\varphi _{\alpha}(p)+t\, v)$ 
is a motion in $Q$ which has constant velocity $v$ and 
is equal to $p$ at $t=0$, which shows that 
$(v,\, p)\in D$. Note that in passing we have proved 
that $\varphi _{\alpha}(p+t\, v)=\varphi _{\alpha}(p)+t\, v$ 
for all $0\leq t\leq 1$, and in particular 
$\varphi _{\alpha}(p+v)=\varphi _{\alpha}(p)+v$

We next claim that $V_{\alpha}-\varphi _{\alpha}(p)$  
is a neighborhood of $0$ in $D_p$.  
If this would not be the case, then there is a sequence 
$v_j$ in $D_p\setminus (V_{\alpha}-\varphi _{\alpha}(p))$ which converges 
to $0$ in $V$ as $j\to\infty$. Because $Q$ is a Hausdorff space, 
there exists an open neighborhood $Q_p$ of $p$ in $Q$ such that the closure 
$\overline{Q_p}$ of $Q_p$ in $Q$ is contained in $Q_{\alpha}$. 
Let $I_j$ be the set of all 
$t\in I_{0,\, p}^{v_j}$ such that $t\geq 0$ and $p+t\, v_j\notin Q_p$,  
and write $t_j:=\op{inf}I_j$. Note that $t_j\geq 0$. 
Because $t\mapsto p+t\, v_j$ 
is continuous from $I^{v_j}_{0,\, p}$ to $Q$, and 
$Q\setminus Q_p$ is closed in $Q$, we have that 
$I_j$ is a closed subset of $I_{0,\, p}^{v_j}$, hence 
$t_j\in I_j$, which implies that $p+t_j\, v_j\notin Q_p$. 
On the other hand we have for every $0\leq t< t_j$ 
that $p+t\, v_j\in Q_p\subset Q_{\alpha}$, hence 
$\varphi _{\alpha}(p+t\, v_j)=
\varphi _{\alpha}(p)+t\, v_j\in V_{\alpha}$. 
Because $v_j\notin D_p\setminus (V_{\alpha}-\varphi _{\alpha}(p))$, 
this cannot happen for $t=1$, which proves that $t_j\leq 1$. 
Because $t\mapsto p+t\, v_j=\gamma _{0,\, p}^{v_j}$ is continuous and 
$p+t\, v_j\in Q_p$ for every $t<t_j$, we have that 
$p+t_j\, v_j\in\overline{Q_p}\subset Q_{\alpha}$, 
hence $\varphi _{\alpha}(p+t_j\, v_j)=\varphi _{\alpha}(p)+t_j\, v_j
\to \varphi _{\alpha}(p)$ in $V$, hence in $V_{\alpha}$. 
Because $\varphi _{\alpha}^{-1}$ is continuous, we conclude 
that $p+t_j\, v_j\to p$ in $Q_{\alpha}$, hence in $Q$, 
in contradiction with the fact that $p+t_j\, v_j\notin Q_p$ 
for every $j$. 

Because $\varphi _{\alpha}(p+v)=\varphi _{\alpha}(p)+v$
for every $v\in V_{\alpha}-\varphi _{\alpha}(p)$,   
the restriction to the neighborhood $V_{\alpha}-\varphi _{\alpha}(p)$ 
of $0$ in $D_p$ of the mapping $v\mapsto p+v$ is equal 
to the continuous mapping 
$v\mapsto {\varphi _{\alpha}}^{-1}(\varphi _{\alpha}(p)+v)$, 
from $V_{\alpha}-\varphi _{\alpha}(p)$ onto the neighborhood 
$Q_{\alpha}$ of $p$ in $Q$, with inverse equal to the 
continuous mapping $q\mapsto\varphi _{\alpha}(q)-\varphi_{\alpha}(p)$ 
from $Q_{\alpha}$ onto $V_{\alpha}-\varphi _{\alpha}(p)$. 
In particular we have that for every $r\in Q$ the mapping 
$z\mapsto r+z$ is a homeomorphism from a neighborhood 
of $0$ in $D_q$ onto a neighborhood of $q$ in $Q$. 

Now let $(v,\, p)\in D$. If $(w,\, q)\in D$ 
such that $q\in Q_{\alpha}$, then $q=p+u$ with 
\[
u:=\varphi _{\alpha}(q)-\varphi _{\alpha}(p)
\in V_{\alpha}-\varphi _{\alpha}(p)
\subset D_p,
\] 
hence $(u+w,\, p)\in D$ and 
$q+w=(p+u)+w=p+(u+w)$ in view of iii) in Lemma \ref{actionlem}. 
Furthermore ii) in Lemma \ref{actionlem} 
implies that $(-v,\, p+v)\in D$ and $p=(p+v)+(-v)$. 
We now have $(-v,\, p+v)\in D$ and $(u+w,\, (p+v)+(-v))\in D$, 
which in view of iii) in Lemma \ref{actionlem} implies that 
$((-v)+(u+w),\, p+v)\in D$ and 
$((p+v)+(-v))+(u+w)=(p+v)+((-v)+(u+w))$. Combining the 
equations we arrive at 
\begin{eqnarray}
q+w&=&p+(u+w)=((p+v)+(-v))+(u+w)=(p+v)+((-v)+(u+w))
\nonumber\\
&=&(p+v)+(-v+\varphi_{\alpha}(q)-\varphi_{\alpha}(p)+w).
\label{qw}
\end{eqnarray}
Because $(w,\, q)\mapsto z=-v+\varphi_{\alpha}(q)-\varphi_{\alpha}(p)+w$ 
is continuous from a neighborhood of $(v,\, p)$ in $D$ to $V$, 
with value $0$ at $(w,\, q)=(v,\, p)$, 
and $z\mapsto r+z$ is continuous from a neighborhood of 
$0$ in $D_r$ to $Q$, with $r=p+v$, we conclude 
that $(w,\, q)\mapsto q+w$ is continuous from a neighborhood 
of $(v,\, p)$ in $D$ to $Q$. Because this holds for every 
$(v,\, p)\in D$, we have proved that the mapping 
$(v,\, p)\mapsto p+v$ is continuous from $D$ to $Q$, 
and therefore it is an $L$\--map from $D$ to $Q$. 
Note also that (\ref{qw}) with $q=p$ implies that 
$v\mapsto p+v$ is a local homeomorphism from $D_p$ onto 
an open subset of $Q$. 

As an $L$\--mapping, the mapping $f:(t,\, v,\,  p)\mapsto p+t\, v$ is uniformly 
continuous from $[0,\, 1]\times D$ to $Q$. If 
$\overline{D}$ denotes the closure of $D$ in the complete space 
$V\times Q$, then $[0,\, 1]\times D$ is dense in the 
complete space $[0,\, 1]\times\overline{D}$.  
Because also $Q$ is complete, 
the mapping $f$  
has a unique extension to 
a continuous mapping $g$ from $[0,\, 1]\times \overline{D}$ 
to $Q$. See Bourbaki 
\cite[Ch. 2, \S 3, No. 6, Th. 2]{bourbaki}. 
By continuity it follows that for each 
$(v,\, p)\in\overline{D}$ the mapping $[0,\, 1]\ni t\mapsto  
g(t,\, v,\, p)$ is a motion in $Q$ with constant velocity $v$, 
equal to $p$ at $t=0$, and the conclusion is that $(v,\, p)\in D$. 
This shows that $\overline{D}=D$, that is, $D$ is closed in $V\times Q$.  
\end{proof}

\begin{lemma}
Let $Q$ be a complete locally convex $V$\--parallel space which 
in addition is connected. Then the action of $V$ on $Q$ is transitive, 
in the sense that for any $p,\, q\in Q$ there exists $v\in V$ 
such that $(v ,\, p)\in D$ and $p+v=q$. 

This implies that for every $p\in Q$ the mapping 
$\Phi _p:v\mapsto p+v$ is a local homeomorphism from 
$D_p$, which is a closed subset of $V$, onto $Q$. 
\label{surjlem}
\end{lemma}
\begin{proof}
We write $p\sim q$ if there exists a $v\in V$ such that 
$(v ,\, p)\in D$ and $p+v= q$. It follows from 
Lemma \ref{actionlem} that $\sim$ is an equivalence relation 
in $Q$. It follows from the last statement in 
Lemma \ref{contactionlem} 
that the equivalence classes 
are open subsets of $Q$. Because $Q$ is connected, the equivalence relation 
$\sim$ has only one equivalence 
class $Q$, which proves the transitivity of the action. 

$D_p$ is a closed subset of $V$, as the preimage of the 
closed subset $D$ of $V\times Q$ under the continuous mapping 
$v\mapsto (v,\, p)$. 
\end{proof}

\begin{lemma}
Let $Q$ be a complete and connected locally convex $V$\--parallel space. 
For any $v\in V$, the following conditions are equivalent. 
\begin{itemize}
\item[a)] There exists a $p\in Q$ such that $I^v_{0,\, p}=\R$. 
\item[b)] For every $p\in Q$ we have $I^v_{0,\, p}=\R$. 
\end{itemize}
Let $N$ denote the set of all $v\in V$ such that a) or b) holds.  
Then $N$ is a linear subspace of $V$, $(N\times Q)\subset D$,  
and the restriction to $N\times Q$ 
of the mapping $(v,\, p)\mapsto p+v$ defines an action 
of the additive group $N$ on $Q$. 

\label{Nlem}
\end{lemma}
\begin{proof}
Assume that a) holds, 
which implies that $(t\, v,\, p)\in D$ for every $t\in\R$. 
Let $q\in Q$. Because of the transitivity in Lemma \ref{surjlem}, 
there exists a $w\in V$ such that $(w,\, q)\in D$ and 
$q+w=p$. Because of iii) in Lemma \ref{actionlem} 
we have for every $t\in\R$ that $(w+t\, v,\, q)\in D$, 
which in view of Definition \ref{Ddef}, the definition of $D$, 
implies that $(s\, (w+t\, v),\, q)\in D$ 
for every $0\leq s\leq 1$. If for any $t\in\R$ and 
$0<s\leq 1$ we replace $t$ by $t/s$ and take the limit 
for $s\downarrow 0$, then we obtain in view of the 
closedness of $D$ that $(t\, v,\, q)\in D$, which implies 
that $t\in I^v_{0,\, q}$. Because 
this holds for every $t\in\R$, it follows that $I^v_{0,\, q}=\R$. 

If $v,\, w\in N$ then we have for every $r,\, s\in\R$ 
that $r\, v,\, s\, w\in N$. Moreover, we have for every 
$p\in Q$ that $(r\, v,\, p)\in D$, $(s\, w,\, p+r\, v)\in D$, 
hence $(r\, v+s\, w,\, p)\in D$ in view of iii) in 
Lemma \ref{actionlem}. If for any $t\in\R$ we 
replace $r$ and $s$ by $t\, r$ and $t\, s$, respectively, we obtain that 
$I^{r\, v+s\, w}_{0,\, p}=\R$, hence $r\, v+s\, w\in N$, 
and we conclude that $N$ is a linear subspace of $V$. 
Furthermore b) implies that $(v,\, p)\in D$ for every $v\in N$ and 
$p\in Q$, and it follows from Lemma \ref{actionlem} 
that the mapping $(v,\, p)\mapsto p+v$ defines an action 
of $N$ on $Q$. 
\end{proof}

\begin{lemma}
Let $Q$ be a complete and connected locally convex $V$\--parallel space. 
For any $v\in V$, the following conditions are equivalent.
\begin{itemize}
\item[i)] There exists $p\in Q$ such that 
$(v,\, p)\in D$ and $p+v=p$. 
\item[ii)] For all $p\in Q$ we have $(v,\, p)\in D$ 
and $p+v= p$. 
\end{itemize}

Let $P$ denote the set of all $v\in V$ such that 
i) or ii) holds. Then $P$ is a discrete additive subgroup 
of the linear subspace $N$ of $V$ which is defined in 
Lemma \ref{Nlem}. 
\label{perlem}
\end{lemma}
\begin{proof}
Assume that i) holds. Lemma \ref{surjlem} 
implies that for any $q\in Q$ there exists 
$u\in V$ such 
that $q=p+u$, and therefore 
\[
q+v=(p+u)+v=p+(u+v)=p+(v+u)=(p+v)+u=p+u=q, 
\]
which proves ii). 

If $p+v=p$, then it follows by induction on $k$ 
that $(k\, v,\, p)\in D$ and $p+k\, v=p$ for every 
positive integer $k$, and using i), ii) in 
Lemma \ref{actionlem} we obtain the same conclusions for 
all $k\in\Z$. This implies that $\Z\subset I^v_{0,\, p}$, 
which in turn implies that $I_{0,\, p}^v=\R$, 
because $I^v_{0,\, p}$ is an interval in $\R$. 
In view of Lemma \ref{Nlem}, we conclude that $v\in N$. 
 
It follows that $P$ is equal to the set of $v\in N$ 
such that the action of $v$ on $Q$ is trivial, and therefore 
$P$ is an additive subgroup of $N$. 
The last statement in Lemma \ref{surjlem} implies that 
$P$ is a discrete subset of $D_p$, hence of the closed subset $N$ 
of $D_p$. 
\end{proof}

\begin{theorem}
Let $Q$ be a complete and connected locally convex $V$\--parallel space. 
Let $P$ be the period group in $V$ defined in Lemma \ref{perlem}. 
Denote by $\R\, P$ the $\R$\--linear span of $P$ in $V$, 
where we note that $(\R\, P)/P$ is a torus. 
Let $C$ be a linear complement in $V$ of the linear span 
$\R\, P$ of $P$ in $V$. Let $p\in Q$. 
Then there is a convex closed subset $\Delta _p$ of  
$C$ such that $D_p=\Delta _p+ \R\, P$. 
Furthermore the mapping $\Phi_p :(v,\ w)\mapsto p+(v+w)$ 
defines an isomorphism of $V$\--parallel spaces 
from $\Delta _p\times ((\R\, P)/P)$ onto $Q$. 
The projection from 
$Q$ onto the $\R\, P$\--orbit space $Q/(\R\, P)$ induces an isomorphism 
from $\Delta _p$, viewed as a $C$\--parallel space, 
onto the $V/(\R\, P)$\--parallel space $Q/(\R\, P)$. 

There is a collection of linear forms $v^*_i$ on $V$ 
and real numbers $c_i$, 
where $i$ runs over some index set $I$, such that 
$D_p$ is equal to the set of all $v\in V$ such that 
$v^*_i(v)\geq c_i$ for all $i\in I$. For every 
such collection $\lambda _i$, $c_i$, 
the linear subspace $N$ of $V$ defined in Lemma \ref{Nlem} 
is equal to the common kernel of the linear forms 
$v^*_i$, $i\in I$, on $V$. 

$Q$ is compact if and only if $\Delta _p$ is compact, 
which implies that $N=\R\, P$, and $P$ is a cocompact discrete subgroup 
of the additive group $N$. 
If $Q$ is a compact connected locally convex polyhedral $V$\--parallel space, 
then $\Delta _p\simeq Q/(\R\, P)$ is a convex polytope in 
$C\simeq V/(\R\, P)$. 
\label{cornerclassthm}
\end{theorem}
\begin{proof}
Let $(v,\, p)\in D$, $(v',\, p)\in D$, and $p+v=p+v'$. 
It follows from Lemma \ref{actionlem} that $(-v' ,\, p+v') 
=(-v' ,\, p+v)\in D$, $(v-v',\, p)=(v+(-v'),\, p)\in D$ 
and $p+(v-v' )=(p+v)+(-v')=(p+v')+(-v')=p$, 
which in view of Lemma \ref{perlem} implies that 
$v-v'\in P\subset N$. 

Define $\Delta _p:=D_p\cap C$. If $v,\, w\in\Delta _p$, then it follows 
from the transitivity in Lemma \ref{surjlem} that 
there exists a $z\in V$ such that $(z,\, p+v)\in D$ 
and $p+w=(p+v)+z=p+(v+z)$, which implies that 
$u:=w-v-z\in P$. Replacing $z$ by $z+u$ we 
therefore can arrange that $z=w-v$. 
Now $(w-v,\, p+v)\in D$ implies that for every 
$t\in [0,\, 1]$ we have $(t\, (w-v),\, p+v)\in D$, which in 
combination with iii) in Lemma \ref{actionlem} implies 
that $(v+t\, (w-v),\,p)\in D$, and therefore 
$v+t\, (w-v)\in D_p\cap C=\Delta _p$. This proves that 
$\Delta _p$ is convex. It is a closed subset 
of $C$ because $D_p$ is a closed subset of $V$, see 
Lemma \ref{surjlem}. 

The fact that the additive group $N$ 
acts on $Q$, cf. Lemma \ref{perlem}, 
and hence its subgroup $\R\, P$ acts on $Q$, implies that 
$D_p=\Delta _p+\R\, P$. 

Because $v,\, v'\in D_p$ and $p+v=p+v'$ imply that 
$v-v'\in P$, and $C$ is complementary to $\R\, P$, the mapping 
$\Phi _p:\Delta _p\times ((\R\, P)/P)\to Q$ 
is injective. On the other hand Lemma \ref{surjlem} implies 
that $\Phi_p$ is a surjective local homeomorphism. 
Because $\Phi _p$ is an $L$\--map with 
$L:C\times (\R\, P)\to V:(w,\, z)\mapsto w+z$, 
it follows that $\Phi _p$ is an isomorphism 
from $\Delta _p\times ((\R\, P)/P)$ onto $Q$. 

The statement about $N$, the $v^*_i$ and the $c_i$ 
follows because $D_p=\Delta _p+(\R\, P)$ is a closed convex 
subset of $V$, and 
$N$ is equal to the lineality of $D_p$, the set of direction 
vectors of lines which are contained in $D_p$, 
cf. Rockafellar \cite[p. 65]{rockafellar}. 

Finally, a compact Hausdorff topological space has a unique 
uniform structure which is compatible with its topology, and 
the compact space is complete with respect to this 
uniform structure, cf. Bourbaki 
\cite[Ch. 2, \S 4, No. 1, Th. 1]{bourbaki}. 
Clearly $Q$ is compact if and only if 
$\Delta _p\times ((\R\, P)/P)$ is compact if and only if 
$\Delta _p$ is compact. In view of $N\cap C\subset\Delta _p$, the 
latter implies that $N\cap C=0$, hence 
$N=\R\, P$, and $P$ is a cocompact discrete subgroup 
of the additive topological group $N$. A convex compact subset 
of a vector space 
$C$ which is a $C$\--parallel manifold with corners is 
a convex polytope in $C$. 
\end{proof}

\begin{remark}
Theorem \ref{cornerclassthm} is a generalization of 
the theorem of Tietze \cite{tietze} and Nakajima \cite{nakajima} 
that any closed, connected, and locally 
convex subset of a finite\--dimensional vector space 
is convex. Our proof of Lemma \ref{actionlem} is close to the 
proof of Klee \cite[(5.2)]{klee} of the generalization of the 
Tietze\--Nakajima theorem to subsets of 
arbitrary topological vector spaces. 
\label{tnrem}
\end{remark}

\section{The symplectic tube theorem}
\label{sympltubes}
In this section we describe the local model 
of Benoist \cite[Prop. 1.9]{benoist} and 
Ortega and Ratiu \cite{ortegaratiu} for a general proper 
symplectic Lie group action. See also 
Ortega and Ratiu \cite[Sec. 7.2--7.4]{orbook} 
for a detailed proof. 
For Hamiltonian actions, such local models 
had been obtained before by Marle \cite{marle} 
and Guillemin and Sternberg \cite[Sec. 41]{gsst}. 

\medskip
Let $(M,\,\sigma )$ be a smooth symplectic manifold 
and $G$ a Lie group which acts smoothly on $(M,\,\sigma )$ 
by means of symplectomorphisms. Furthermore assume that the 
action is {\em proper}, which means that for any compact 
subset $K$ of $M$ the set of all $(g,\, m)\in G\times M$ 
such that $(m,\, g\cdot m)\in K$ is compact in $G\times M$. 
The action of $G$ is certainly proper if $G$ is 
compact. 

For every $g\in G$ we will write $g_M:x\mapsto g\cdot x$ 
for the action of $g$ on $M$. For every element $X$ in 
the Lie algebra $\got{g}$ of $G$, the infinitesimal action 
on $M$ will be denoted by $X_M$. It is a smooth vector field $X_M$ on $M$, 
the flow of which leaves $\sigma$ invariant.  

It follows from the properness of the action, 
that for every $x\in M$ the stabilizer subgroup 
$H:=G_x:=\{ g\in G\mid g\cdot x =x\}$ of $x$ in $G$ is a compact, 
Lie subgroup of $G$, and the mapping $A_x:g\mapsto g\cdot x:G\to M$ 
induces a $G$\--equivariant 
smooth embedding 
\begin{equation}
\alpha _x:g\, H\mapsto g\cdot x:G/H\to M 
\label{Axdef}
\end{equation} 
from $G/H$ into $M$, with closed image, 
equal to the orbit $G\cdot x$ of $G$ through the point $x$. 
Here $g\in G$ acts on $G/H$ by sending $g'\, H$ to 
$(g\, g')\, H$. 
The Lie algebra $\got{h}:=\got{g}_x$ of $H:=G_x$ is equal to the set 
of $X\in\got{g}$ such that $X_M(x)=0$. 
The linear mapping 
$\op{T}_1\! A_x:\got{g}\to\op{T}_x\! M$ induces induces a linear 
isomorphism from $\got{g}/\got{h}=\got{g}/\got{g}_x$ onto 
$\got{g}_M(x):=\op{T}_x(G\cdot x)$.

For the description of the symplectic form in the local model, 
we begin with the closed two\--form 
\begin{equation}
\sigma ^{G/H}:=(\alpha _x)^*\sigma 
\label{sigmaGdef}
\end{equation} 
on $G/H$, which represents the ``restriction'' of $\sigma$ to the 
orbit $G\cdot x\simeq G/H$ through the point $x$. 
Here $\alpha _x:G/H\to M$ is defined in (\ref{Axdef}). 
The $G$\--invariance 
of $\sigma$ and the $G$\--equivariance of $\alpha _x$ 
imply that $\sigma ^{G/H}$ is 
$G$\--invariant. 

If we identify $\op{T}_{1\, H}(G/H)$ 
with $\got{g}/\got{h}$, and $p:\got{g}\to\got{g}/\got{h}$ 
denotes the projection $X\mapsto X+\got{h}$, 
then $\sigma ^{G/H}_{1\, H}$, and therefore the 
$G$\--invariant two\--form $\sigma ^{G/H}$ on $G/H$, 
is determined by the antisymmetric bilinear form 
\begin{equation}
\sigma ^{\got{g}}:=(\op{T}_1\! A_x)^*\,\sigma _x=p^*\sigma ^{G/H}_{1\, H}
\label{sigmagdef}
\end{equation}
on $\got{g}$, which 
is invariant under the 
adjoint action of $H$ on $\got{g}$. It follows that   
the kernel 
\begin{equation}
\got{l}:=\op{ker}\sigma ^{\got{g}}=\op{ker}((\op{T}_1\! A_x)^*\,\sigma _x)
\label{ldef}
\end{equation}
is an $\op{Ad}H$\--invariant linear subspace of $\got{g}$. 
We have $\got{h}\subset\got{l}$, because of 
(\ref{sigmagdef}) and the fact that $p^*\sigma ^{G/H}_{1\, H}$ 
vanishes on the kernel $\got{h}$ of $p$. 

If $L$ is a linear subspace of a symplectic vector space 
$(V,\,\sigma )$, then the symplectic orthogonal complement 
$L^{\sigma }$ of $L$ in $V$ is defined as the set of 
all $v\in V$ such that $\sigma (l,\, v)=0$ for every $l\in L$. 
The restriction to $\got{g}_M(x)^{\sigma _x}$ of 
$\sigma _x$ defines a symplectic form $\sigma ^W$ on the 
vector space 
\begin{equation}
W:=\got{g}_M(x)^{\sigma _x}/(\got{g}_M(x)^{\sigma _x}\cap\got{g}_M(x)),
\label{Wdef}
\end{equation}
and the mapping 
$X+\got{h}\mapsto X_M(x)$ defines a linear isomorphism from 
$\got{l}/\got{h}$ onto $\got{g}_M(x)^{\sigma _x}\cap\got{g}_M(x)$. 
The linearized action $H\ni h\mapsto\op{T}_x\! h_M$ 
of $H$ on $\op{T}_x\! M$ is symplectic and leaves 
$\got{g}_M(x)\simeq \got{g}/\got{h}$ invariant, acting on it 
via the adjoint representation. That is, 
$h\cdot X_M(x)=(\op{Ad}h)(X)_M(x)$, if 
$X\in\got{g}$ and $h\in H$. It therefore also 
leaves $\got{g}_M(x)^{\sigma _x}$ invariant and induces 
an action of $H=G_x$ on the symplectic vector space 
$(W,\,\sigma ^W)$ by means of symplectic linear transformations. 

With $\got{l}$ as in (\ref{ldef}), we ``enlarge''  
the vector space $W$ to the vector space 
\begin{equation}
E:=(\got{l}/\got{h})^*\times W,
\label{Edef}
\end{equation}
on which $h\in H$ acts by sending $(\lambda ,\, w)$ 
to $(((\op{Ad}h)^*)^{-1}(\lambda ),\, h\cdot w)$. 

For any action by linear transformations of a compact Lie group 
$K$ on a vector space $V$, any $K$\--invariant linear 
subspace $L$ of $V$ has an $K$\--invariant linear complement 
$L'$ in $V$. For instance, if $\beta$ is an inner product 
on $V$, then the average $\overline{\beta}$ of $\beta$ 
over $K$ is a $K$\--invariant inner product on $V$, 
and the $\overline{\beta}$\--orthogonal complement 
$L'$ of $L$ in $V$ has the desired properties. 
Choose $\op{Ad}H$\--invariant linear complements 
$\got{k}$ and $\got{c}$ of $\got{h}$ and $\got{l}$ in 
$\got{g}$, respectively.  
Let $X\mapsto X_{\got{l}}:\got{g}\to\got{l}$ and 
$X\mapsto X_{\got{h}} :\got{g}\to\got{h}$ denote the 
linear projection from $\got{g}$ onto $\got{l}$ and $\got{h}$ 
with kernel equal to $\got{c}$ and $\got{k}$, respectively.  
Then these linear projections are $\op{Ad}H$\--equivariant. 

If $g\in G$, then we denote by $\op{L}_g:g'\mapsto g\, g':G\to G$ 
the multiplication from the left by means of $g$. 
Define the smooth one\--form $\eta ^{\#}$ 
on $G\times E$ by 
\begin{equation}
\eta ^{\#}_{(g,\, (\lambda ,\, w))}((\op{T}_1\op{L}_g)(X),\, 
(\delta\lambda ,\,\delta w))
:=\lambda (X_{\got{l}})+\sigma ^W(w,\,\delta w+X_{\got{h}}\cdot w)/2
\label{etadef}
\end{equation}
for all $g\in G$, $\lambda\in (\got{l}/\got{h})^*$, $w\in W$, 
and all $X\in\got{g}$, $\delta\lambda\in 
(\got{l}/\got{h})^*$,  
$\delta w\in W$. Here we identify the tangent 
spaces of a vector space with the vector space itself 
and $(\got{l}/\got{h})^*$ with the space of linear forms on 
$\got{l}$ which vanish on $\got{h}$.  

Let $E$ be defined as in (\ref{Edef}), with 
$\got{l}$ and $W$ as in (\ref{ldef}) and (\ref{Wdef}), 
respectively. 
Let $G\times _HE$ denote the orbit space of $G\times E$ 
for the proper and free action of $H$ on $G\times E$, 
where $h\in H$ acts on $G\times E$ by sending 
$(g,\, e)$ to $(g\, h^{-1},\, h\cdot e)$. 
The action of $G$ on $G\times _HE$ is induced 
by the action $(g,\, (g', e))\mapsto (g\, g',\, e)$ 
of $G$ on $G\times E$. Let $\pi :G\times _HE\to G/H$ 
denote the mapping which is induced by the projection 
$(g,\, e)\mapsto g:G\times E\to G$ 
onto the first component. Because $H$ acts on $E$ by means 
of linear transformations, this projection 
exhibits $G\times _HE$ as a $G$\--homogeneous vector bundle 
over the homogeneous space $G/H$, which fiber $E$ and 
structure group $H$. 
With these notations, we have the following 
local normal form for the symplectic $G$\--space  
$(M,\,\sigma ,\, G)$. 
\begin{theorem}
If $\pi _H:G\times E\to G\times _HE$ denotes $H$\--orbit mapping, 
then there is a unique smooth one\--form $\eta$ on 
$G\times_HE$, such that $\eta ^{\#}={\pi _H}^*\,\eta$. 
Here $\eta ^{\#}$ is defined in {\em (\ref{etadef})}. 

Furthermore, there exists an open $H$\--invariant 
neighborhood $E_0$ of the origin in $E$ 
and a $G$\--equivariant diffeomorphism $\Phi$ from 
$G\times _HE_0$ onto an open $G$\--invariant neighborhood 
$U$ of $x$ in $M$, such that $\Phi (H\cdot (1,\, 0))=x$ and 
\begin{equation}
\Phi ^*\sigma =\pi ^*\sigma ^{G/H}+\op{d}\!\eta .
\label{Gsigma}
\end{equation}
Here the mapping $\pi :G\times _HE_0\to G/H$ is 
induced by the projection onto the first component, 
and $\sigma ^{G/H}$ is defined in 
{\em (\ref{sigmaGdef})}. 
\label{Gthm}
\end{theorem}

\noindent
J.J. Duistermaat\\
Mathematisch Instituut, Universiteit Utrecht\\
P.O. Box 80 010, 3508 TA Utrecht, The Netherlands\\
e\--mail: duis@math.uu.nl

\bigskip\noindent
A. Pelayo\\
Department of Mathematics, University of Michigan\\
2074 East Hall, 530 Church Street, Ann Arbor, MI 48109--1043, USA\\
e\--mail: apelayo@umich.edu


\begin{thebibliography}{99}

\bibitem{atiyah}M. Atiyah: Convexity and commuting Hamiltonians. 
{\em Bull. London Math. Soc.} {\bf 14} (1982) 1--15. 

\bibitem{audin}M. Audin: {\em Torus Actions on Symplectic Manifolds}. 
Second revised edition. Progress in Mathematics, 93. 
Birkh\"{a}user Verlag, Basel, 2004. 

\bibitem{am}L. Auslander and L. Markus: 
Holonomy of flat affinely connected manifolds. 
{\em Ann. of Math.} {\bf 62} (1955) 139--151. 

\bibitem{auslander}L. Auslander: The structure of complete 
locally affine manifolds. {\em Topology} {\bf 3} (1964) 131--139. 

\bibitem{bg}C. Benson and C.S. Gordon: K\"ahler and symplectic structures 
on nilmanifolds. {\em Topology} {\bf 27} (1988) 513--518. 

\bibitem{benoist}Y. Benoist: Actions symplectiques de groupes compacts. 
{\em Geometriae Dedicata} {\bf 89} (2002) 181--245.  

\bibitem{benoistcorr}Y. Benoist: 
Correction to ``Actions symplectiques de groupes compacts''.\\ 
http://www.dma.ens.fr/~benoist.  

\bibitem{bt}R. Bott and L.W. Tu: {\em Differential Forms 
in Algebraic Topology}. Springer\--Verlag, New York, Heidelberg, 
Berlin, 1982.  

\bibitem{bourbaki}N. Bourbaki: {\em \'El\'ements de Math\'ematique, 
Premi\`ere Partie, Livre III: Topologie G\'en\'erale, 
Chap. 1-2}. Hermann, Paris, 1961.  

\bibitem{cartan}\'E. Cartan: Sur les nombres de Betti des 
espaces de groupes clos. {\em C.R. Acad. Sc.} {\bf 187} (1928) 
196--198 = {\em {\OE}uvres}, partie I, vol. 2, 999--1001. 

\bibitem{delzant}T. Delzant: Hamiltoniens p\'eriodiques et image convex 
de l'application moment. {\em Bull. Soc. Math. France} {\bf 116} (1988) 
315--339. 

\bibitem{duis}J.J. Duistermaat: Equivariant cohomology 
and stationary phase. {\em Contemp. Math.} {\bf 179} 
(1994) 45--62.  

\bibitem{heatkernel}J.J. Duistermaat; 
{\em The Heat Kernel Lefschetz Fixed Point Formula 
for the Spin\--c Dirac Operator}. Birkh\"auser Boston, 1996. 

\bibitem{dh}J.J. Duistermaat and G.J. Heckman: 
On the variation in the cohomology of the symplectic form 
of the reduced phase space. {\em Invent. Math.} {\bf 69} (1982) 259--268. 
Addendum in: {\em  Invent. Math.} {\bf 72} (1983) 153--158. 

\bibitem{dk}J.J. Duistermaat and J.A.C. Kolk: {\em Lie Groups}. 
Springer\--Verlag, Berlin, Heidelberg, 2000. 

\bibitem{danilov}V.I. Danilov: The geometry of toric varieties. 
{\em Russ. Math. Surveys} {\bf 33}:2 (1978) 97--154 = from 
{\em Uspekhi Mat. Nauk SSSR} {\bf 33}:2 (1978) 85--134. 

\bibitem{giacobbe}A. Giacobbe: Convexity of multi\--valued 
momentum maps. {\em Geometriae Dedicata} {\bf 111} (2005) 1--22.  

\bibitem{ginzburg}V.L. Ginzburg: 
Some remarks on symplectic actions of compact groups.
{\em Math. Z.} {\bf 210} (1992) 625--640.

\bibitem{greenberg}M. Greenberg: {\em Lectures on 
Algebraic Topology}. W.A. Benjamin, New York, Amsterdam, 1967. 

\bibitem{gls}V. Guillemin, E. Lerman and S. Sternberg; 
{\em Symplectic fibrations and multiplicity diagrams.} 
Cambridge University Press, Cambridge, 1996.

\bibitem{gs}V. Guillemin and S. Sternberg: Convexity properties 
of the moment mapping. {\em Invent. Math.} {\bf 67} (1982) 491--513. 

\bibitem{multfree}V. Guillemin and S. Sternberg: Multiplicity\--free 
spaces. {\em J. Diff. Geom.} {\bf 19} (1984) 31--56. 

\bibitem{gsst}V. Guillemin and S. Sternberg: {\em Symplectic 
Techniques in Physics}. Cambridge University Press, 
Cambridge, etc., 1984. 

\bibitem{guillemin}V. Guillemin: {\em Moment Maps and 
Combinatorial Invariants of Hamiltonian $T^n$\--spaces.} 
Birkh\"auser, Boston, Basel, Berlin, 1994.  

\bibitem{hs}A. Haefliger and \'E. Salem: 
Actions of tori on orbifolds.
{\em Ann. Global Anal. Geom.} {\bf 9} (1991) 37--59.

\bibitem{hungerford}T.W. Hungerford: {\em Algebra}. 
Springer\--Verlag, New York, 1974. 

\bibitem{K} Y. Karshon:  
{\em Periodic Hamiltonian flows on four-dimensional manifolds. }
Memoirs Amer. Math. Soc. No. 672  {\bf 141} (1999), viii+71 pp. 

\bibitem{KT} Y. Karshon and S. Tolman: Centered complexity one 
Hamiltonian torus actions, 
{\em Trans. Amer. Math. Soc.} {\bf 353} (2001), no. 12, 4831--4861.  

\bibitem{klee}V. Klee: Convex sets in linear spaces. 
{\em Duke Math. J.} {\bf 18} (1951) 443--466. 

\bibitem{Ko} M. Kogan, On completely integrable systems with 
local torus actions. 
{\em Ann. Global Anal. Geom.} {\bf 15} (1997), no. 6, 543--553. 

\bibitem{koszul}J.L. Koszul: Sur certains groupes de transformations de Lie. 
pp. 137-141 in: {\em G\'eom\'etrie Diff\'erentielle.} 
Coll. Int. du C.N.R.S., Strasbourg, 1953. 

\bibitem{ls}N.C. Leung and M. Symington: Almost toric symplectic 
four-manifolds. \\arXiv:math.SG/0312165v1, 8 Dec 2003. 

\bibitem{maclane}S. MacLane: {\em Categories for the working mathematician.} 
Graduate Texts in Mathematics, 5. Springer-Verlag, New York, 1971, 1998. 

\bibitem{malcev}A.I. Malcev: On a class of homogeneous spaces. 
{\em Izv. Akad. Nauk SSSR Ser. Mat.} {\bf 13}, 1 (1949) 9-32 
(MR 0028842 (10, 507 d)) = {\em Amer. Math. Soc. Translations} 
(1951) No. 39, 33 pp. 

\bibitem{marle}C.-M. Marle: Classification des actions 
hamiltoniennes au voisinage d'une orbite. 
{\em C. R. Acad. Sci. Paris S\'er. I Math.}  
{\bf 299} (1984) 249--252. 
Mod\`ele d'action hamiltonienne d'un groupe de Lie sur une vari\'et\'e 
symplectique. {\em Rend. Sem. Mat. Univ. Politec. Torino} 
{\bf 43} (1985) 227--251. 

\bibitem{mather}J.N. Mather: Stability of $C^{\infty}$ 
mappings: II. Infinitesimal stability implies stability. 
{\em Ann. of Math.} {\bf 89} (1969) 254--291.

\bibitem{M} D. McDuff: The moment map for circle actions on 
symplectic manifolds. {\em J. Geom. Phys.} 
{\bf 5} (1988), no. 2, 149--160. 

\bibitem{MS} D. McDuff and D.A. Salamon:   
{\em Introduction to Symplectic Topology}, 2nd edition. OUP, 1998. 

\bibitem{nakajima}S. Nakajima: \"{U}ber konvexe Kurven und Fl\"{a}chen. 
{\em T\^{o}hoku Math. J.} {\bf 29} (1928) 227--230.

\bibitem{novikov}S.P. Novikov: The Hamiltonian formalism and a 
multivalued analogue of Morse theory. 
{\em Russ. Math. Surveys} {\bf 37}:5 (1982) 1--56. 

\bibitem{OR}P. Orlik and  F. Raymond: Actions of the torus on 4-manifolds, I  
{\em Trans. Amer. Math. Soc.} {\bf 152} (1970), 531-559, 
II {\em Topology} {\bf 13} (1974), 89-112. 

\bibitem{ortegaratiu}J.-P Ortega and T.S. Ratiu: A symplectic slice theorem.  
Lett. Math. Phys. {\bf 59} (2002) 81--93. 

\bibitem{orbook}J.-P Ortega and T.S.Ratiu: 
{\em Momentum maps and Hamiltonian reduction.}
Progress in Mathematics, 222.
Birkh\"auser Boston, Boston, MA, 2004.
 
\bibitem{ps}R.S. Palais and T.E. Stewart: Torus bundles over a torus. 
{\em Proc. Amer. Math. Soc.} {\bf 12} (1961) 26--29.  

\bibitem{pao}P. S. Pao: The topological structure of 4-manifolds with 
effective torus actions, I 
{\em Trans. Amer. Math. Soc.} {\bf 227} (1977), 279-317, 
II {\em Ill. J.Math.} {\bf  21} (1977), 883-894. 

\bibitem{rockafellar}R.T. Rockafellar: {\em Convex Analysis}. 
Princeton University Press, Princeton, New Jersey, 1970. 

\bibitem{steenrod}N. Steenrod: {\em The Topology of Fibre Bundles}. 
Princeton University Press, Princeton, New Jersey, 1951, 1972. 

\bibitem{s}M. Symington: Four dimensions from two in symplectic topology.
pp. 153--208 in {\em Topology and geometry of manifolds (Athens, GA, 2001)}. 
Proc. Sympos. Pure Math., 71, Amer. Math. Soc., Providence, RI, 2003. 

\bibitem{thurston}W.P. Thurston: Some examples of symplectic 
manifolds. {\em Proc. Amer. Math. Soc.} {\bf 55} (1976) 467--468. 

\bibitem{tietze}H. Tietze: \"{U}ber Konvexheit im kleinen und 
im gro{\ss}en und \"{u}ber gewisse den Punkten einer Menge 
zugeordete Dimensionszahlen. {\em Math. Z.} {\bf 28} (1928) 697--707. 
\bibitem{whitney}H. Whitney: Differentiable even functions. 
{\em Duke Math. J.} {\bf 10} (1943). 159--160. 

\end{thebibliography}
\end{document}